\documentclass[10pt,reqno]{amsart}
\usepackage{amsmath,amssymb,latexsym,soul,cite,mathrsfs}
\usepackage{subfigure}

\usepackage{color,enumitem,graphicx}
\usepackage[colorlinks=true,urlcolor=blue,
citecolor=red,linkcolor=blue,linktocpage,pdfpagelabels,
bookmarksnumbered,bookmarksopen]{hyperref}
\usepackage[english]{babel}
\usepackage{enumitem}

\usepackage[left=2.7cm,right=2.7cm,top=2.9cm,bottom=2.9cm]{geometry}

\usepackage[hyperpageref]{backref}

\makeatletter
\newcommand{\leqnomode}{\tagsleft@true}
\newcommand{\reqnomode}{\tagsleft@false}
\makeatother

\numberwithin{equation}{section}

\newtheorem{theorem}{Theorem}[section]
\newtheorem{lemma}[theorem]{Lemma}
\newtheorem{corollary}[theorem]{Corollary}
\newtheorem{proposition}[theorem]{Proposition}

\newtheorem{definition}[theorem]{Definition}

\newtheorem{example}[theorem]{Example}

\renewcommand{\rightarrow}{\to}

\title[Weingarten Surfaces Associated to Laguerre Minimal Surfaces ]{Weingarten Surfaces Associated to Laguerre Minimal Surfaces}

\author[Laredo Rennan Pereira Santos]{Laredo Rennan Pereira Santos}
\author[Armando Mauro Vasquez Corro]{Armando V. Corro}

\address[Laredo Rennan Pereira Santos]{Instituto Federal de Educa\c{c}\~{a}o, Ci\^{e}ncia e Tecnologia de Goi\'{a}s
	\newline\indent
Rua 64, Expans\~{a}o Parque Lago, Formosa-Brasil.
	\newline\indent
	73813816, Goi\'{a}s-GO, Brasil}
\email{\href{mailto:laredo.santos@ifg.edu.br}{laredo.santos@ifg.edu.br}}

\address[Armando V. Corro]{Instituto de Matem\'{a}tica e Estat\'{i}stica,
	Universidade Federal de Goi\'{a}s
	\newline\indent
	74001-970, Goi\'{a}s-GO, Brasil}
\email{\href{mailto:corro@mat.ufg.br}{corro@mat.ufg.br}}

\thanks{{\textbf{Corresponding author}: Laredo Rennan Pereira Santos\{laredo.santos@ifg.edu.br\}: Instituto Federal de Educa\c{c}\~{a}o, Ci\^{e}ncia e Tecnologia de Goi\'{a}s, 
	Formosa-Goi\'{a}s, Brasil. ORCID number: 0000-0001-5216-2026}}

\thanks{{\textbf{Armando V. Corro} \{corro@mat.ufg.br\}: Instituto de Matem\'{a}tica e Estat\'{i}stica,
	Universidade Federal de Goi\'{a}s, Goi\^ania-Goi\'{a}s-GO, Brasil.    }}

\thanks{  }


\begin{document}
	

\begin{abstract}
In the work \cite{Laredo} the author shows that every hypersurface in Euclidean space is locally associated to the unit sphere by a sphere congruence, whose radius function $R$ is a geometric invariant of hypersurface. In this paper we define for any surface $\Sigma$ its spherical mean curvature $H_S$ which depends on principal curvatures of $\Sigma$ and the radius function $R$. Then we consider two classes of surfaces: the ones with $H_S = 0$, called $H_1$-surfaces, and the surfaces with spherical mean curvature of harmonic type, named $H_2$-surfaces. We provide for each these classes a Weierstrass-type representation depending on three holomorphic functions and we prove that the $H_1$-surfaces are associated to the minimal surfaces, whereas the $H_2$-surfaces are related to the Laguerre minimal surfaces. As application we provide a new Weierstrass-type representation for the Laguerre minimal surfaces - and in particular for the minimal surfaces - in such a way that the same holomorphic data provide examples in  $H_1$-surface/minimal surface classes or in $H_2$-surface/Laguerre minimal surface classes.  We also characterize the rotational cases, what allow us finding a complete rotational Laguerre minimal surface.

 \vspace{0.5cm}
\begin{center}
 \textbf{keywords}: Generalized Weingarten Surfaces, Laguerre Minimal Surfaces, Weierstrass-Type Representation
\end{center} 
 
 
\end{abstract}

\maketitle

	

\section{Introduction}

For long time the plane, the helicoid and the catenoid were the only examples known from minimal surfaces. The next example was given by Scherk in 1835, whose surface bears his name, and in 1864 Enneper exhibits the simplest example of minimal surface found until that moment. The classical theory of minimal surfaces experiences a great advance in its so-called first golden age, from 1855 to 1890, approximately, and this thanks to the connection between the theory of minimal surfaces and the Analysis Complex. The Weierstrass representation for minimal surfaces is a big sign of this phenomenon, allowing to obtain examples of such surfaces from a pair of holomorphic functions. This representation was obtained locally by Weierstrass in 1866 and its key point is to provide a recipe for defining a multitude of examples of minimal surfaces. One of his versions can be seen in \cite{Carmo}.

In the literature there exist Weierstrass-type representations for some classes of surfaces, among which are certain classes of Weingarten surfaces, object of interest of this work.

An oriented surface $S \subset \mathbb{R}^3$ is said to be a Weingarten surface if there exists a differentiable relation $W$ between its Gaussian curvature $K$ and its mean curvature $H$ such that $W(H,K) \equiv 0 $. They were introduced by Weingarten (\hspace{-0.01cm}\cite{W 1}, \cite{W 2}) in an attempt to find a class of isometric surfaces to a given surface of revolution. Surfaces of constant Gaussian curvature and surfaces of constant mean curvature (in particular, minimal surfaces) are examples of Weingarten surfaces.

The general classification of Weingarten surfaces is still an open question. Classification of certain classes of Weingarten surfaces has been reported in the literature and a great number of them arise in many situations. Papantoniou \cite{P} classified the Weingarten surfaces of revolution whose principal curvatures $k_i$ satisfy a linear relation  $Ak_1 + Bk_2 = 0$, with $A$ and $B$ not simultaneously null. Schief \cite{Schief} studied \textit{generalized Weingarten surfaces}, which accept the relation $1 + \mu H + (\mu^2 \pm \rho^2)K = 0$
where the functions $\rho$ and $\mu$ are harmonics in a certain sense.  These surfaces have been shown to be integrable. 

In \cite{Dias 1}, Dias introduced a class of oriented Weingarten surfaces in $\mathbb{R}^3$ that satisfies a relation of the form 
$$
A(\Psi_\nu, \Lambda_\nu) + B(\Psi_\nu, \Lambda_\nu)H + C(\Psi_\nu, \Lambda_\nu)K = 0,
$$
where $A, B, C: \mathbb{R}^2 \rightarrow \mathbb{R}$ are differentiable functions that depend on  the support function $\Psi_\nu$ and the quadratic distance function $\Lambda_\nu$, for some fixed point $\nu \in \mathbb{R}^3$. These surfaces are referred to as   \textit{generalized Weingarten surfaces that depend on the support and distance functions} (in short, DSGW surfaces). Some of the main classes of Weingarten surfaces studied in the literature are classes of DSGW surfaces, such as the linear Weingarten, Appell, and Tzitz\'{e}ica's surfaces.

In 1888, Appell \cite{Appell} studied a class of oriented surfaces in $\mathbb{R}^3$ associated with  area preserving transformation in the sphere. Later, Ferreira and Roitman \cite{F e R} found that these surfaces satisfy the Weingarten relation $H+\Psi_\nu K=0$, for a fixed point $\nu \in \mathbb{R}^3$. 

A typical method to characterize classes of Weingarten surfaces is to provide them with a Weierstrass-type representation by which they can be parameterized in terms of holomorphic functions. In this sense we can highlight the papers \cite{Dias 1}, \cite{Corro} and \cite{Caroline}.

Considering surfaces $\Sigma$ for which $\Psi(p) \neq 1 $, for all $p \in \Sigma$, we introduce its \textit{radial curvatures $s_i$ associated to sphere} and \textit{mean curvature $H_S$ associated to sphere} as follows:
\begin{equation*}
s_{i}= \frac{1+k_i}{1-k_iR}, \hspace{1.0cm} H_S = \frac{1}{n}\sum_{i=1}^{n}s_i,
\end{equation*}
where  $k_i$ are the principal curvatures of $\Sigma$, \hspace{0.1cm} $1\leqslant i \leqslant n$, and $R$ is a geometric invariant of $\Sigma$ given by the radius function of a sphere congruence. In this way, the surface $\Sigma$ is named  \textit{surface of null spherical mean curvature} (in short, $H_1$-surface) if $H_S \equiv 0$ and  $\Sigma$ is called a \textit{surface with spherical mean curvature of harmonic type} (in short, $H_2$-surface) if it holds

\begin{equation*}
\Delta_\sigma \left[\frac{H_S}{\Psi - 1} \right] = 0,
\end{equation*}
where $\sigma = I + 2R \hspace{0.1cm}II + R^2 \hspace{0.1cm}III$, with $I, \hspace{0.1cm} II, \hspace{0.1cm} III$  the fundamental forms of $\Sigma$.

Considering a function  $h: \mathbb{R}^2 \rightarrow \mathbb{R}$, the two-dimensional \textit{Helmholtz equation} for $h$ is defined by

$$
\Delta h(z) + c(K(z))^2 h(z) = 0,
$$
and the two-dimensional \textit{generalized Helmholtz equation} for $h$ is given by

$$
\Delta \left( \frac{1}{(K(z))^2}\big(\Delta h(z) + c(K(z))^2 h(z)\big)\right) = 0
$$
where $K(z)$ is a function and $c$ is a real non-zero constant. We show that the $H_1$-surfaces (resp. $H_2$-surfaces) are associated to the solutions of the two-dimensional Helmholtz equation (resp. two-dimensional generalized Helmholtz equation). In fact, we prove that both classes has a parameterization of the form

\begin{equation}\label{X}
X(u) = Y(u) -2 \left(\frac{h(u)+c}{S(u)}\right)\eta(u), \hspace{0.5 cm}u \in U\subset \mathbb{R}^2,
\end{equation}
\vspace{0.3cm}
where  $c$ is a nonzero real constant, $Y$ is an orthogonal parameterization of the unit sphere, $h$ is a solution for a two-dimensional generalized Helmholtz equation e 
\begin{equation*}
\eta = \nabla_{L}h +hY, \hspace{1.0cm}
S = \left< \eta, \eta\right> =  \big|\nabla_Lh\big|^2 + h^2,
\end{equation*}	
with $ L_{ij} = \left < Y_{,i}, Y_{,j} \right >$. We find that when the parameterization $X$ in \eqref{X} defines a $H_1$-surface (resp. $H_2$-surface), then $\eta$ is an immersion that defines a minimal surface (resp. Laguerre minimal surface). Thus, after taking a suitable parameterization $Y$ for $\mathbb{S}^n$, the immersion $\eta$ can be expressed in terms of three holomorphic functions and becomes an alternative Weierstrass representation for the Laguerre minimal (in particular, minimal) surfaces.

Laguerre minimal surfaces $ \Sigma$ in  $ \mathbb{R}^3$ are critical points of the functional
\begin{equation*}
\int_\Sigma\frac{(H^2-K)}{K} \hspace{0.1cm} dA,
\end{equation*}
where $H$ and $K$ denote the mean and Gaussian curvatures of $\Sigma$, respectively, and $dA$ is the $\Sigma$ area element. The Euler-Lagrange equation of the Laguerre minimal surfaces is given by 

\begin{equation*}
\Delta_{III}\left(\frac{H}{K}\right)=0,
\end{equation*}
where $\Delta_{III}$ is the Laplacian operator with respect to the third fundamental form $III$ of $\Sigma$. These surfaces have been extensively studied and as example we cite  \cite{Wesley 2} in which is presented a classification for the Laguerre minimal surfaces with planar curvature lines.

Finally we consider the rotational cases for the $H_1$ and $H_2$-surfaces, giving as application a characterization for the rotational  Laguerre minimal surfaces and exhibiting for them a complete example.

The  paper is organized as follows. Section $2$ is devoted to certain classical definitions and theorems in Differential Geometry and Complex Analysis and to the presentation of the results concerning Helmholtz equation and sphere congruence used in the text. In Section $3$, we define and discuss the $H_1$ and $H_2$-surfaces and we establish a link between them and the solutions for certain Helmholtz equations. We also show that for each $H_1$-surface (resp. $H_2$-surface) there exists a correspondent minimal surface (resp. Laguerre minimal surface), providing for them Weierstrass-type representations. In section $4$ we consider the rotational $H_1$-surfaces and section $5$ deals with the rotational cases for $H_2$-surfaces and  Laguerre minimal surfaces.

\section{Preliminaries}
In this section we give the notation and the main classical results in the literature that will be used in the work.

\subsection{Hypersurfaces in the Euclidean Space}
 Throughout this paper $f_{,i}$ indicates the partial derivative of a differentiable function $f: \mathbb{R}^n \rightarrow \mathbb{R}^m$ with respect to $i$-th variable, $U$ denotes an open subset of $\mathbb{R}^n$ and $\Sigma$ a hypersurface in $\mathbb{R}^{n+1}$ with   normal Gauss map $N$. In this sense, if $X:U \rightarrow \Sigma$ is a local parameterization of $\Sigma$, the matrix $W = (W_{ij})$ such that
$$
N_{,i}= \sum_{j=1}^{n}W_{ij}X_{,j}, \hspace{0.3cm} 1\leq i \leq n,
$$ 
is called the\textit{ Weingarten matrix} of $\Sigma$. The vector $X_{,ij}$, $ 1\leqslant i, j \leqslant n$, can be written as 
\begin{equation}\label{eq Xij}
X_{,ij}= \sum_{k=1}^{n}\Gamma_{ij}^{k}X_{,k}+b_{ij}N,
\end{equation}
and if the parameterization  $X$ is such that the metric  $g_{ij} = \left < X_{,i}, \hspace{0.1cm}X_{,j} \right> $ is diagonal, the Christoffel symbols satisfy
\begin{eqnarray}\label{eq Chris}
	\Gamma_{ij}^k &=& 0, \hspace{0.5cm} for \hspace{0.2cm} i,j,k \hspace{0.2cm} distincts, \nonumber\\
	\Gamma_{ij}^j &=& \frac{g_{jj,i}}{2g_{jj}}, \hspace{0.5cm}  for \hspace{0.1cm}all \hspace{0.2cm} i,j;\\
	\Gamma_{ii}^j &=& -\frac{g_{ii,j}}{2g_{jj}} = -\frac{g_{ii}}{g_{jj}}\Gamma_{ji}^{i}, \hspace{0.5cm} for \hspace{0.1cm} i\neq j. \nonumber
\end{eqnarray}

The first fundamental form $I$ of $\Sigma$ is the standard scalar product of $ \mathbb{R}^{n+1}$ restricted to the tangent hyperplanes   $T_p\Sigma$, whereas the second and third forms of $\Sigma$, denoted by $II$ and $III$ respectively, are defined as 
$$
II_p\big(w_1, w_2\big) = \left< -dN_p(w_1), w_2 \right>, \hspace{0.5cm} w_1, w_2 \in T_p\Sigma,
$$
$$
III_p\big(w_1, w_2\big) = \left< -dN_p(w_1), -dN_p(w_2) \right>, \hspace{0.5cm} w_1, w_2 \in T_p\Sigma,
$$
where $p \in \Sigma$ and  $dN_p$ is the differential of the normal Gauss map in $p$.

Take $\Sigma$ oriented by its normal Gauss map $N$. Given $  \nu  \in \mathbb{R}^n$, the functions $\Psi_\nu, \Lambda_\nu: \Sigma \rightarrow \mathbb{R}$ given by
\begin{equation}\label{fun sup dis}
\Psi_\nu(p)= \big< p-\nu, N(p) \big>, \hspace{0.5cm} \Lambda_\nu(p)= \big| p-\nu\big|^2, \hspace{0.5cm} p \in \Sigma,
\end{equation}
where $ \left < \hspace{0.1cm}, \hspace{0.1cm}\right >$ denotes the Euclidean scalar product in $\mathbb{R}^n$, are called the \textit{support function }and \textit{quadratic distance function} with respect to $\nu \in \mathbb{R}^n$, respectively. Geometrically, $\Psi_\nu(p)$ measures the signed distance from $\nu$ to the tangent plane $T_p\Sigma$ and $\Lambda_\nu(p)$ measures the square of the distance from $p$ to $\nu$. If $\nu$ is the origin, we write $\Psi_0 = \Psi$ and $\Lambda_0 = \Lambda$.

\subsection{Helmholtz Equation}

The reduced Helmholtz equation is an elliptic differential equation describing phisical phenomena related to oscillatory problems. Considering a function  $h: \mathbb{R}^2 \rightarrow \mathbb{R}$, the two-dimensional \textit{Helmholtz equation} for $h$ is defined by
$$
\Delta h(z) + c(K(z))^2 h(z) = 0,
$$
where $ K $ is a function and $c$ is a real non-zero constant. In \cite{Generalized}, the authors introduce the generalized Helmholtz equation for a function $h$ as 
$$
\Delta \left( \frac{1}{(K(z))^2}(\Delta h(z) + c(K(z))^2 h(z)\right) = 0
$$
where $K$ is a function and $c$ is a real non-zero constant. Note that every solution of the Helmholtz equation  is a solution for the generalized Helmholtz equation. 

The next result is the Theorem $1$ from \cite{Generalized} which provides explicit solutions for the generalized Helmholtz equation depending on three holomorphic functions.

\begin{theorem}\label{Helm G}
Let $g$ be a holomorphic function, $c$ a real non-zero constant and $K = (2\sqrt{2}|g'|)/(c + |g|^2)$. In this case, the functions
 
$$
h = \frac{\left <1, A \right > + \left <g, B \right > }{1 + |g|^2}	
$$

\hspace{-0.3cm}are solutions of the generalized Helmholtz equation, where $A$, $B$ are holomorphic functions.

\end{theorem}

The following corollary comes from the aforementioned work and it brings conditions for a solution of the generalized Helmholtz equation to be solved from the Helmholtz equation.

\begin{corollary}\label{Helm}
Let $g$ be a holomorphic function, $c$ a real non-zero constant and $K = (2\sqrt{2}|g'|)/(c + |g|^2)$. In this case, the functions

\vspace{0.3cm}
$$
h = \frac{\left <1, A \right > + \left <g, B \right > }{1 + |g|^2}	
$$

\hspace{-0.3cm} are solutions of the Helmholtz equation when $A$  is a holomorphic function and $B$ is a holomorphic function such that $B = \frac{1}{c} \int (A'g - Ag' + ic_1g') \hspace{0.1cm}dz$, where $c_1$ is a real constant.

\end{corollary}

In \cite{n dimensional}, Corro and Rivero introduce the $n$-dimensional generalized Helmholtz and $n$-dimensional Helmholtz equation in the same way as above, considering now the function $h: \mathbb{R}^n \rightarrow \mathbb{R}$ in $n$-variables. It is provided a class of solutions to them in terms of biharmonic functions and they are used to describe classes of generalized Weingarten hypersurfaces.

For the $2$-dimensional case, they define the \textit{harmonic generalized Weingarten surface depending on support function and the radius function} (in short, RSHGW-surface) as the ones satisfing

$$
\Delta \left(\frac{R}{d}\left(\frac{H}{K} - R - \Psi\right)\right) = 0,
$$
and the \textit{generalized Weingarten surface depending on support function and the radius function} (in short, RSGW-surface) as those surfaces satisfing

$$
\frac{H}{K} - R - \Psi = 0.
$$

\vspace{0.3cm}
Using the result in Corollary \eqref{Helm G}, the authors characterizes the RSHGW-surfaces in terms of functions $h$ which are solutions for the generalized Helmholtz equation for  $K = (2\sqrt{2}|g'|)/(c + |g|^2)$. In the same way, using the solutions $h$ in Corollary \eqref{Helm} for the Helmholtz equation  when $K = (2\sqrt{2}|g'|)/(c + |g|^2)$, the RSGW-surfaces are characterized. For the rotational RSHGW-surfaces, they conclude that the solutions $h$ assume the form

$$
h = \frac{a_2 + c_1u + e^{2u}(a_3 + c_2u)}{1+e^{2u}}, \hspace{1.0cm}a_2, a_3, c_1, c_2 \in \mathbb{R}
$$

\hspace{-0.3cm} and for the RSGW-surfaces of rotation, the solutions $h$ are

$$
h = c_1 - \big(a_2 +c_1(u-1)\big)\tanh u, \hspace{0.5cm}c_1, a_2 \in \mathbb{R}
$$

\vspace{1.0cm}

\subsection{Holomorphic Functions}

The identification of the complex plane $\mathbb{C}$ with $\mathbb{R}^2$ naturally induces the notion of inner product in the space of  holomorphic functions. For $f,g : U \subset \mathbb{C} \rightarrow \mathbb{C}$ holomorphic functions, the inner product $ \left< f,g \right >$ is a real function defined in $U$, given by
$$
\left< f,g \right > = \left< 1, f \right >\left< 1,g \right > + \left< i, f \right >\left< i, g \right >,
$$
where $ \left< 1, f \right >= Re(f) $ and $  \left< i , f \right > = Im(f)$ denote the real and imaginary parts of $f$, respectively. Moreover, the norm of a holomorphic function $f : U \subset \mathbb{C} \rightarrow \mathbb{C}$ is defined as
$$
|f|= \sqrt{\left< f,f \right>}.
$$

This inner product satisfies the following properties for holomorphic functions $f, g$ and $h$:
\begin{enumerate}\label{prop hol}
\item $\left< f,g \right>_{,1} = \left< f',g \right> + \left< f,g' \right>$. \label{item 1}
\item $ \left< f,g \right>_{,2} = \left< if',g \right> + \left< f, ig' \right>$.\label{item 2}
\item $ \left< fh, g \right> = \left< f, \bar{h}g \right>$.\label{item 3}
\item $\bar{f}g =  \left< f, g \right> + i\left< if, g \right>$.\label{item 7}
\end{enumerate}
where $f'$ denotes the complex derivative of $f$. Using the notation settled in the beginning of the section, the relationship between the real and complex derivatives of a holomorphic function $f$ is
\begin{equation*}
f' = f_{,1} = -if_{,2}.
\end{equation*}

Here we present some results from the theory of holomorphic functions which later will be useful in our work.
\begin{theorem}\label{eq rfh}
Every real harmonic function  defined in an open simply connected set of $ \mathbb{C}$ is the real part of a holomorphic function defined in this set.
\end{theorem}

\begin{proposition}\label{pro cte}
	
Let $f,g,h$ be holomorphic functions. Then the equality
\begin{equation}\label{eq fgh}
\left<1,f \right> + \left< g, h \right> = 0,
\end{equation}	
is valid if and only if there exist real constants $c_1, c_2$ and a complex constant $z_1$ such that

\begin{equation}
\left \{
\begin{array}{ccc}
h & = & \hspace{-0.3cm}ic_1g + z_1 \\
f & = & \hspace{-0.1cm} -\bar{z}_1g + ic_2
\\
\end{array}
\right.
\end{equation}
\end{proposition}


\vspace{1.0cm}
\subsection{Sphere Congruence}
A \textit{sphere congruence} in $\mathbb{R}^{n+1}$ is a $n$-parameter family of spheres whose centers lie on a hypersurface $\Sigma_0$ contained in $\mathbb{R}^{n+1}$ with a differentiable radius function. In other words, if we  consider  $\Sigma_{0}$  locally parameterized by  $X_0: U \subset \mathbb{R}^{n} \rightarrow \mathbb{R}^{n+1}$, then for each point $ u \in U$ there exists a sphere centered at $X_0(u)$ with radius  $R(u)$, where $R$ is a differentiable real function, named \textit{radius function}.

An \textit{envelope of a sphere congruence} is a hypersurface $\Sigma \subset \mathbb{R}^{n+1}$ such that each point of $\Sigma$ is tangent to a sphere of the sphere congruence. Two hypersurfaces $\Sigma$ and $\widetilde{\Sigma}$ are said to be (\textit{locally})\textit{associated by a sphere congruence} if there is a (local) diffeomorphism  $ \varphi: \Sigma \rightarrow \widetilde{\Sigma} $ such that at corresponding points $p$ and $\varphi(p)$ the hypersurfaces are tangent to the same sphere of the sphere congruence. It follows that the normal lines at corresponding points intersect at an equidistant point on the hypersurface $\Sigma_0$. If, moreover, the diffeomorphism $\varphi$ preserves lines of curvature, we say that $\Sigma$ and $\widetilde{\Sigma}$ are associated by a \textit{Ribaucour transformation}.

In the work \cite{Laredo} is established that for a hypersurface $\Sigma$ in $\mathbb{R}^{n+1}$ satisfying
\begin{equation}\label{cond cong}
\left< p - \nu, N(p)\right > \neq 1, \hspace{0.5cm}\forall p \in \Sigma
\end{equation}  
and $\nu \in \mathbb{R}^{n+1}$ fixed, there exists a sphere congruence for which $\Sigma$ and the unit sphere $E$ centered in $\nu$ are envelopes. In this case, the radius function is given by the expression

\begin{equation}\label{eq raio}
R(p) = \frac{1 - |p-\nu|^2 }{2\big( \big< p-\nu, N(p)\big> - 1\big)},
\end{equation}

\vspace{0.3cm}
\hspace{-0.4cm}which shows that $R$ is a geometric invariant of $\Sigma$, in the sense it doesn't depend on the parameterization of hypersurface.  Futhermore, it is proved that a such hypersurface $\Sigma$ can be locally parameterized from a local parameterization of $E$ in a way  described below.

\begin{theorem}\label{Teo Pri}
Let $ \Sigma$ be  a hypersurface in  $ \mathbb{R}^{n+1}$ such that $   \left < p, N(p)\right > \neq 1$, for all $p \in \Sigma$, where $N$ is its  normal Gauss map. For each orthogonal local parameterization $ Y:U\subset \mathbb{R}^n \rightarrow \mathbb{S}^n$ of $\mathbb{S}^{n}$, there is a differentiable function $ h: U \subset \mathbb{R}^{n}  \rightarrow \mathbb{R}$, associated to this parameterization, such that $\Sigma$ can be locally parameterized by 
\begin{equation}\label{eq X}
X(u) = Y(u) -2 \left(\frac{h(u)+c}{S(u)}\right)\eta(u), \hspace{0.5 cm}u \in U,
\end{equation}
\vspace{0.3cm}
where the function $h$ satisfy \hspace{0.1cm}$h(u) \neq 0$, for all $ u \in U$, $c$ is a nonzero real constant and
\begin{equation}\label{ eq S}
\eta = \nabla_{L}h +hY, \hspace{1.0cm}
S = \left< \eta, \eta\right> =  \big|\nabla_Lh\big|^2 + h^2,
\end{equation}	
with $ L_{ij} = \left < Y_{,i}, Y_{,j} \right >$.

In these coordinates, the  Gauss map $ {N}$ of $ \Sigma $ is given by
\begin{equation}\label{eq N}
N(u) = Y(u) - 2\frac{h(u)}{S(u)}\eta(u), \hspace{0.5 cm}u \in U.
\end{equation}
Moreover, the Weingarten matrix $W$ of $\Sigma$ is 
\begin{equation}\label{eq W}
W=  \big[SI - 2hV\big] \big[SI - 2\big(h +c\big)V\big]^{-1},
\end{equation}
where $ V= (V_{ij})$ is given by
\begin{equation}\label{eq V}
V_{ij} = \frac{1}{L_{jj}}\left( h_{,ij} - \sum_{k}^{n}{h_{,k}\Gamma_{ij}^{k}}+ hL_{ij}\delta_{ij} \right ),    
\end{equation}
with $\Gamma_{ij}^{k}$  the Christoffel symbols of the metric $L$ and $I$ is the identity matrix $n \times n$. We also have that $X$ is regular if and only if, 
\begin{equation}\label{eq P}
P = det \hspace{0.1cm}\big[SI -2\big(h +c\big)V \big] \neq 0, \hspace{0.3cm} \mbox{for all} \hspace{0.2cm} u \in U.
\end{equation}

Conversely, if  $Y:U\subset \mathbb{R}^n \rightarrow \mathbb{S}^n$ is an orthogonal local parameterization and $h: U \rightarrow \mathbb{R}$ is a differentiable function that doesn't vanish in any point and satisfies \eqref{eq P}, then \eqref{eq X} defines an immersion in $\mathbb{R}^{n+1}$ with normal Gauss map $N$ given by \eqref{eq N}, Weingarten matrix described by \eqref{eq W} and $\left< X, N \right> \neq 1 $ in all point.
	
\end{theorem}

\vspace{1.0cm}
\textbf{Remark 1.}\label{rotation}
A hypersurface  $\Sigma$ in the conditions above is locally associated to $\mathbb{S}^n$  by a sphere congruence and the function $  h: U \subset \mathbb{R}^{n}  \rightarrow \mathbb{R}$ of which the theorem refers is given by

\begin{equation}\label{eq h}
h(u) = - \frac{c}{R(u)+1}, \hspace{0.5 cm}u \in U,
\end{equation}

\hspace{-0.3cm} where $c$ is a nonzero constant and $R$ is the radius function. 



In the case the matrix $V$ is diagonal, the hypersurface $\Sigma$ is parameterized by lines of curvature and it is  associated to $\mathbb{S}^n$ by a Ribaucour transformation. Finally, if we take $Y = \pi_{-}^{-1}$, where  $\pi_{-}: \mathbb{S}^{n}\backslash\{-e_{n+1}\} \rightarrow \mathbb{R}^{n}$ is the stereographic projection,  $\Sigma$ is a hypersurface of rotation if, and only if, the function $h$ is radial.

\vspace{1.0cm}
Below it follows some properties of hypersurface $\eta$ given in \eqref{ eq S}.

\vspace{1.0cm}
\textbf{Remark 2.}\label{sup eta}
Observe that, in the conditions of  Theorem \eqref{Teo Pri},  $Y$ is a unit normal vector field to the hypersurface $\eta: U \subset \mathbb{R}^{n}  \rightarrow \eta(U)$ given by 
\begin{equation*}
\eta(u)= \nabla_{L}h(u) +h(u)Y(u), \hspace{0.5cm} u \in U,
\end{equation*}
and $  \left< \eta, Y \right> = h $, so that  $h$ is the support function of hypersurface $\eta$. Futhermore, $\eta$ has regularity condition equal to $det \hspace{0.1cm} V \neq 0$ and its Weingarten matrix is $V^{-1}$, with $V$  as in \eqref{eq V}. Indeed, it is valid that
\begin{equation}\label{caju}
\eta_{,j} = \sum_{k}V_{jk}Y_{,k} \hspace{0.5cm}1\leq j \leq n,
\end{equation}
which means that the matrix $V^T$ is the coefficient matrix of $d\eta$ in the base $\{Y_{,i}\}$, so that $\eta$ is an immersion iff $det \hspace{0.1cm} V^T = det \hspace{0.1cm} V \neq 0$. In addition, the equality \eqref{caju} implies that 
$$
\sum_{j}^{n}(V^{-1})_{ij}\eta_{,j} = Y_{,i}    \hspace{0,5cm}1\leq i \leq n,
$$
which is to say that $V^{-1}$ is the Weingarten matrix of $\eta$.

Now, we are going to show a fact that will be useful later in our discussion.

\begin{lemma}\label{eq eta radial}
When $Y$ is the inverse of stereographic projection, the hypersurface $\eta(U)$ is rotational if and only if $h$ is a radial function.
\end{lemma}


\begin{proof}

In fact, if  $\eta(U)$ is rotational, the ortogonal sections to the axis of rotation determine on the surface $(n-1)$-dimensional spheres centered on this axis. Note that along these spheres both $|\eta|^2$ and  the angle between $\eta$ and $Y$ must be constant. Since 
\begin{equation*}
\left< \eta, Y \right> = h, \hspace{1.0cm} \left< \eta, \eta \right> = |\nabla_Lh|^2 + h^2, \hspace{0.5cm}L_{ij} = \left < Y_{,i}, Y_{,j} \right >,
\end{equation*}
we conclude that $h$ and so $|\nabla_Lh|^2$ are constant along these spheres. Taking $Y$ as the inverse of stereographic projection, we get $ L_{ij} = 4\delta_{ij}/(1+|u|^2)^2$, so that
\begin{equation*}
|\nabla_Lh|^2 = \left(\frac{1+|u|^2}{2}\right)^2|\nabla h|,
\end{equation*}
what says that $|u|^2$ is constant as one goes around  the ortogonal sections. Therefore the function $h$ is constant along $(n-1)$-dimensional  spheres in $\mathbb{R}^n$ centered in the origin and, therefore, is a radial function.

On the other hand, if we suppose that $h$ is a radial function, we can write $h(u) = J(|u|^2), u \in U$, for some differentiable function $J$. We set $|u|^2 = t$ and we denote the derivative of $J$ with respect to $t$ as $J'(t)$. Therefore $h_{,i} = 2J'u_i$ and taking the parameterization $Y$ as the inverse of stereographic projection, we get
\begin{equation*}
\eta = \Bigg(\left(J'(1-t)+ \frac{2J}{1+t}\right)u, -2tJ'+ J\left(\frac{(1-t)}{1+t}\right)\Bigg).
\end{equation*}
If $-2tJ'+ J\left(\frac{(1-t)}{1+t}\right)$ is constant, then
$$
\Bigg| \left(J'(1-t)+ \frac{2J}{1+t}\right)u \Bigg|^2 =  \left(J'(1-t)+ \frac{2J}{1+t}\right)^2t,
$$
which means that the ortogonal sections to the axis $x_{n+1}$ determine on $\eta(U)$ $(n-1)$-dimensional spheres centered on this axis, so that  $\eta(U)$ is rotational.

\end{proof}

\vspace{1.0cm}
The next proposition follows the same steps as Theorem \eqref{Teo Pri}, but now in the context of Riemann surfaces. This allows us to work with holomorphic functions, which will enable us to construct  Weierstrass-type representations.
\begin{theorem}\label{teo g}
Let $\Sigma$ be a  Riemann surface and $ X: \Sigma \rightarrow \mathbb{R}^3$ an immersion such that $ \big< X(p), N(p) \big> \neq 1$, for all $p \in \Sigma$, where $N$ is the normal Gauss map of $X$. Consider also a parameterization $Y:U \subset \mathbb{R}^2 \rightarrow \mathbb{S}^2$ of the unit sphere given by $Y = \pi_{-}^{-1}\circ g$, where $g : \mathbb{C} \rightarrow \mathbb{C}_{\infty}$ is a holomorphic function such that $g' \neq 0$ and 
$\pi_{-}^{-1}: \mathbb{C} \rightarrow \mathbb{S}^2 \setminus {\{-e_3\}}$ is the inverse of  stereographic projection.	Then there is a differentiable function $ h: U \subset \mathbb{R}^{2}  \rightarrow \mathbb{R}$ associated to this parameterization, such that $\Sigma$ can be locally parameterized by 
\begin{equation}\label{eq X g}
X=  \frac{1}{T}\Big(2g, 2-T\Big) - \frac{2(h+c)}{S}\eta,
\end{equation}
\vspace{0.3cm}
where $c$ is a nonzero real constant, $T= 1+|g|^2$  and
\begin{equation*}
\eta = \nabla_{L}h +hY, \hspace{1.0cm}
S = \left< \eta, \eta\right> =  \big|\nabla_Lh\big|^2 + h^2,
\end{equation*}	
with 
\begin{equation}
L_{ij} = 
\left < Y_{,i}, Y_{,j} \right > = \frac{4|g'|^2}{T^2}\delta_{ij}, \hspace{0.5cm} T = 1+|g|^2, \hspace{0.5cm}   1\leq i,j \leq 2, 	
\end{equation}
	
In these coordinates, the  Gauss map $ {N}$ of $ \Sigma $ is given by
\begin{equation}\label{eq N g}
N = \frac{1}{T}\Big(2g, 2-T\Big) - 2\frac{h}{S}\eta.
\end{equation}
Moreover, the Weingarten matrix $W$ of $\Sigma$ is $W=  \big[SI - 2hV\big] \big[SI - 2\big(h +c\big)V\big]^{-1} $, where the matrix $V$  is such that
\begin{equation}\label{eq V11}
V_{11}= \frac{1}{L_{11}}\left[h_{,11}- \left< \frac{g''}{g'}- \frac{2}{T}g' \bar{g} \hspace{0.1cm}, \hspace{0.1cm} \nabla h \right> +hL_{11}\right], 
\end{equation}
	
\begin{equation}\label{eq V12}
V_{12}= \frac{1}{L_{22}}\left[h_{,12}- \left< i\Big(\frac{g''}{g'}- \frac{2}{T}g' \bar{g}\Big) \hspace{0.1cm}, \hspace{0.1cm} \nabla h \right> \right],
\end{equation}
	
\begin{equation}\label{eq V22}
V_{22}= \frac{1}{L_{22}}\left[h_{,22} + \left< \frac{g''}{g'}- \frac{2}{T}g' \bar{g} \hspace{0.1cm}, \hspace{0.1cm} \nabla h \right> +hL_{22}\right].
\end{equation}
	
We also have that $X$ is regular if and only if 
\begin{equation}\label{eq P g}
P = S^2 -2(h+c) \hspace{0.1cm} S\hspace{0.1cm} trV + 4(h+c)^2 detV \neq 0.
\end{equation}
	
\end{theorem}

\begin{proof}
Taking $Y:U \subset \mathbb{R}^2 \rightarrow \mathbb{S}^2$ as in the statement, we have, for $ u \in U$,

\begin{equation}\label{eq Y g}
Y(u)  =   \frac{1}{1+|g(u)|^2}\Big(2g(u), 1-|g(u)|^2\Big).
\end{equation} 
	
Thus,  from Theorem \eqref{Teo Pri}  there is a differentiable function $ h: U\subset \mathbb{R}^2 \rightarrow \mathbb{R}$ such that $X(\Sigma)$ can be locally parameterized by \eqref{eq X g} with Gauss map \eqref{eq N g}.
	
Theorem \eqref{Teo Pri} also ensures that the Weingarten matrix of $X(\Sigma)$ is  $W=  \big[SI - 2hV\big] \big[SI - 2\big(h +c\big)V\big]^{-1} $, with  regularity condition given by $P \neq 0 $, where $P = det \hspace{0.1cm}\big[SI -2\big(h +c\big)V \big]$ with $V$  as in \eqref{eq V}. For $n=2$, this regularity condition may be rewritten as the equation \eqref{eq P g}.
	
In order to make explicit the $V$'s entries, let us find the Christoffel symbols of the metric $L_{ij}= \left< Y_{,i}\hspace{0.1cm}, \hspace{0.1cm} Y_{,j} \right>, \hspace{0.1cm} 1 \leq i, j \leq 2$. From the expression \eqref{eq Y g} for $Y$, we have that the metric $L$ is diagonal given by
\begin{equation}\label{eq L}
L_{ij} = 
\frac{4|g'|^2}{T^2}\delta_{ij}, \hspace{0.5cm} T = 1+|g|^2, \hspace{0.5cm}   1\leq i,j \leq 2, 	
\end{equation}
and from equations \eqref{eq Chris},  we have that its Christoffel symbols are 
$$
\Gamma_{11}^{1} =  \frac{T\left<g', g''\right>- 2|g'|^2\left<g, g'\right>}{T|g'|^2},  \hspace{1.5cm} \Gamma_{22}^{2} = \frac{T\left<g', ig''\right>- 2|g'|^2\left<g, ig'\right>}{T|g'|^2}. 
$$	
	
From \eqref{eq V}, it follows that
\begin{eqnarray}\label{V}
V_{11} &= &\frac{1}{L_{11}}\left[h_{,11}- \left< \Gamma_{11}^{1} +i\Gamma_{11}^{2} \hspace{0.1cm}, \hspace{0.1cm} \nabla h \right> +hL_{11}\right], \nonumber\\
V_{12} &=& \frac{1}{L_{22}}\left[h_{,12}- \left< \Gamma_{12}^{1} + i \Gamma_{12}^{2} \hspace{0.1cm}, \hspace{0.1cm} \nabla h \right> \right], \\
V_{22} &=& \frac{1}{L_{22}}\left[h_{,22} - \left< \Gamma_{22}^{1} + i\Gamma_{22}^{2}\hspace{0.1cm}, \hspace{0.1cm} \nabla h \right> +hL_{22}\right]. \nonumber
\end{eqnarray}
	
Using the property \eqref{item 3} for holomorphic functions, we have
\begin{eqnarray}\label{V11}
\Gamma_{11}^{1} +i\Gamma_{11}^{2} &=& \left< \frac{1}{\bar{g'}} , g'' \right> - \frac{2}{T}\left<g , g' \right> + i\left( \frac{2}{T}\left<g , ig' \right> - \left< \frac{1}{\bar{g'}} , ig'' \right> \right) \nonumber \\ \nonumber
&=& \left< 1 , \frac{g''}{g'} - \frac{2}{T}g'\bar{g} \right> +i\left< i , \frac{g''}{g'} - \frac{2}{T}g'\bar{g} \right>  \nonumber \\
&=& \frac{g''}{g'} - \frac{2}{T}g'\bar{g}
\end{eqnarray}
Similarly,
\begin{equation}\label{acerola}
\Gamma_{12}^{1} + i \Gamma_{12}^{2} = i\left( \frac{g''}{g'} - \frac{2}{T}g'\bar{g} \right), \hspace{1.0cm} \Gamma_{22}^{1} + i\Gamma_{22}^{2} = \frac{2}{T}g'\bar{g} - \frac{g''}{g'}.
\end{equation}
Substituing \eqref{V11} and \eqref{acerola} in the expressions \eqref{V}  for the $V$'s entries , we are done.	
\end{proof}

\vspace{1.0cm}
\section{$H_1$ and $H_2$ Surfaces}

Since the radius function is a geometric invariant for a hypersurface  $\Sigma$ which is associated to a unit sphere by a sphere congruence, we can define certain curvatures for it. Thus,  for such a  $\Sigma$, we define its \textit{ spherical radial curvatures $s_i$} and \textit{ spherical mean curvature $H_S$} as
\begin{equation}\label{eq curv esf}
s_{i}= \frac{1+k_i}{1-k_iR}, \hspace{1.0cm} H_S = \frac{1}{n}\sum_{i=1}^{n}s_i,
\end{equation}
\vspace{0.3cm}
\hspace{-0.4cm}
where  $k_i$ are the principal curvatures of $\Sigma$, \hspace{0.1cm} $1\leqslant i \leqslant n$, and $R$ is the radius function given by \eqref{eq raio}.

Next we define special classes of surfaces which elements are all envelopes associated to $\mathbb{S}^2$ by a sphere congruence.

\begin{definition}
Let $\Sigma$ be a  surface and $ X: \Sigma \rightarrow \mathbb{R}^3$ an immersion such that $ \big< X(p), N(p) \big> \neq 1$, for all $p \in \Sigma$, where $N$ is the normal Gauss map of $X$. The surface  $\Sigma$ is called a \textit{surface of null spherical mean curvature}, in short $H_1$-surface, if holds
$$
H_S = 0,
$$
and $\Sigma$ is called a \textit{surface with spherical mean curvature of harmonic type}, in short $H_2$-surface,  if it satisfies  

\begin{equation}
\Delta_\sigma \left[\frac{H_S}{\Psi - 1} \right] = 0,
\end{equation}

\vspace{0.3cm}
\hspace{-0.4cm}
where $H_S$ is the spherical mean curvature of $\Sigma$ and $\sigma = I + 2R \hspace{0.1cm}II + R^2 \hspace{0.1cm}III$, with $I, II, III$  the fundamental forms of $\Sigma$.

\end{definition}

\begin{lemma}
Let $X$ be a hypersurface associated to the unit sphere $Y$ by a sphere congruence. Then the quadratic form $\sigma =  = I + 2R \hspace{0.1cm}II + R^2 \hspace{0.1cm}III$, with $I, II, III$  the fundamental forms of $X$, is conformal to the first fundamental form of $Y$.

\begin{proof}
By sphere congruence, it holds 

$$
X + RN = Y + RY = (1+R)Y
$$	
where $R$ is the radius function. Differentiating	the equality above, we have

$$
dX + dR N + R dN = dR Y + (1+R)dY
$$

Taking the norm squared on both sides, we get

$$
I +  2R \hspace{0.1cm}II + R^2 \hspace{0.1cm}III = (1+R^2) \left<dY, dY \right>
$$
which completes the proof.
	
\end{proof}
	
\end{lemma}
Because of the result above and considering $Y$ the parameterization given by the inverse of stereographic projection,  a function is harmonic with respect to the quadratic form $\sigma$ if and only if is harmonic in the metric of sphere, what in turn is conform to the Euclidean metric.

The next result characterizes the $H_1$-surfaces  in terms of solutions of a Helmholtz equation and the $H_2$-surfaces in terms of solutions of a generalized Helmholtz equation.

\begin{theorem}\label{H1 H2}
Let $ X: \Sigma \rightarrow \mathbb{R}^3$ an immersion as in Theorem \eqref{teo g}.  Then 

\vspace{0.3cm}
\begin{enumerate}
\item $ \Sigma $ is a $H_1$-surface if and only if $h$ satisfies the Helmholtz equation

\begin{equation}\label{eq H1}
\Delta h + \frac{8|g'|^2}{T^2} h = 0, \hspace{0.5cm} T = 1+|g|^2,
\end{equation}
 
\vspace{0.3cm}
\item $\Sigma$ is a $H_2$-surface if and only if $h$ satisfies the generalized Helmholtz equation

\begin{equation}\label{morango}
\Delta \left[ \frac{1}{\left(\frac{4 |g'|^2}{T^2}\right)}\left( \Delta h + \frac{8|g'|^2}{T^2} h \right) \right] = 0, \hspace{1.0cm} T = 1+|g|^2
\end{equation}
\end{enumerate}

\vspace{0.1cm}
\hspace{-0.2cm}
where in both cases $g$ is a holomorphic function such that $g' \neq 0$.

\end{theorem}


\begin{proof}
Consider the immersion $X$ described as in equation \eqref{eq X g}. We are going to show that $tr V$, for $V$ as in \eqref{eq V}, can be expressed in two distinct ways.

At first, as $W = \big[SI - 2hV\big] \big[SI - 2\big(h +c\big)V\big]^{-1}$ is the Weingarten matrix for $X$, it follows that the principal curvatures $k_i$ of $X$ are given by

$$
- k_i = \frac{S -2h \sigma_i}{S - 2(h + c)\sigma_i}, \hspace{0.5cm}
$$
\vspace{0.3cm}

where $\sigma_i$, $ i=1,2$, are the eigenvalues of matrix $V$. Thus

$$
\sigma_i = \frac{S}{2h}\left( \frac{1+k_i}{1- k_i R}\right)
$$

where $R$ is the radius function. Therefore, 

$$
tr V = \frac{S}{h}H_S
$$
\vspace{0.3cm}

Note that from equations \eqref{eq X g} and \eqref{eq N g}, it holds

$$
|X-N|^2 = \frac{4c^2}{S^2} \left < \eta, \eta \right > = \frac{4c^2}{S}, \hspace{0.5cm} c \in \mathbb{R}^{*},
$$

and we get

\begin{equation*}
tr V = \left(\frac{4c^2}{h |X-N|^2}\right) H_S
\end{equation*}

\vspace{0.3cm}
Now, noting that $|X-N|^2 = \Lambda -2\Psi +1$ and recalling that $ h = - c/(R+1)$ and $ R = (1-\Lambda)/\big(2(\Psi-1)\big)$, we find lastly

\begin{equation}\label{tr V 1}
tr V = \left(\frac{2c}{\Psi - 1}\right) H_S
\end{equation}

\vspace{0.3cm}
On the other hand, by expressions \eqref{eq V11} and \eqref{eq V22} for matrix $V$, it follows

\vspace{0.3cm}
\begin{equation}\label{tr V 2}
 tr V  = \frac{\Delta h}{L} + 2 h, \hspace{1.0cm} L = \frac{4|g'|^2}{(1+|g|^2)^2}
\end{equation}

\vspace{0.3cm}
Thus, considering the equalities \eqref{tr V 1} and \eqref{tr V 2},  the surface $X$ is a $H_1$-surface or a $H_2$-surface  if and only if the function $h$ satisfies the Helmholtz equation \eqref{eq H1} or 
the generalized Helmholtz equation \eqref{morango}, respectively.

\end{proof}


Next we have a characterization for the $H_1$-surfaces that relies on the Corollary \eqref{Helm} which gives solutions for the Helmholtz equation .

\begin{corollary}\label{H1_h}
Let $\Sigma$ be a surface as in Theorem \eqref{teo g}.  Then $\Sigma$ is a $H_1$-surface if and only if

\begin{equation}\label{h H1}
h(u,v) = \frac{\big<1, A(z) \big> + \big< g(z), B(z) \big> }{1 + |g(z)|^2}, \hspace{0.5cm}z = u+iv,
\end{equation}

\vspace{0.3cm}
\hspace{-0.5cm}
where $A$ is a holomorphic function and $B$ is a holomorphic function such that $ B(z) = \int(A'(z)g(z) - A(z)g'(z) + ic_1g'(z)) \hspace{0.1cm}dz $, for $c_1$ a real constant.
\end{corollary}

The next Corollary is a slight adaptation of Theorem 1 from \cite{Generalized}. It provides solutions for the generalized Helmholtz equation for special functions $K$.

\begin{corollary}\label{H2 h}
Let $\Sigma$ be a surface as in Theorem \eqref{teo g}.  Then $\Sigma$ is a $H_2$-surface if and only if

\begin{equation}\label{h sol}
h(u,v) = \frac{\big<1, A(z) \big> + \big< g(z), B(z) \big> }{1 + |g(z)|^2}, \hspace{0.5cm}z = u+iv,
\end{equation}
where  $A$, $B$ are holomorphic functions. 

\end{corollary}

\begin{proof}

Consider
$$
h = \frac{f}{T}, \hspace{0.5cm} T = 1+|g|^2
$$	
\vspace{0.5cm}
In this case, the Laplacian of $h$ is given by the expression below
	
$$
\Delta h = \frac{\Delta f}{T} + 2 \left <\nabla f, \nabla \left( \frac{1}{T} \right) \right > + f \Delta \left( \frac{1}{T} \right) 
$$
	
As $ T = 1 + |g|^2$, we get
	
\begin{eqnarray*}
\Delta h & = & \frac{\Delta f}{T} - 4 \left <\nabla f, \frac{g \bar{g'}}{T^2} \right > + f \left( \frac{8 |gg'|^2}{T^3} - \frac{4 |g'|^2}{T^2}\right) \\
& = & \frac{\Delta f}{T} - 4 \left <\nabla f, \frac{g \bar{g'}}{T^2} \right > + 4 f |g'|^2 \left( \frac{1}{T^2} - \frac{2}{T^3}\right).
\end{eqnarray*}

Now this equation can be rewrite as
	
\begin{equation}
\Delta h + \frac{8 |g'|^2}{T^2}h = \frac{|g'|^2}{T^2} \left(T \frac{\Delta f}{|g'|^2} - 4 \left < \nabla f, \frac{g}{g'} \right > + 4f \right)
\end{equation}
\vspace{0.5cm}
	
Therefore, the function $ h = f/T$ is a solution of generalized Helmholtz equation 
	
\begin{equation}
\Delta \left\{ \frac{1}{\left(\frac{4 |g'|^2}{T^2}\right)}\left( \Delta h + \frac{8|g'|^2}{T^2} h \right) \right\} = 0
\end{equation}
if and only if 
	
$$
\Delta \left\{ T \frac{\Delta f}{|g'|^2} - 4 \left < \nabla f, \frac{g}{g'} \right > + 4f \right \} \hspace{0.1cm} = \hspace{0.1cm} T \Delta \left( \frac{\Delta f}{|g'|^2}\right) = 0
$$
\vspace{0.5cm}
	
On the other hand, the solutions of equation $\Delta \left( \frac{\Delta f}{|g'|^2}\right) = 0$ are given by $ f = \left < 1, A \right > + \left < g, B \right >$, with $A$, $B$ holomorphic functions. Thus, we are done.
	
\end{proof}


\begin{proposition}\label{Minimal}
In the conditions of Theorem \eqref{teo g}, $\Sigma$ is a $H_1$-surface if and only if $\eta$ is a minimal surface.
\end{proposition}

\begin{proof}
From Remark $2$, we have $\widetilde{W} = V^{-1}$ is the Weingarten matrix of surface $\eta$. Let $-\tilde{k_i}$ be the eigenvalues of $\widetilde{W}$. Then

$$
-\tilde{k_i} = \frac{1}{\sigma_i} \hspace{0.5cm} \Rightarrow \hspace{0.5cm} \sigma_i = - \frac{1}{\tilde{k_i}}
$$
for $\sigma_i$ the eigenvalues	of $V$. Thus,

$$
tr V = \sigma_1 + \sigma_2 = -2 \frac{\widetilde{H}}{\widetilde{K}}
$$
where $\widetilde{H}$ and $\widetilde{K}$ are the mean and Gaussian curvatures of $\eta$, respectively. From equality \eqref{tr V 2}, we conclude

$$
\frac{\Delta h}{L} + 2 h = -2\frac{\widetilde{H}}{\widetilde{K}}, \hspace{1.5cm} L = \frac{4|g'|^2}{(1 + |g|^2)^2}
$$	

\vspace{0.3cm}
\hspace{-0.4cm}	
so that $h$ is a solution of Helmholtz equation \eqref{eq H1} if and only if $\widetilde{H} = 0$.	
\end{proof}


\begin{proposition}\label{Laguerre}
In the conditions of Theorem \eqref{teo g}, $\Sigma$ is a $H_2$-surface if and only if $\eta$ is a Laguerre minimal surface.
\end{proposition}

\begin{proof}
	
From the last demonstration, we have

$$
\frac{\Delta h}{L} + 2 h = -2\frac{\widetilde{H}}{\widetilde{K}}, \hspace{1.5cm} L = \frac{4|g'|^2}{(1 + |g|^2)^2}.
$$	
	
\vspace{0.3cm}
Let $I_Y$ be the first fundamental form  in local coordinates of sphere. Considering the parameterization $Y = \pi_{-}^{-1}\circ g$, where $g : \mathbb{C} \rightarrow \mathbb{C}_{\infty}$ is a holomorphic function such that $g' \neq 0$ and 
$\pi_{-}^{-1}: \mathbb{C} \rightarrow \mathbb{S}^2 \setminus {\{-e_3\}}$ is the inverse of  stereographic projection, we have
	
$$
I_Y = \frac{4|g'|^2}{T^2}\delta_{ij}, \hspace{1.0cm} T = 1+|g|^2.
$$
	
Now recall  that $Y$ is a unit normal vector field to the surface $\eta$. Thus if $\widetilde{III}$ stands for the third fundamental form for $\eta$, then
	
$$
\widetilde{III} = \left < dY, dY \right > = I_Y
$$
	
Therefore
	
$$
X \hspace{0.3cm}\mbox{is} \hspace{0.3cm} H_2\mbox{-surface} \hspace{0.3cm}\Leftrightarrow \hspace{0.3cm}\Delta_{I_Y}\left[ \frac{\Delta h}{L} + 2 h\right] = 0 \hspace{0.3cm}\Leftrightarrow \hspace{0.3cm} \Delta_{I_Y}\left(\frac{\widetilde{H}}{\widetilde{K}}\right) = 0 \hspace{0.3cm}\Leftrightarrow \hspace{0.3cm} \Delta_{\widetilde{III}}\left(\frac{\widetilde{H}}{\widetilde{K}}\right) = 0
$$
this last condition indicating $\eta$ be Laguerre minimal. 
	
\end{proof}

\textbf{Remark 3} 
\hspace{0.2cm} From Theorem \eqref{teo g} and Corolaries \eqref{H1_h} and \eqref{H2 h}, we get a Weierstrass type representation for the  $H_1$-surfaces and $H_2$-surfaces depending on three holomorphic functions given by
 
\begin{equation}
X=  \frac{1}{T}\Big(2g, 2-T\Big) - \frac{2(h+c)}{S}\eta,
\end{equation}

\vspace{0.3cm}
\hspace{-0.4cm}
where $c$ is a nonzero real constant, $T= 1+|g|^2$, \hspace{0.1cm} $ S = \left< \eta, \eta\right> =  \big|\nabla_Lh\big|^2 + h^2 $  and

\begin{eqnarray}\label{eta WR}
\eta \hspace{0.1cm} = \hspace{0.1cm} \nabla_{L}h +hY \hspace{0.1cm} = \hspace{0.1cm} \left(\frac{T}{2}\frac{\nabla h}{|g'|^2}g' - g\left< \nabla h, \frac{g}{g'} \right> + \frac{2h}{T}g \hspace{0.1cm}, \hspace{0.1cm} \frac{(2-T)}{T}h - \left< \nabla h, \frac{g}{g'} \right>  \right).
\end{eqnarray}

\vspace{0.3cm}
Thus, for $h$ given as in Corollary \eqref{H1_h}, $X$ is a Weierstrass type representation for the $H_1$-surfaces, whereas for $h$ given as in Corollary \eqref{H2 h}, $X$ is a Weierstrass type representation for the $H_2$-surfaces. These surfaces generically have singularities given by the expression \eqref{eq P}.

\vspace{0.5cm}
\textbf{Remark 4} \hspace{0.2cm} On the other hand, from Proposition \eqref{Minimal}, the expression \eqref{eta WR} above is an alternative Weierstrass representation for the minimal surfaces when the function $h$ is given as in Corollary \eqref{H1_h}. In the same way, from Proposition \eqref{Laguerre} we conclude that the expression \eqref{eta WR} is a Weierstrass representation for the Laguerre minimal surfaces when $h$ is given as in Corollary \eqref{H2 h}. In both cases the regularity condition is expressed by $det \hspace{0.1cm}V \neq 0$, where $V$ is the matrix which entries are given by \eqref{eq V11}, 
\eqref{eq V12} and \eqref{eq V22}.


\subsection{Examples of $H_1$-surfaces and Minimal Surfaces}

\begin{example}
Considering $A(z) = e^z$, $g(z) = z$ in Corollary  \eqref{H1_h}, we get $B(z)= e^z(z-2) + ic_1z $. The correspondent $H_1$-surface and   $\eta$-minimal surface are drawn below 

\begin{figure}[h]
	\centering
	\subfigure{\includegraphics[scale=0.30]{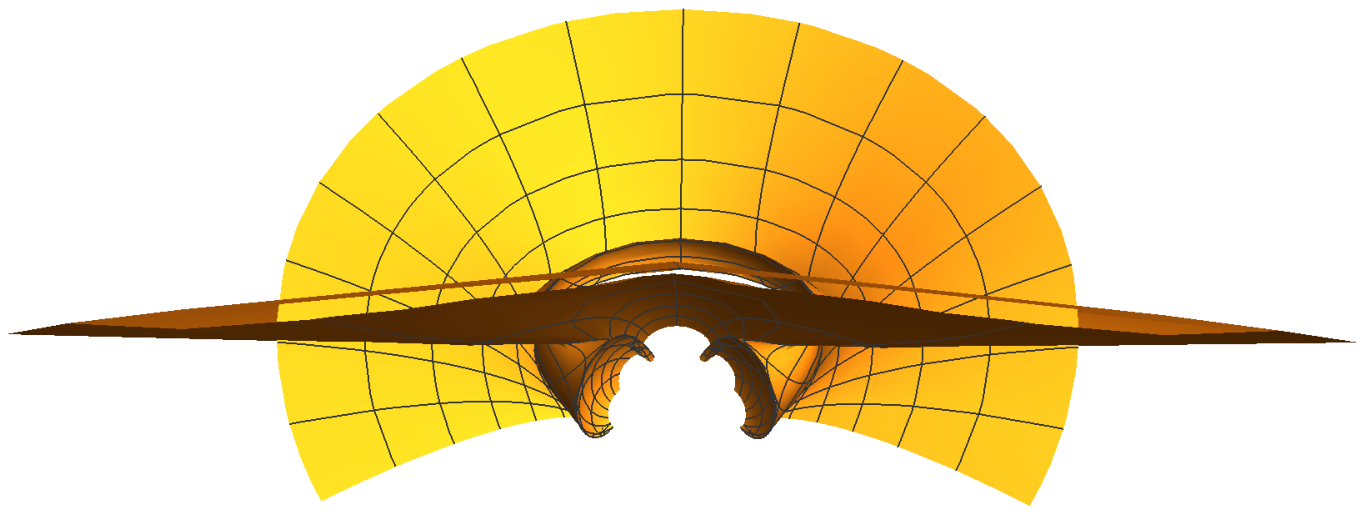}}
	\subfigure{\includegraphics[scale=0.30]{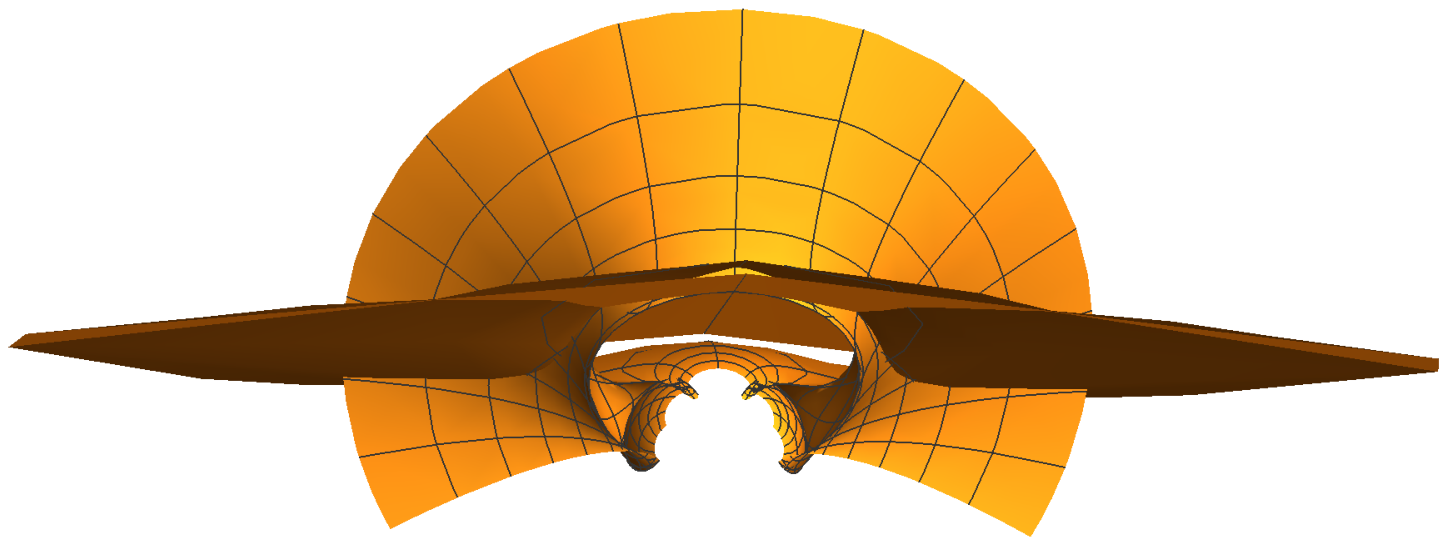}}
\caption{$H_1$-surface for $g(z) = z$ and $A(z) = e^z$ }
\end{figure}
\vspace{-0.5cm}
\begin{figure}[h]
\centering
\subfigure{\includegraphics[scale=0.55]{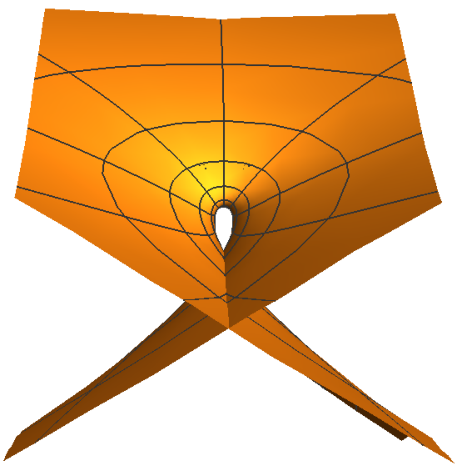}}\hspace{1.5cm}
\subfigure{\includegraphics[scale=0.50]{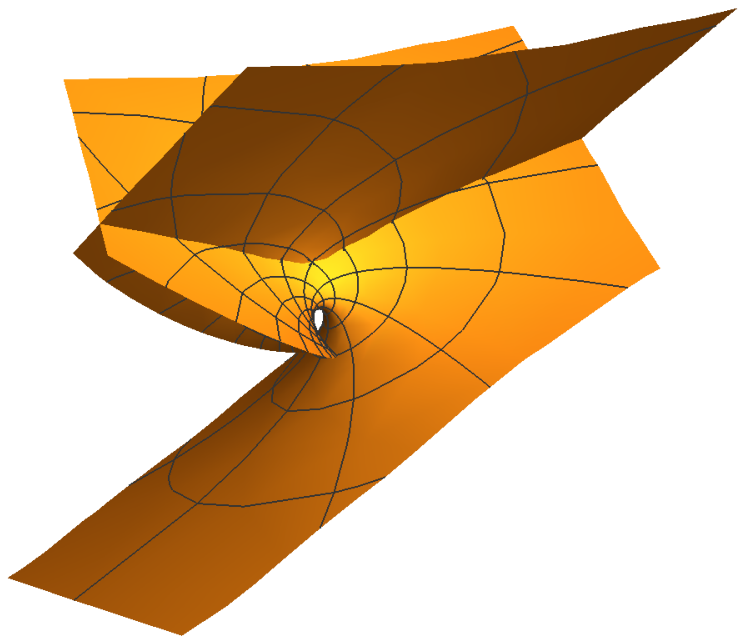}}
\caption{Minimal surface for $g(z) = z$ and $A(z) = e^z$ }
\end{figure}

\end{example}

\newpage
\begin{example}
Considering $A(z) = z$, $g(z) = z^2$, then $B(z)= - \frac{1}{3}z^3 + i c_1z^2$ in Corollary  \eqref{H1_h}. The correspondent $H_1$-surface and $\eta$ minimal surface are given next.
	
\begin{figure}[h]
\centering
\subfigure{\includegraphics[scale=0.25]{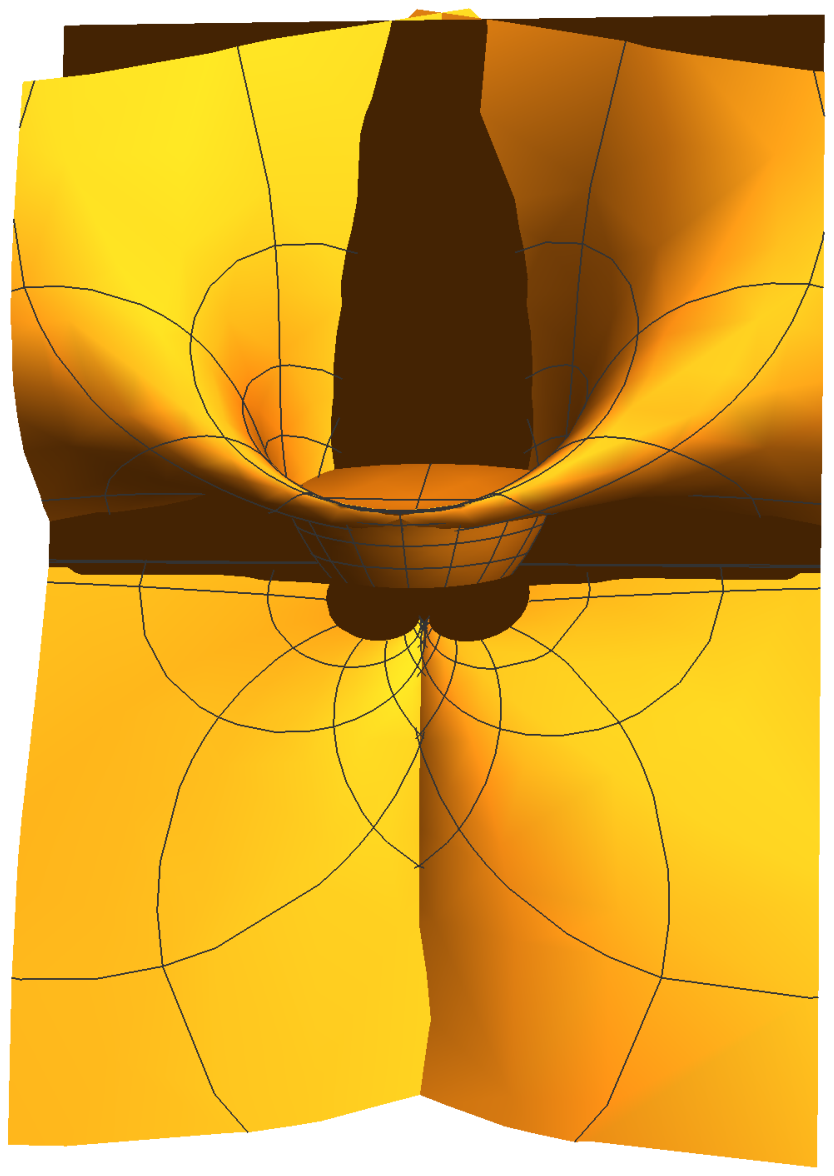}}\hspace{1.5cm}
\subfigure{\includegraphics[scale=0.30]{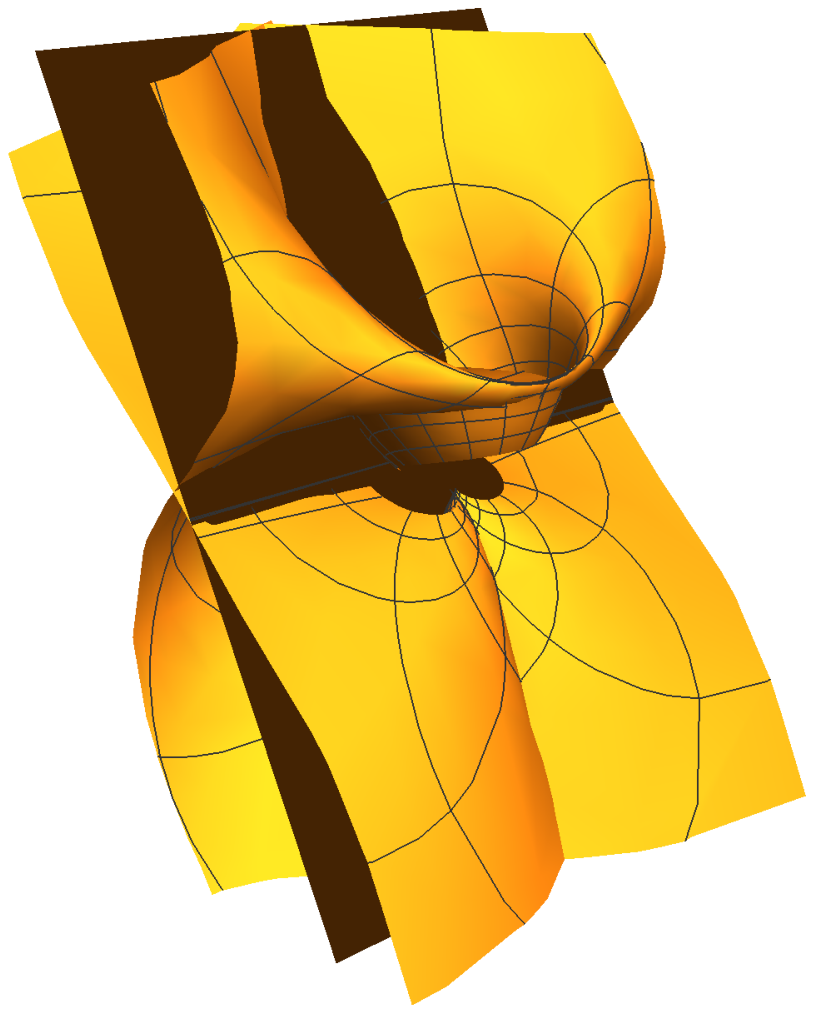}}
\caption{$H_1$-surface for $g(z) = z^2$ and $A(z) = z$ }
\end{figure}
	
\begin{figure}[h]
\centering
\subfigure{\includegraphics[scale=0.30]{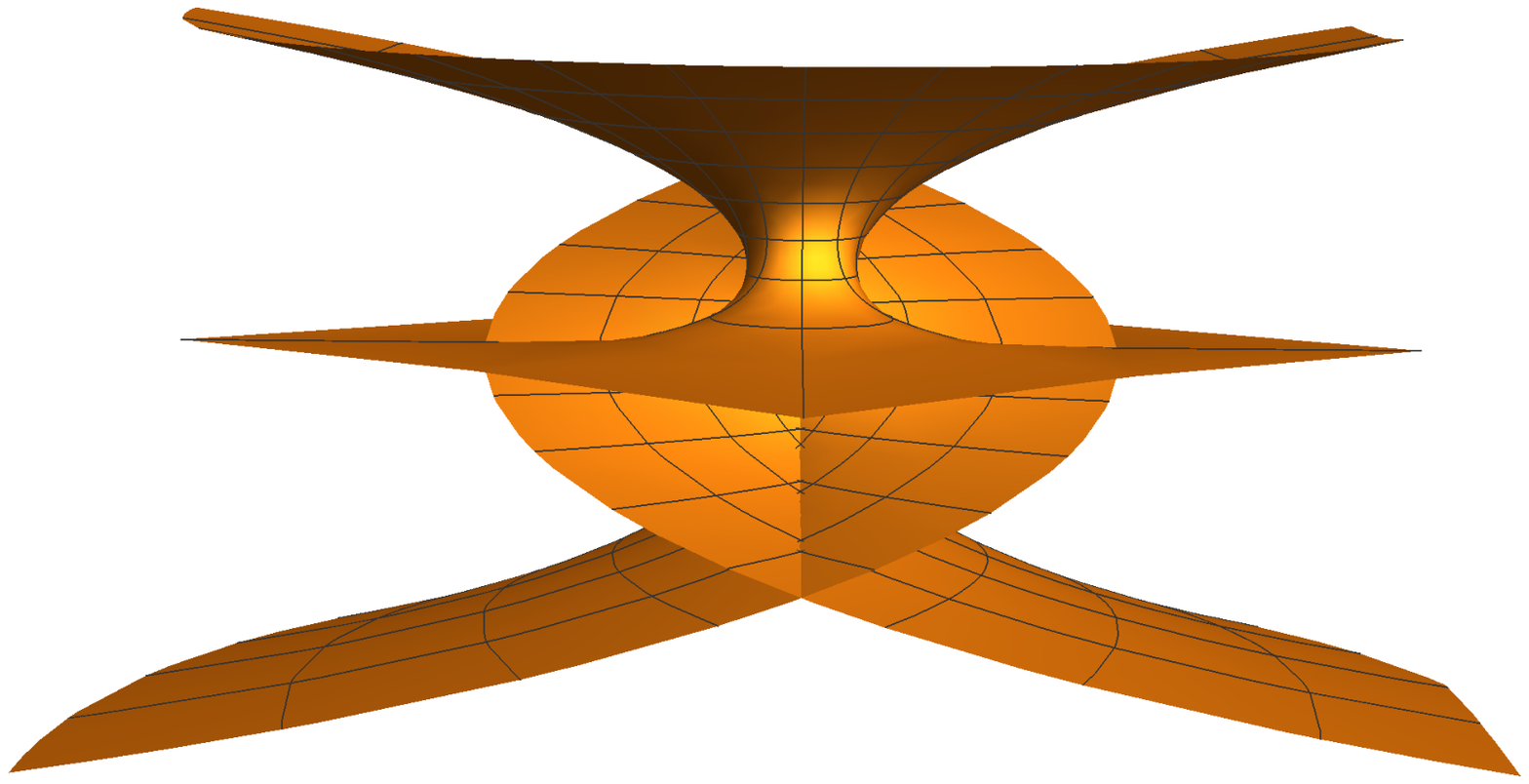}}
\subfigure{\includegraphics[scale=0.30]{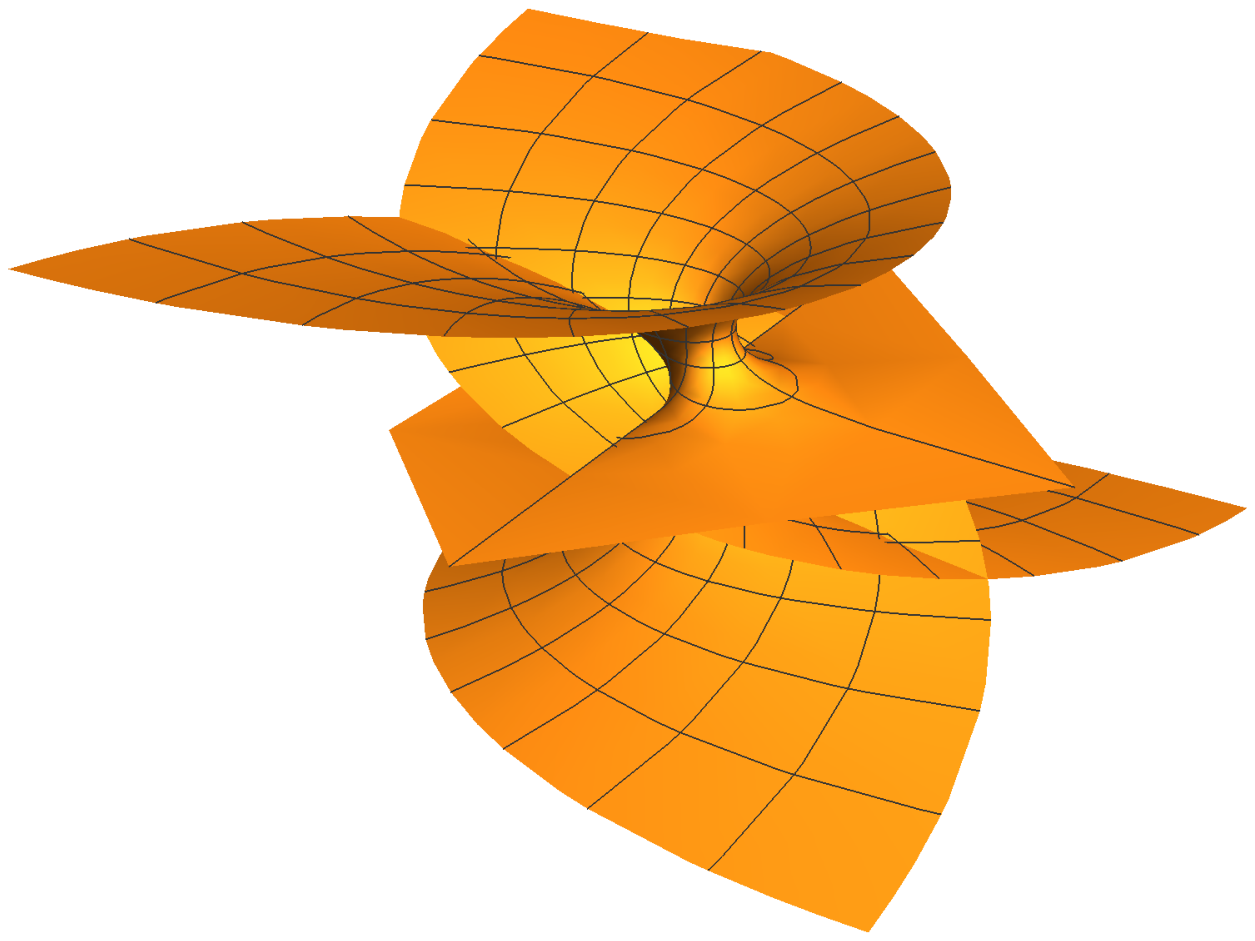}}
\caption{Minimal surface for $g(z) = z^2$ and $A(z) = z$ }
\end{figure}
	
\end{example}

\newpage
\begin{example}
Considering $A(z) = z^2$, $g(z) = z^{-1}$ in Corollary  \eqref{H1_h}, it  follows that $B(z)= 3z +ic_1z^{-1}$. The correspondent $H_1$ and $\eta$ surfaces are drawn below 
		
\begin{figure}[h]
\centering
\subfigure{\includegraphics[scale=0.30]{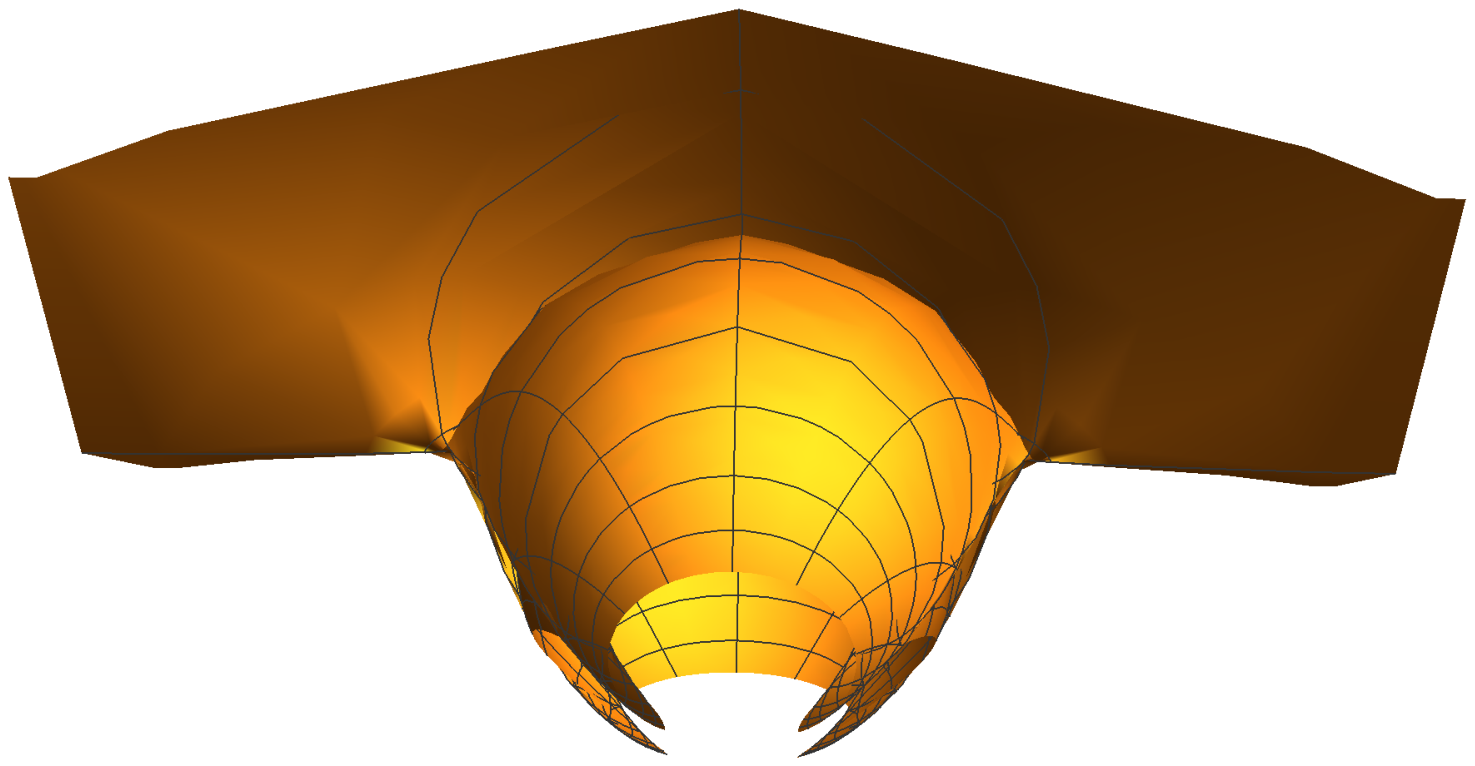}}
\subfigure{\includegraphics[scale=0.35]{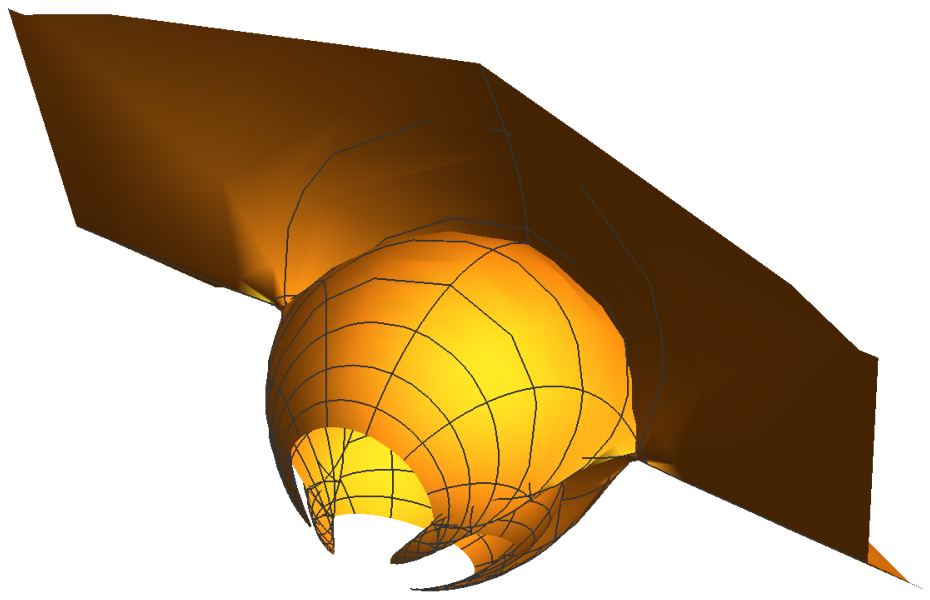}}
\caption{$H_1$-surface for $A(z) = z^2$ and $g(z) = z^{-1}$}
\end{figure}

\begin{figure}[h]
	\centering
	\subfigure{\includegraphics[scale=0.30]{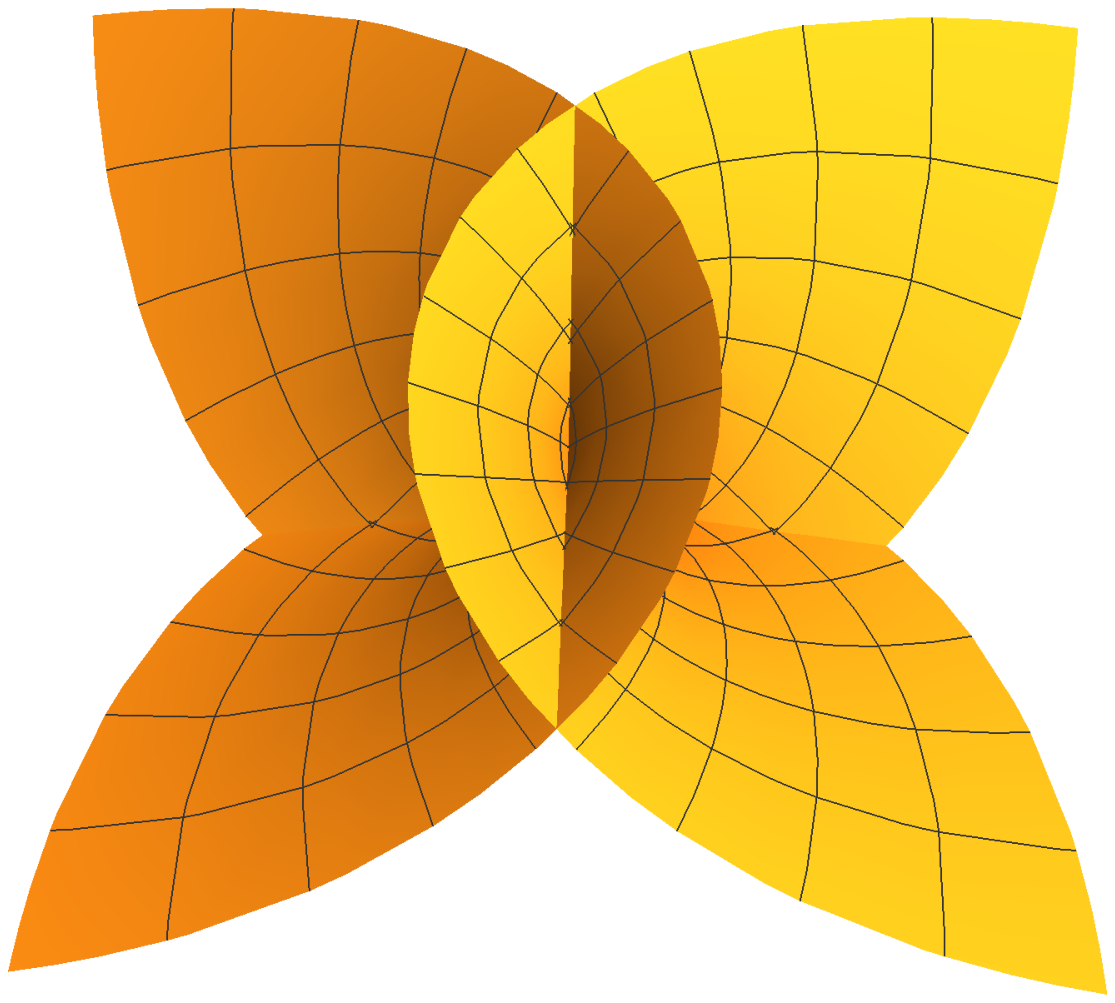}}\hspace{1.5cm}
	\subfigure{\includegraphics[scale=0.40]{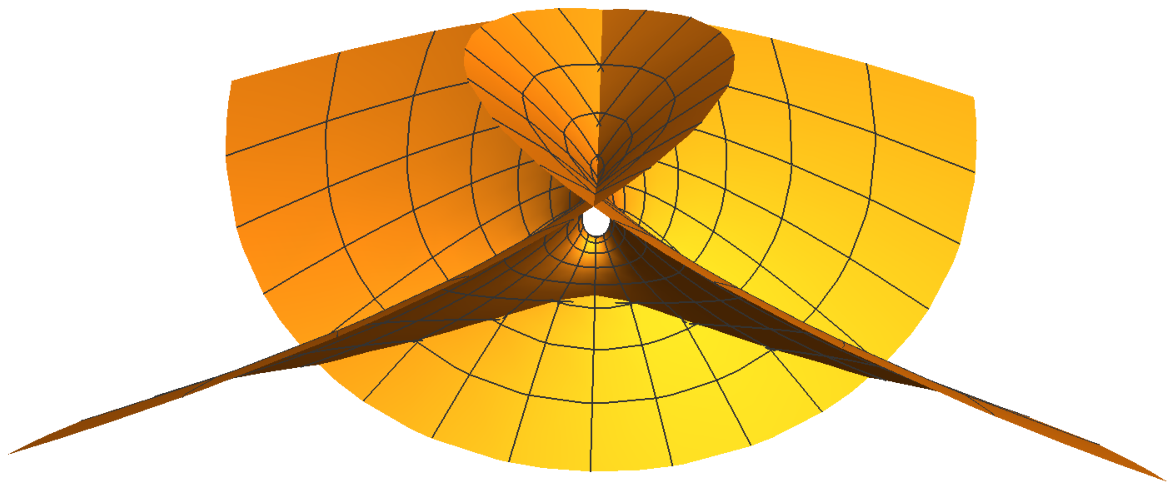}}
	\caption{ Minimal surface for $A(z) = z^2$ and $g(z) = z^{-1}$ }
\end{figure}

\end{example}
\newpage
\begin{example}
Considering $g(z) = z$, $A(z) = \sin z$  in Corollary  \eqref{H1_h}, we have $B(z)= z \sin z + 2 \cos z + ic_1 z$. The correspondent $H_1$ and $\eta$ surfaces are drawn below. 
	
\begin{figure}[h]
\centering
\subfigure{\includegraphics[scale=0.25]{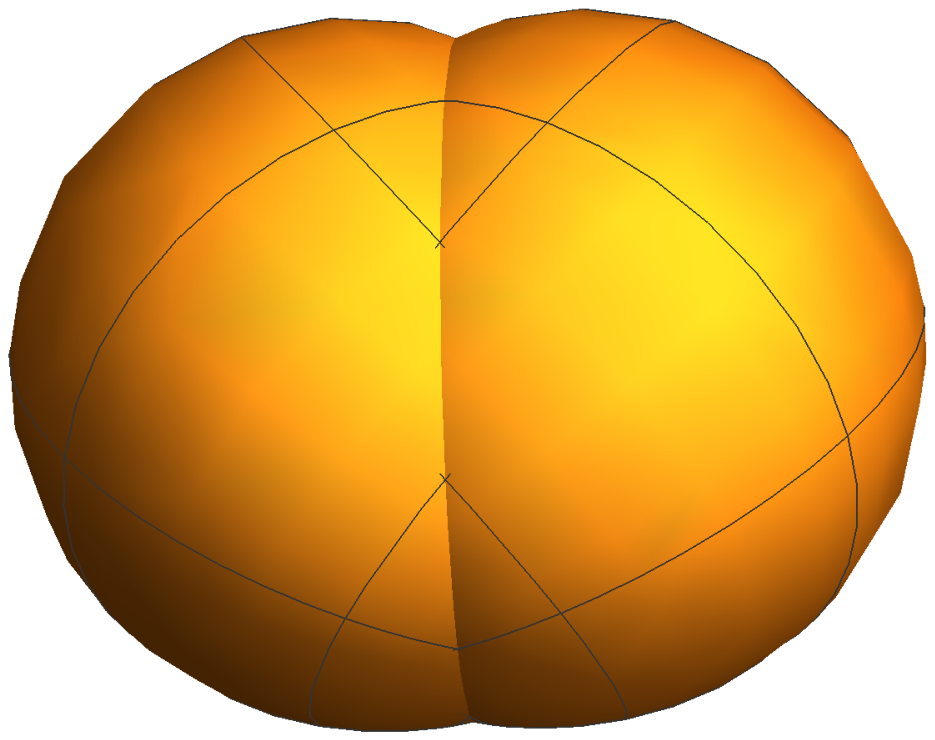}}\hspace{1.0 cm}
\subfigure{\includegraphics[scale=0.30]{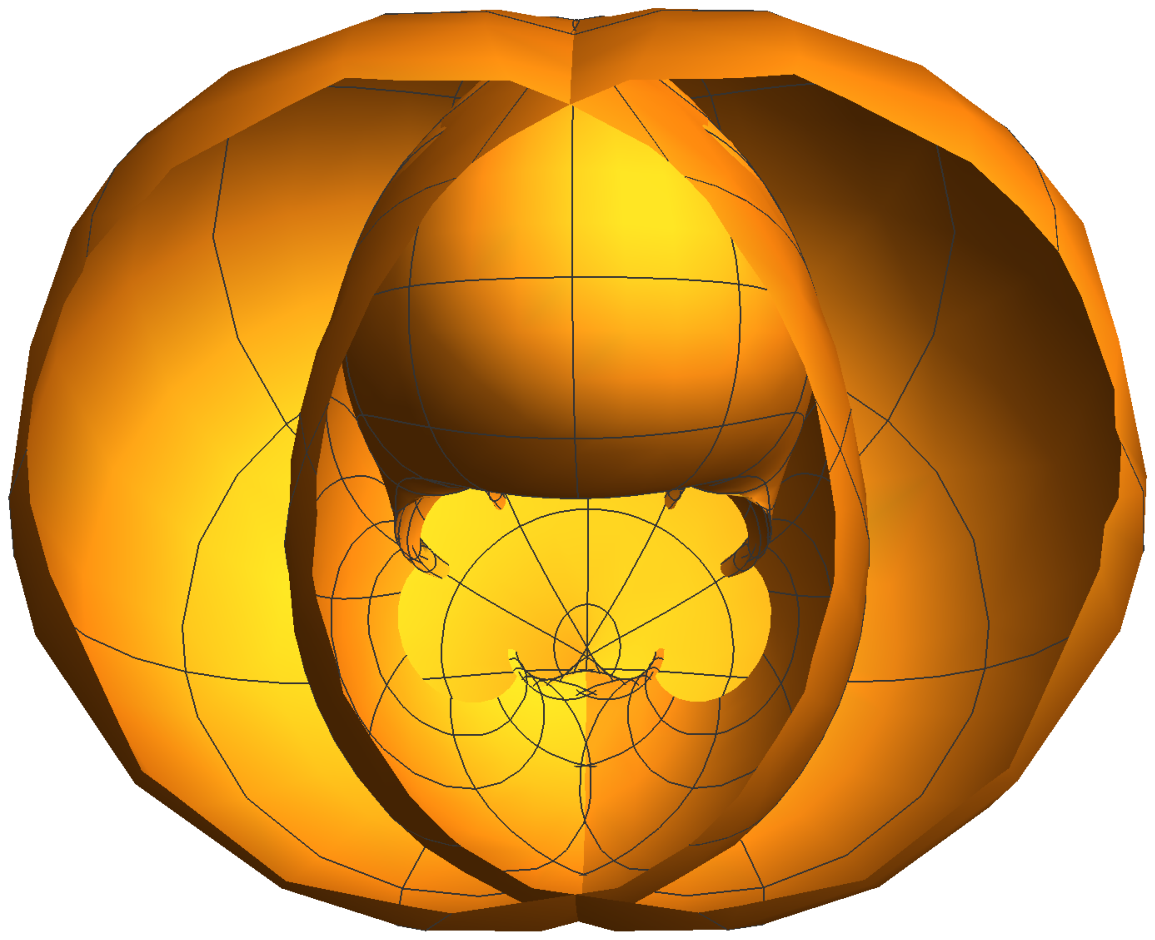}}
\caption{$H_1$-surface for $g(z) = z$ and $A(z) = \sin z$}
\end{figure}
	
\begin{figure}[h]
\centering
\subfigure{\includegraphics[scale=0.30]{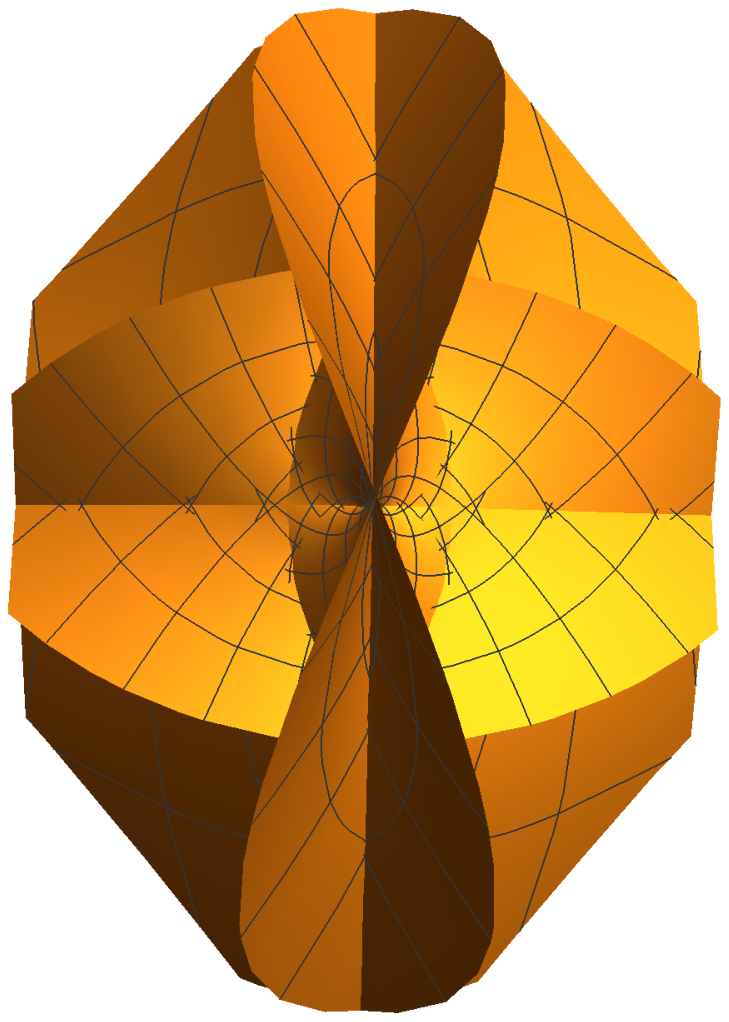}}\hspace{1.0 cm}
\subfigure{\includegraphics[scale=0.30]{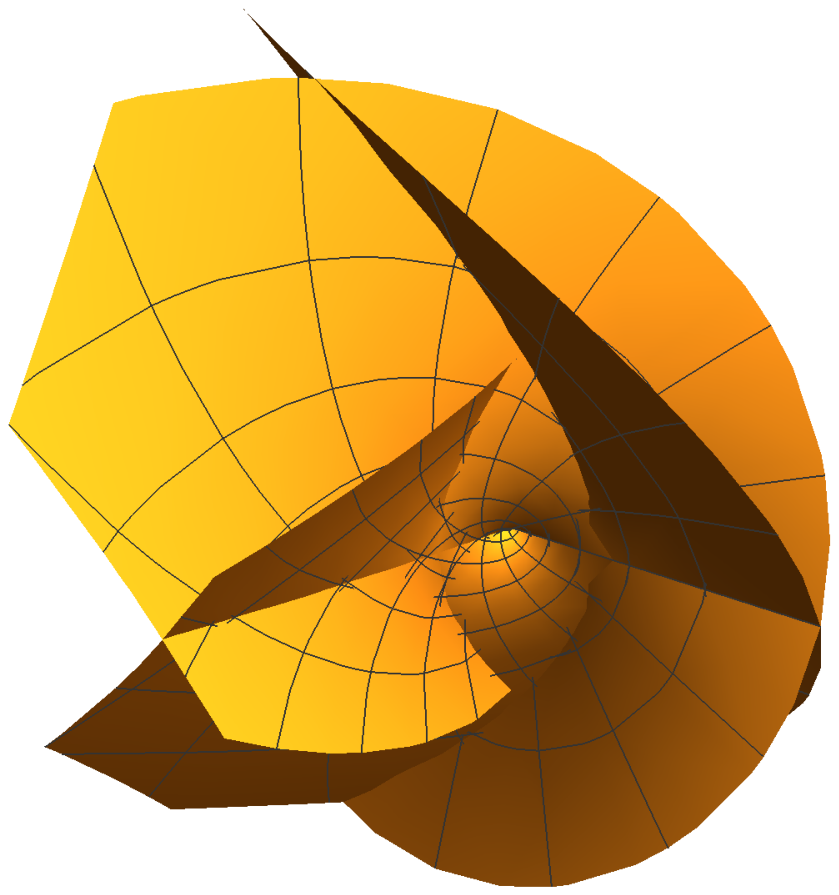}}
\caption{ Minimal surface for $g(z) = z$ and $A(z) = \sin z$ }
\end{figure}
	
\newpage	
\end{example}
\subsection{Special Examples of Minimal Surfaces}

We can construct interesting examples of minimal surfaces by looking at the function $h$ in Corollary \eqref{H1_h} in different ways.

\vspace{0.5cm}
\textbf{Remark 5} The function $h$ in Corollary \eqref{H1_h} can be expressed only in terms of functions $A$ and $g$  as
\vspace{0.5cm}
$$
h(u,v) = \left< 1 \hspace{0.1cm},\hspace{0.1cm} A(z) \right> - 2\frac{\left< g(z) \hspace{0.1cm}, \hspace{0.1cm} \int A(z)g'(z) dz \right>}{1+|g(z)|^2}, \hspace{0.5cm}z=u+iv,
$$

\vspace{0.3cm}
\hspace{-0.4cm}
by integrating by parts  the function $B(z) = \int \Big[ \big(A(z)/g(z)\big)'g(z)^2 + ic_1g'(z) \Big] dz$.

\vspace{0.3cm}
This remark allow us to take the function $h$ assuming a special form. 

\begin{proposition}\label{prop f}
For the function $h$ in Theorem \eqref{teo g} given as

$$
h(u,v) = \left< 1, f'(g(z)) \right> - 2\frac{\left< g(z), f(g(z)) \right>}{1+|g(z)|^2}, \hspace{0.5cm}z=u+iv,
$$

with $f$ holomorphic, the surface $\Sigma$ is a $H_1$-surface.

\end{proposition}

\begin{proof}
Just use Remark $5$ for $A(z) = f'(g(z))$.	
	
\end{proof}
\begin{example}
Considering $f(z) = z^a$, $ a \in \mathbb{R}$, and $g(z) = e^z$  in Proposition  \eqref{prop f}, we have 
\vspace{0.5cm}
$$
h(u,v) = \frac{e^{u(a-1)}(a - 2e^{2u} + ae^{2u}) \cos(v - av)}{1 + e^{2u}}, \hspace{0.5cm} z = u + iv
$$
and the correspondent $\eta$-minimal surfaces are periodic with respect to the second variable. Below it follows some examples for different values of $a$.

\vspace{0.5cm}
\begin{enumerate}
\item For $ a = 3/2 $, the $\eta$-minimal surface has the second variable  periodic with period $2$.

\begin{figure}[h]
\centering
\subfigure{\includegraphics[scale=0.30]{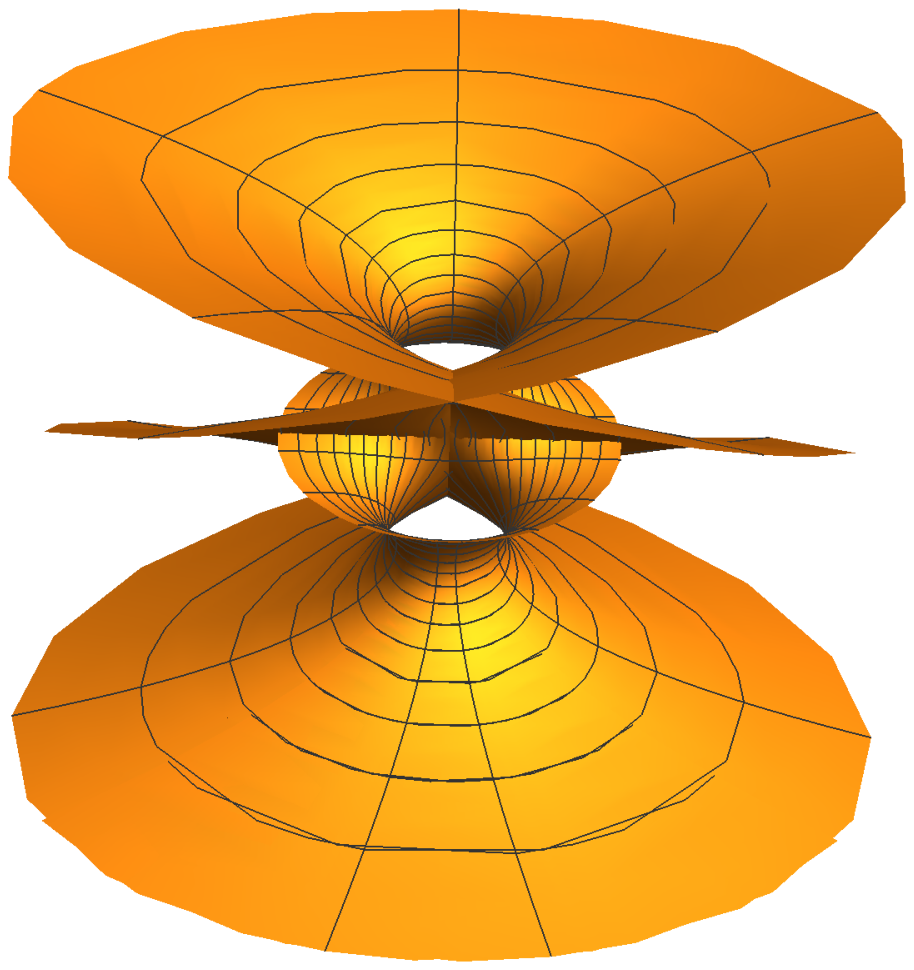}}\hspace{1.0 cm}
\subfigure{\includegraphics[scale=0.30]{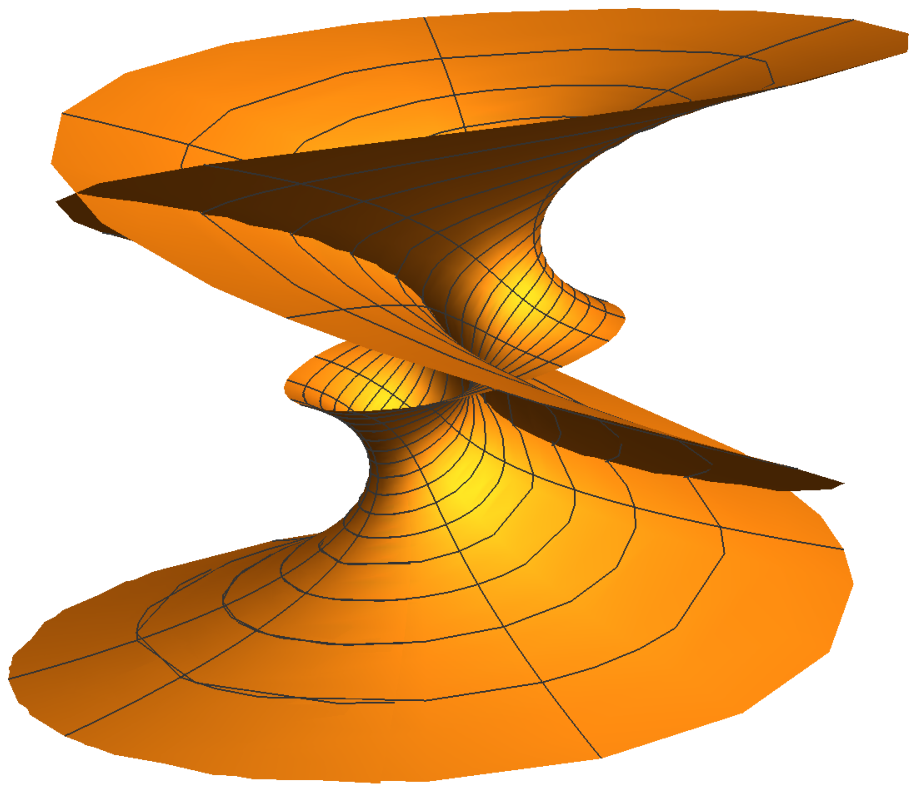}}\hspace{0.5 cm}
\subfigure{\includegraphics[scale=0.30]{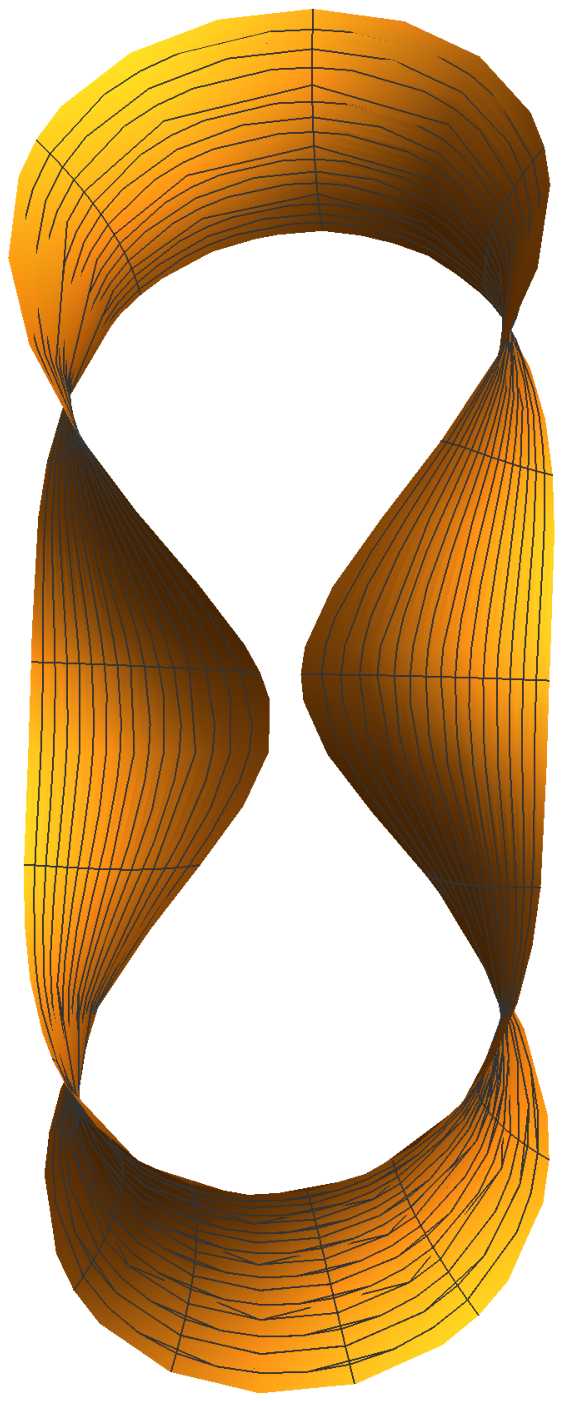}}
\caption{ Minimal surface for $g(z) = e^z$ and $f(z) = z^{3/2}$}
\end{figure}
\newpage
\item For $ a = 4/3 $, the $\eta$-minimal has the second variable  periodic with period $3$.

\begin{figure}[h]
\centering
\subfigure{\includegraphics[scale=0.30]{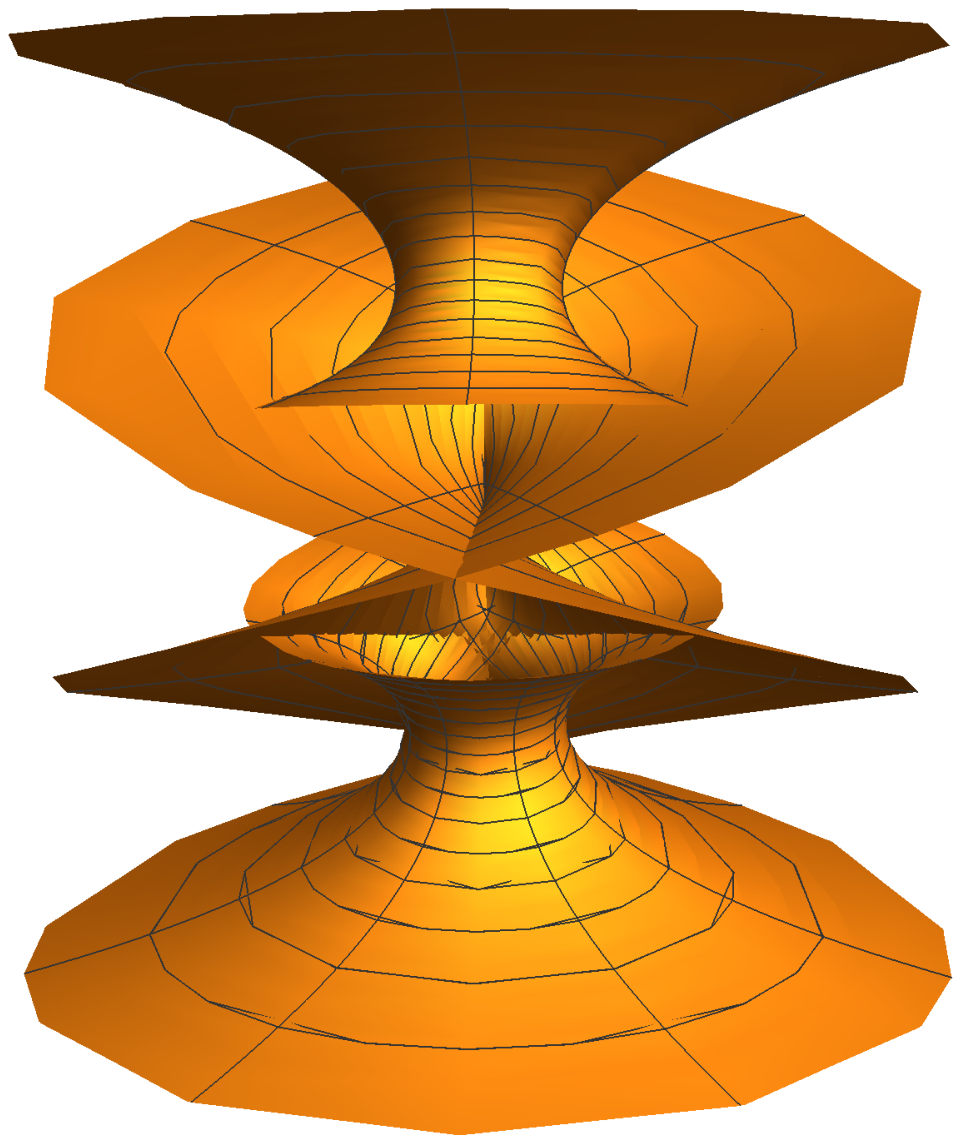}}\hspace{1.0 cm}
\subfigure{\includegraphics[scale=0.30]{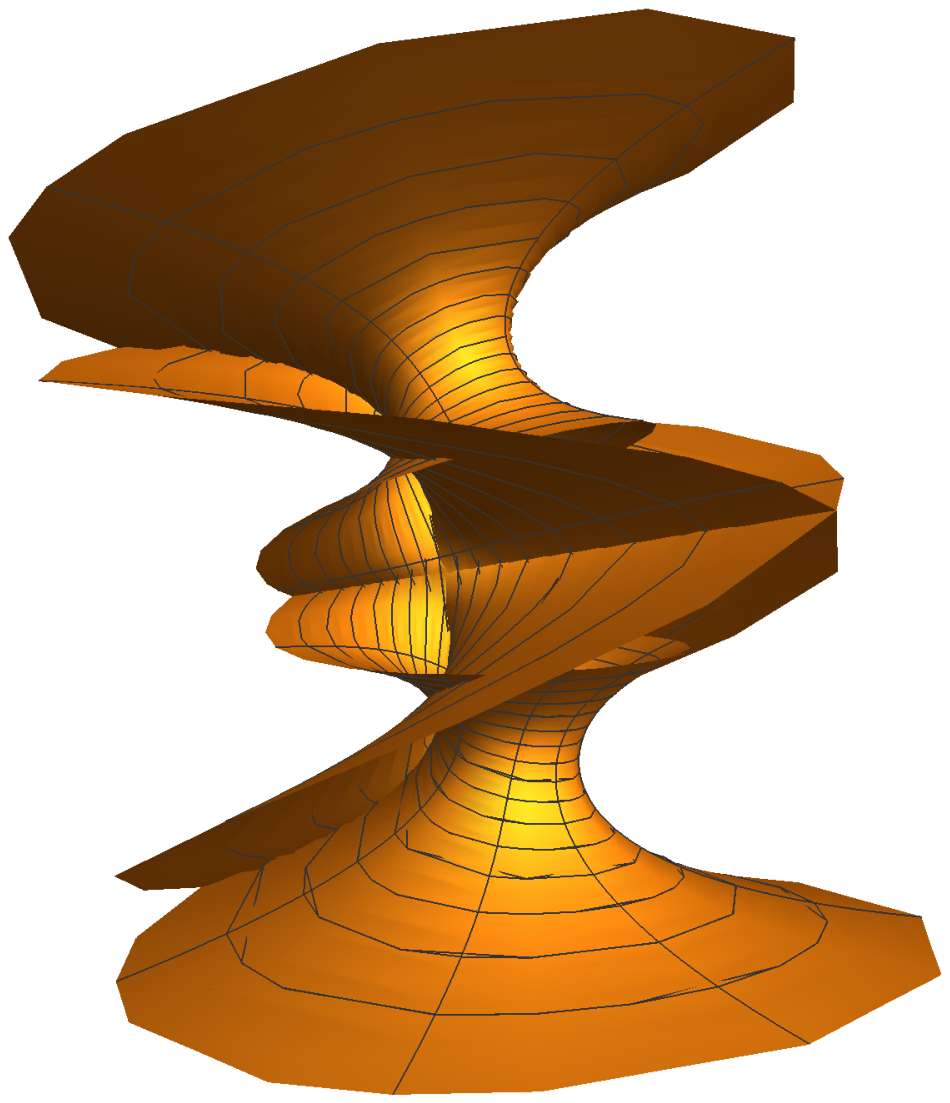}}\hspace{0.5 cm}
\subfigure{\includegraphics[scale=0.30]{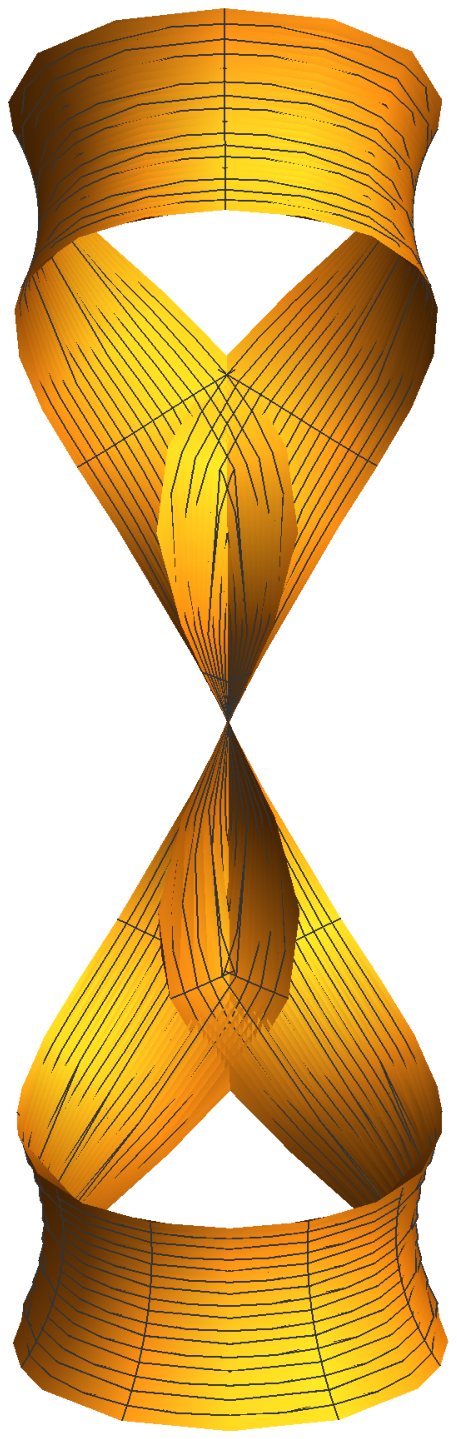}}
\caption{ Minimal surface for $g(z) = e^z$ and $f(z) = z^{4/3}$}
\end{figure}

\item For $ a = 5/4 $, the $\eta$-minimal surface has the second variable  periodic with period $4$.

\begin{figure}[h]
\centering
\subfigure{\includegraphics[scale=0.40]{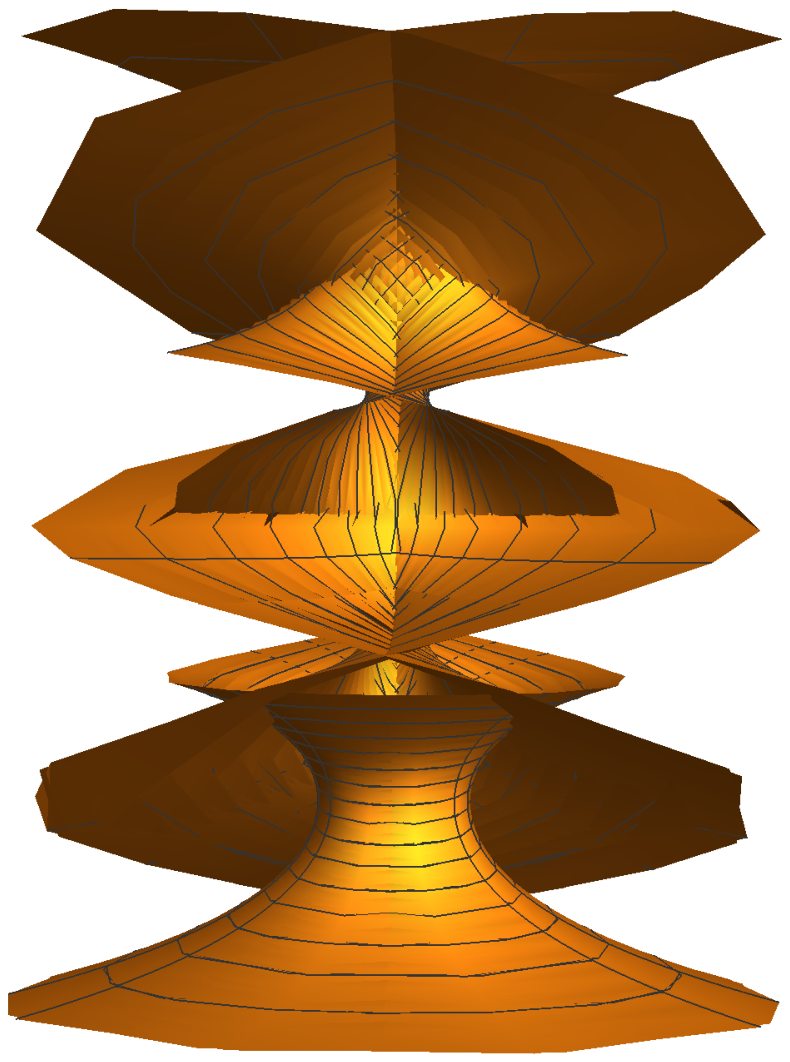}}\hspace{1.0 cm}
\subfigure{\includegraphics[scale=0.40]{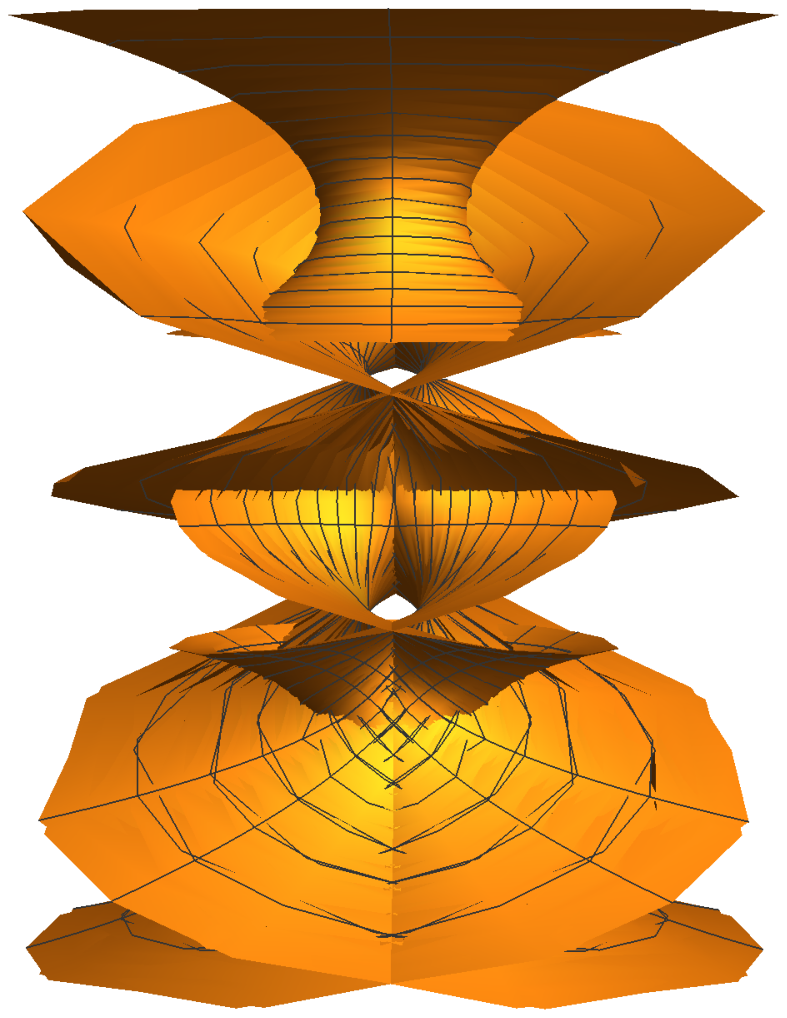}}\hspace{1.0 cm}
\subfigure{\includegraphics[scale=0.40]{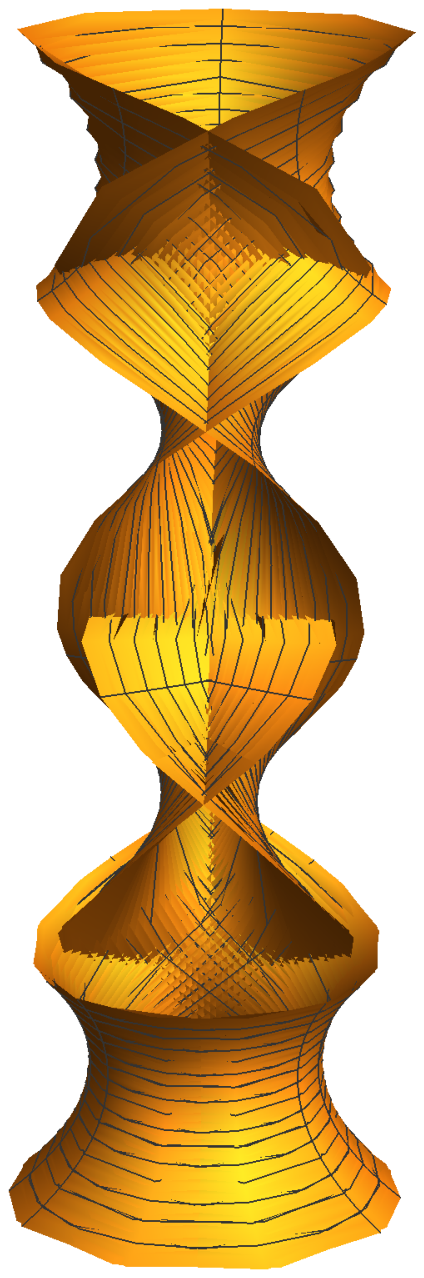}}
\caption{ Minimal surface for $g(z) = e^z$ and $f(z) = z^{5/4}$}
\end{figure} 

\newpage
\item For $ a = 5/3 $, the $\eta$-minimal surface is sketched below.
\begin{figure}[h]
\centering
\subfigure{\includegraphics[scale=0.35]{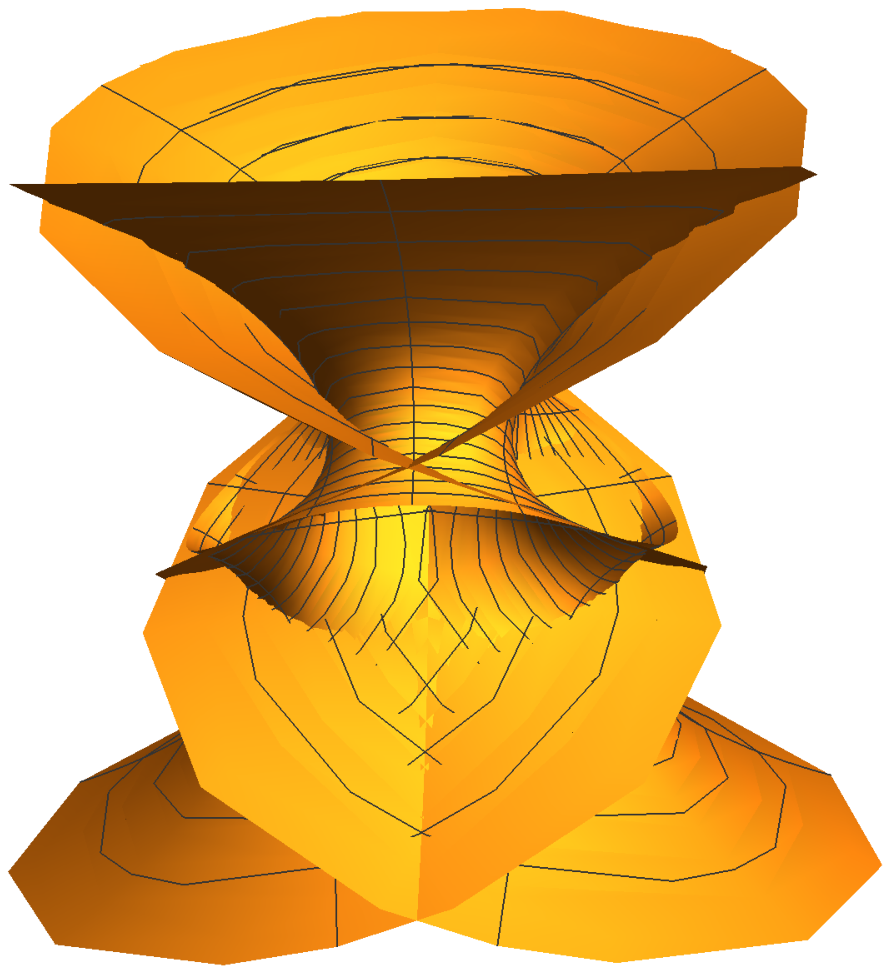}}\hspace{1.0cm}
\subfigure{\includegraphics[scale=0.35]{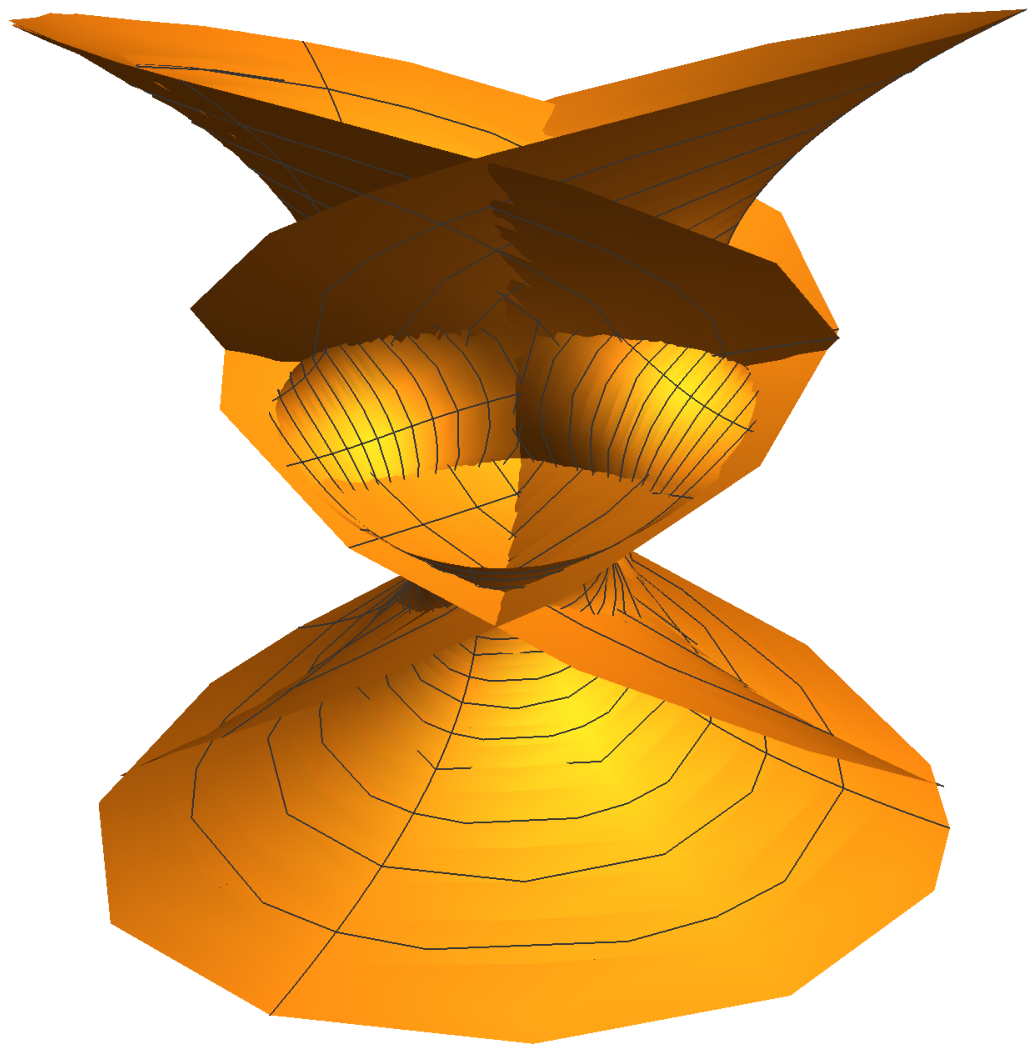}}\hspace{2.0cm}
\subfigure{\includegraphics[scale=0.35]{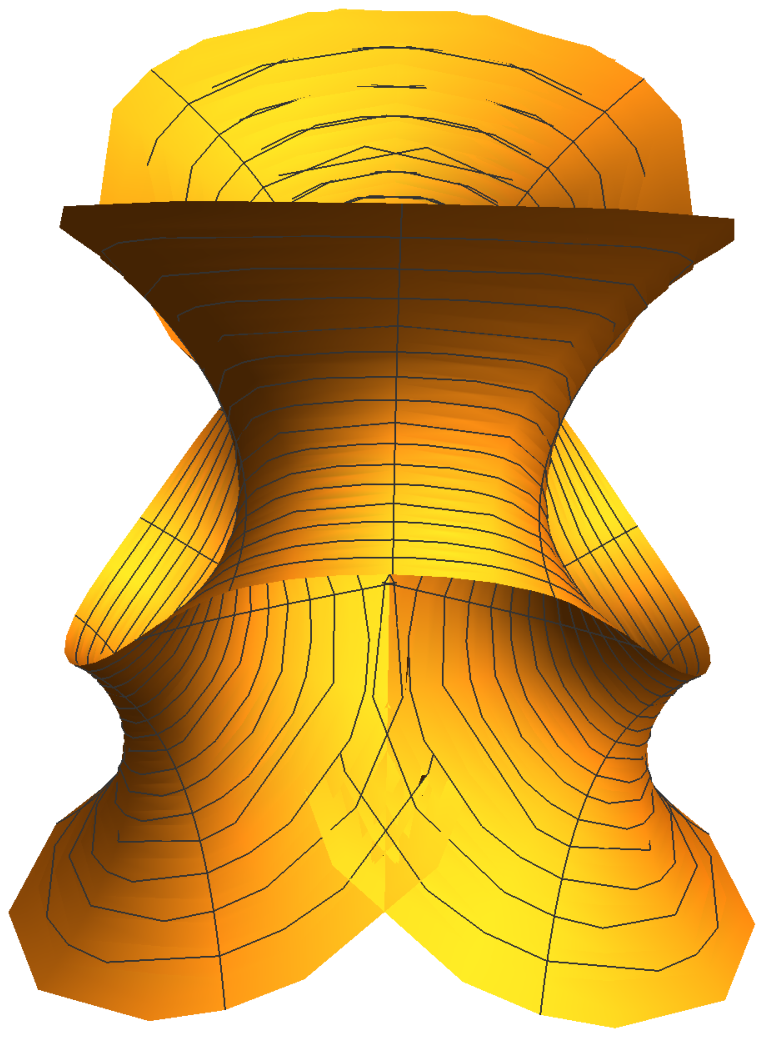}}\hspace{1.0cm}
\subfigure{\includegraphics[scale=0.35]{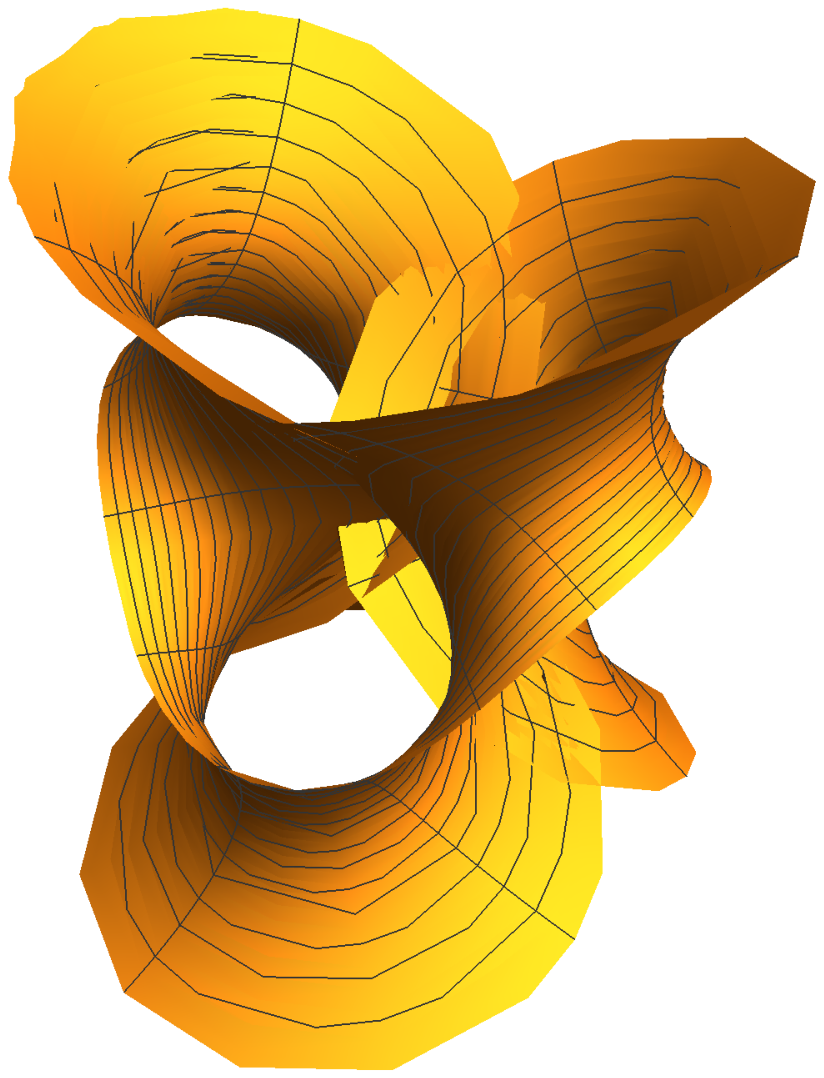}}
\caption{ Minimal surface for $g(z) = e^z$ and $f(z) = z^{5/3}$}
\end{figure}  
\item For $ a = 7/4 $, the $\eta$-minimal surface is like follows.

\begin{figure}[h]
\centering
\subfigure{\includegraphics[scale=0.45]{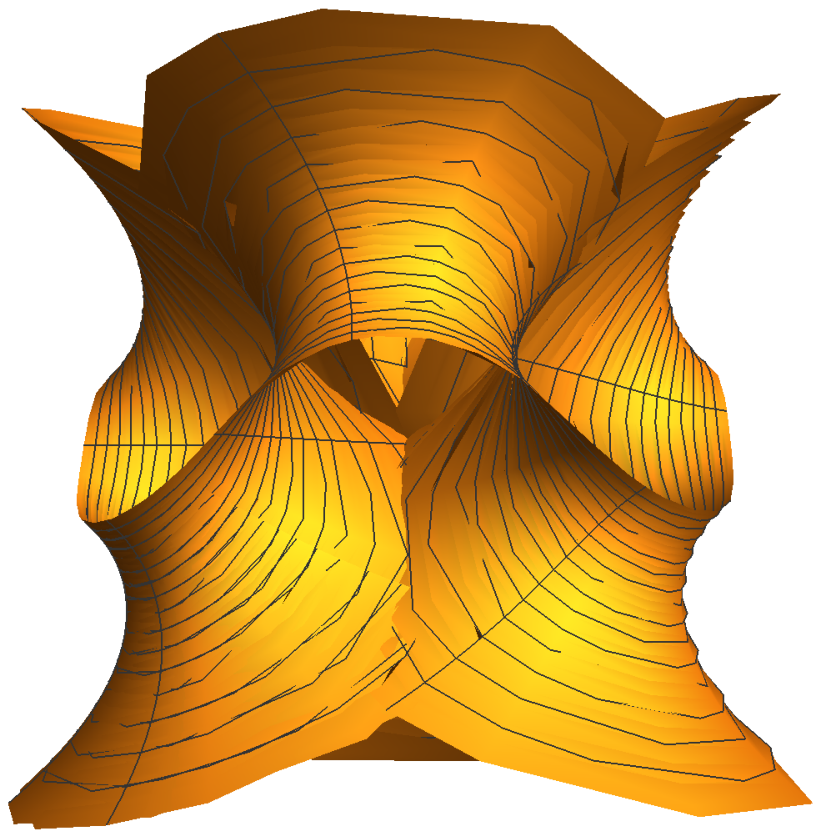}}\hspace{0.5cm}
\subfigure{\includegraphics[scale=0.35]{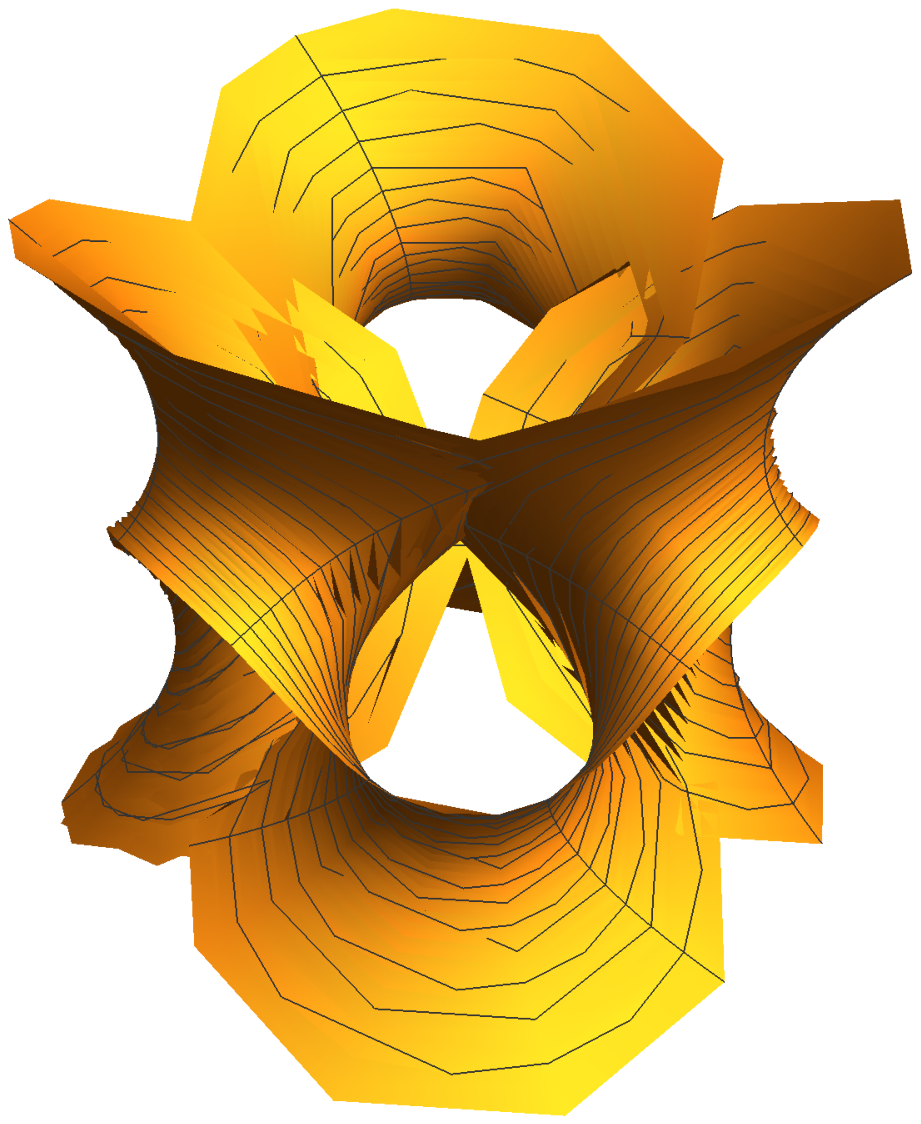}}\hspace{0.5cm}
\subfigure{\includegraphics[scale=0.45]{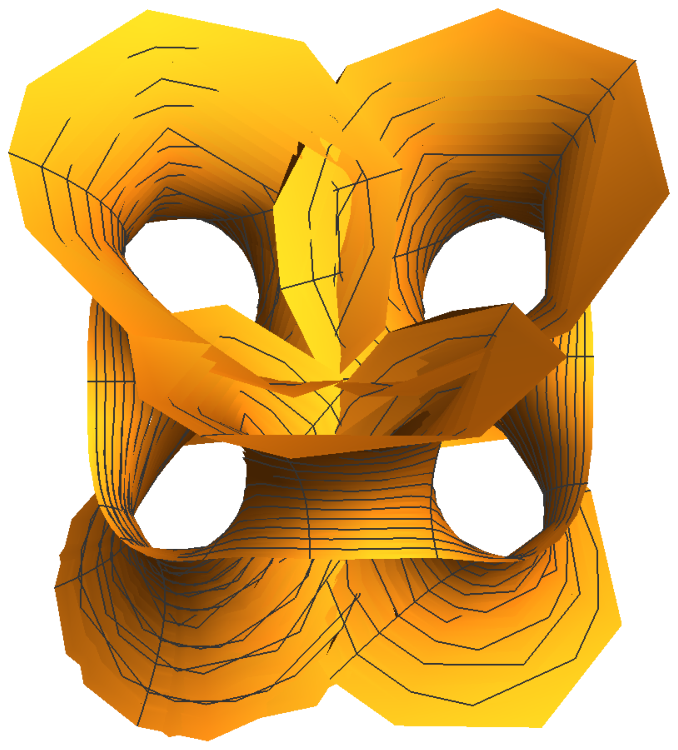}}\hspace{0.5cm}
\caption{ Minimal surface for $g(z) = e^z$ and $f(z) = z^{7/4}$}
\end{figure}  
 
\end{enumerate}
	
\newpage	
\end{example}

\subsection{Examples of $H_2$-Surfaces and Laguerre Minimal Surfaces}

\begin{example}
Considering $g(z) = \sin z$, $A(z) = z$ and $B(z) = z$ in Corollary  \eqref{H2 h}, we obtain the correspondent $H_2$ and $\eta$-Laguerre minimal surfaces in the figure below. 
	
\begin{figure}[h]
\centering
\subfigure{\includegraphics[scale=0.40]{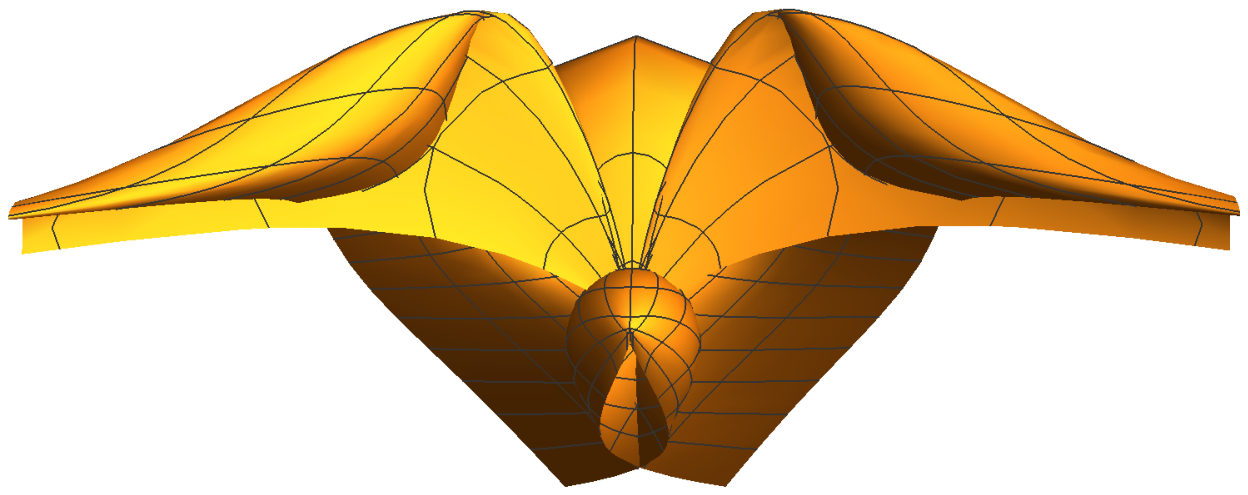}}\hspace{-1.0 cm}
\subfigure{\includegraphics[scale=0.40]{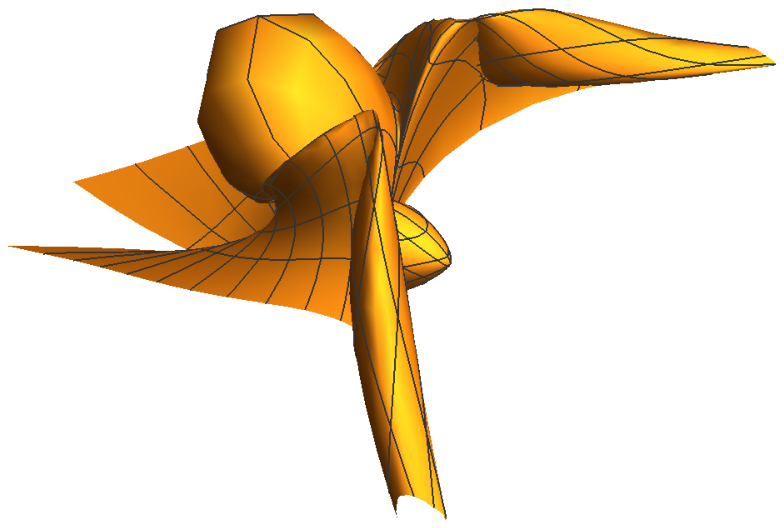}}
\caption{$H_2$-surface for $g(z) = \sin z$, $A(z) = z$ and $B(z) = z$ }
\end{figure}
	
\begin{figure}[h]
\centering
\subfigure{\includegraphics[scale=0.40]{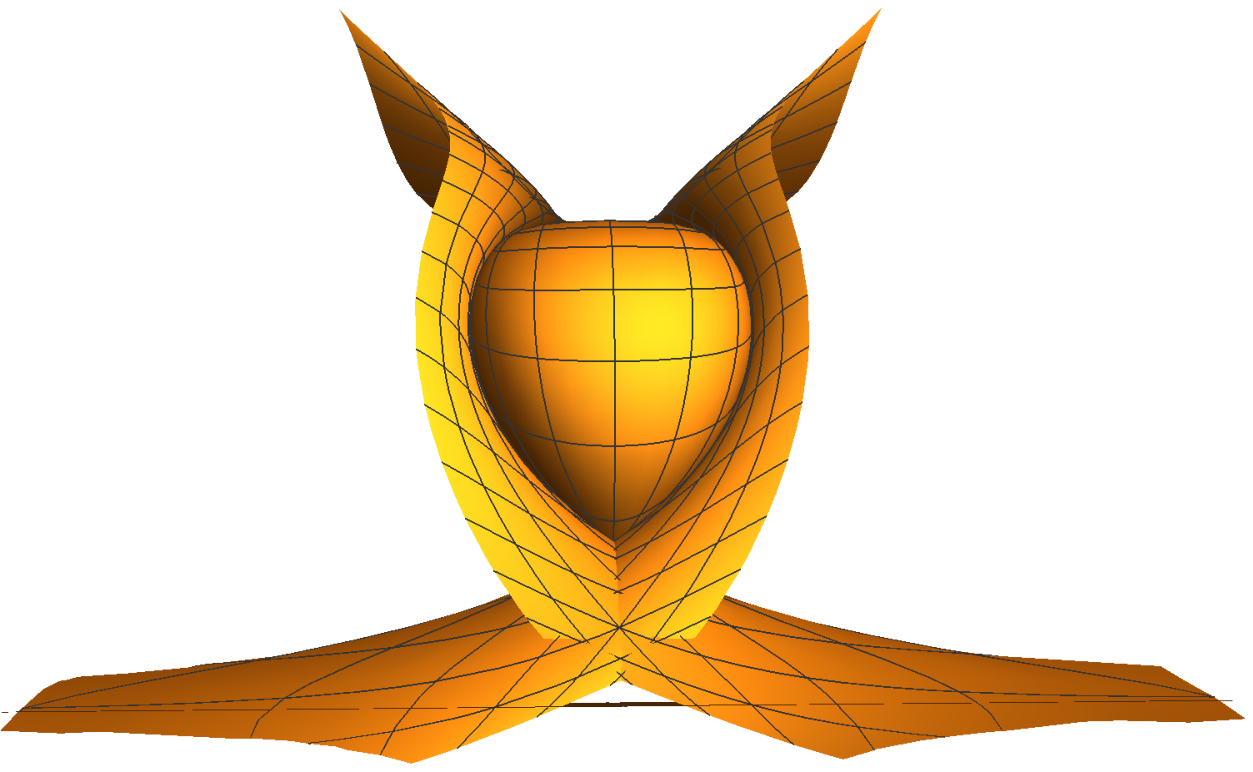}}\hspace{1.0 cm}
\subfigure{\includegraphics[scale=0.40]{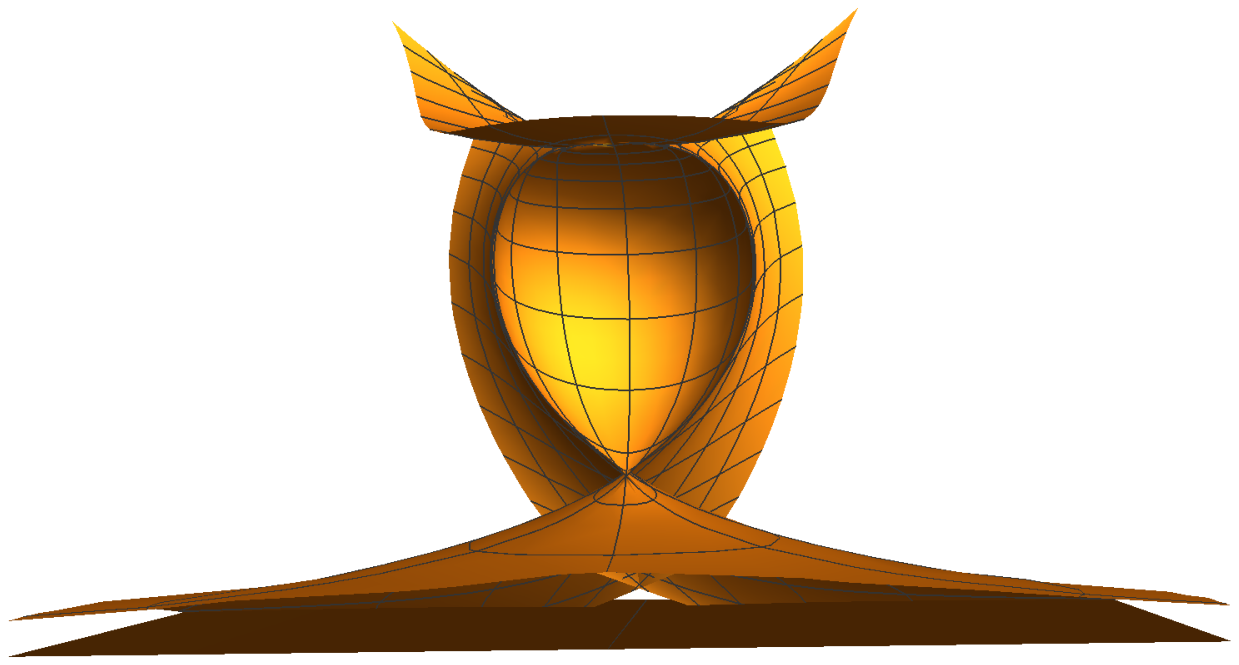}}
\caption{ Laguerre minimal surface for $g(z) = \sin z$, $A(z) = z$ and $B(z) = z$  }
\end{figure}
	
\newpage	
\end{example}

\begin{example}
Considering $g(z) = \sinh z$, $A(z) = \cosh z$ and $B(z) = z^2$ in Corollary  \eqref{H2 h}, we get the correspondent $H_2$ and $\eta$-Laguerre minimal surfaces drawn here. 
	
\begin{figure}[h]
\centering
\subfigure{\includegraphics[scale=0.40]{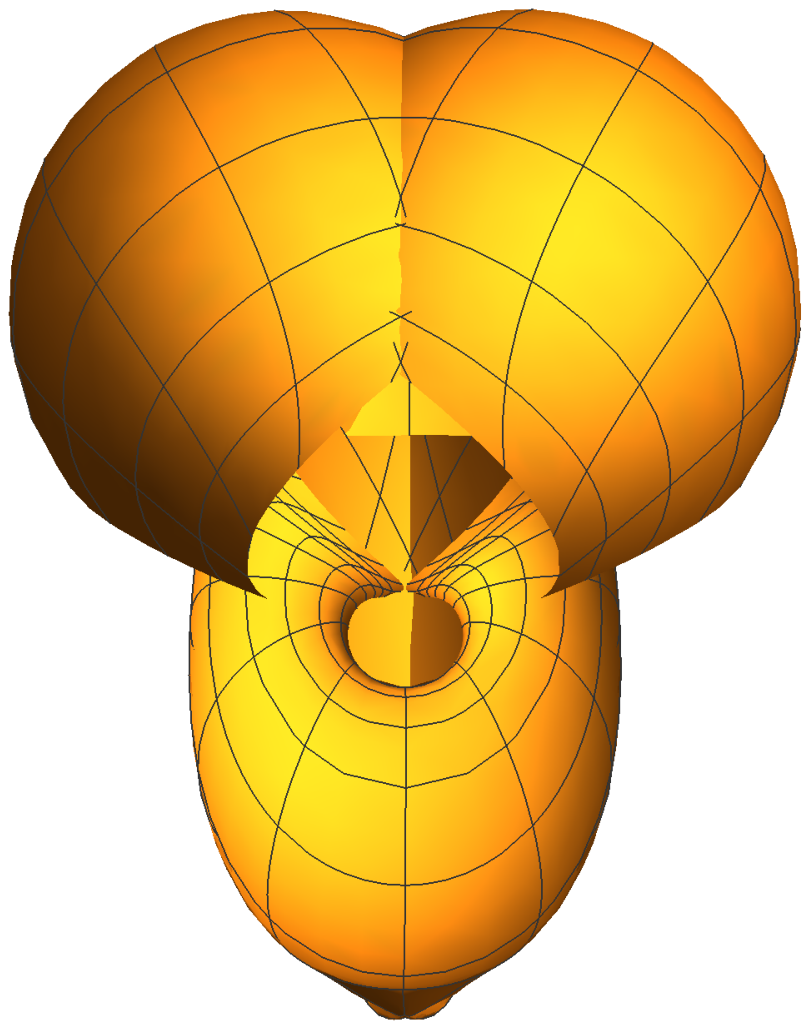}}\hspace{1.0cm}
\subfigure{\includegraphics[scale=0.40]{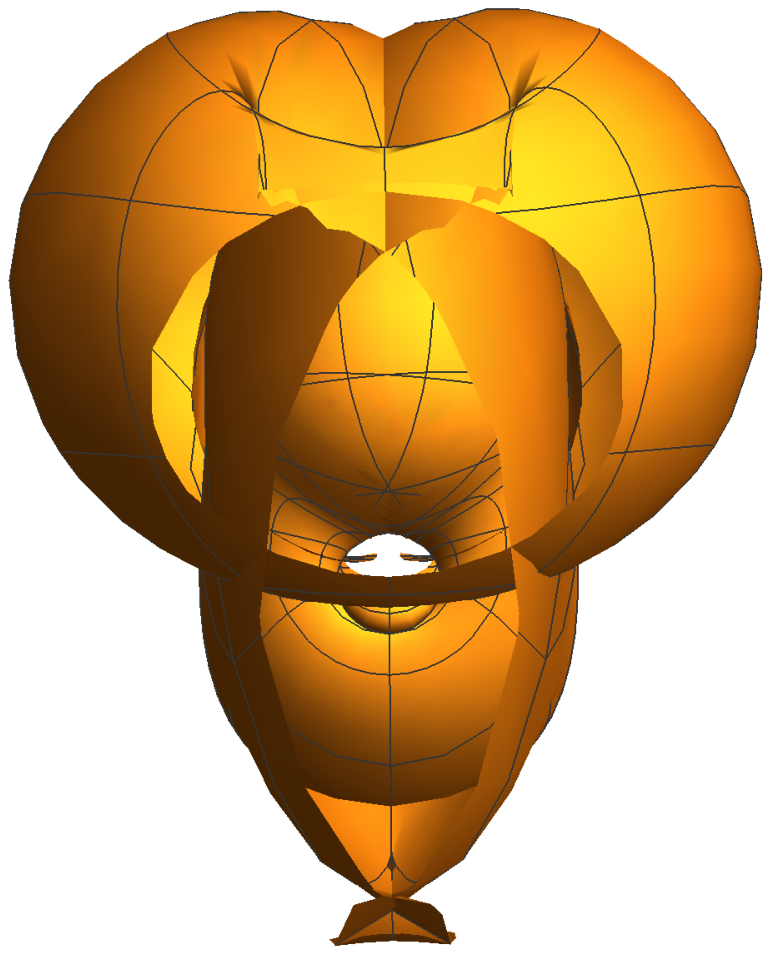}}
\caption{$H_2$-surface for $g(z) = \sinh z$, $A(z) = \cosh z$ and $B(z) = z^2$}
\end{figure}
	
\begin{figure}[h]
\centering
\subfigure{\includegraphics[scale=0.40]{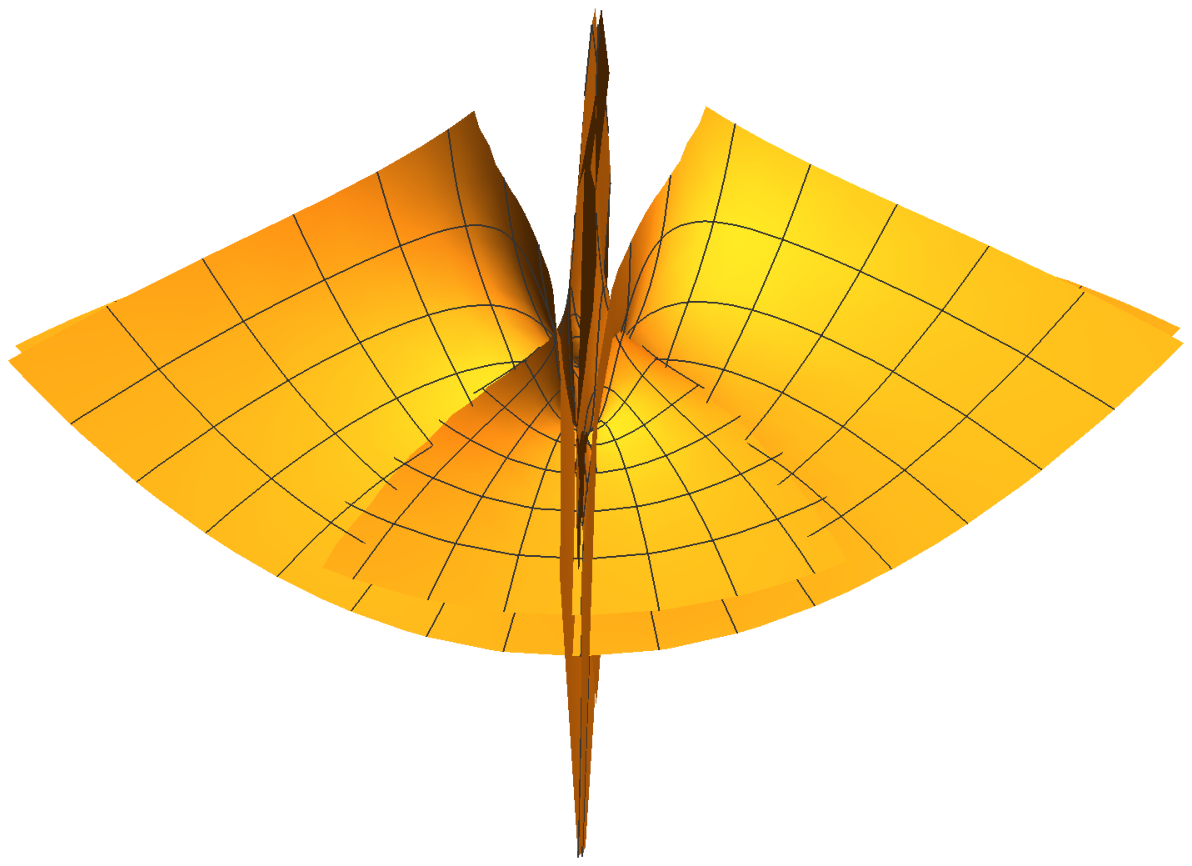}}\hspace{1.0cm}
\subfigure{\includegraphics[scale=0.40]{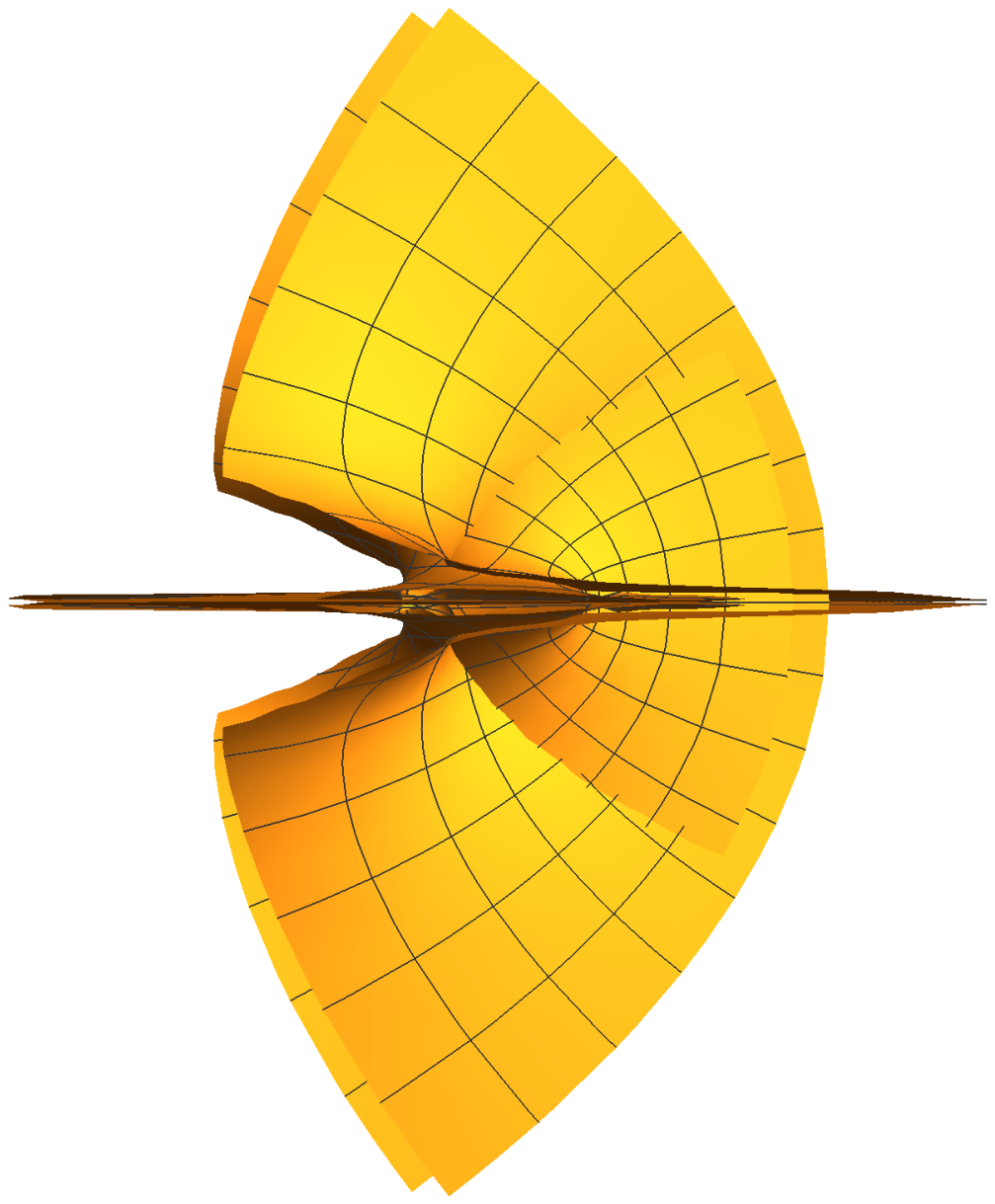}}
\caption{ Laguerre minimal surface for $g(z) = \sinh z$, $A(z) = \cosh z$ and $B(z) = z^2$ }
	\end{figure}
		
\end{example}
\newpage

\begin{example}
Considering $g(z) = z$, $A(z) = e^z$ and $B(z) = \cos z$ in Corollary  \eqref{H2 h}, the correspondent $H_2$ and $\eta$-Laguerre minimal surfaces are like next. 
	
\begin{figure}[h]
\centering
\subfigure{\includegraphics[scale=0.40]{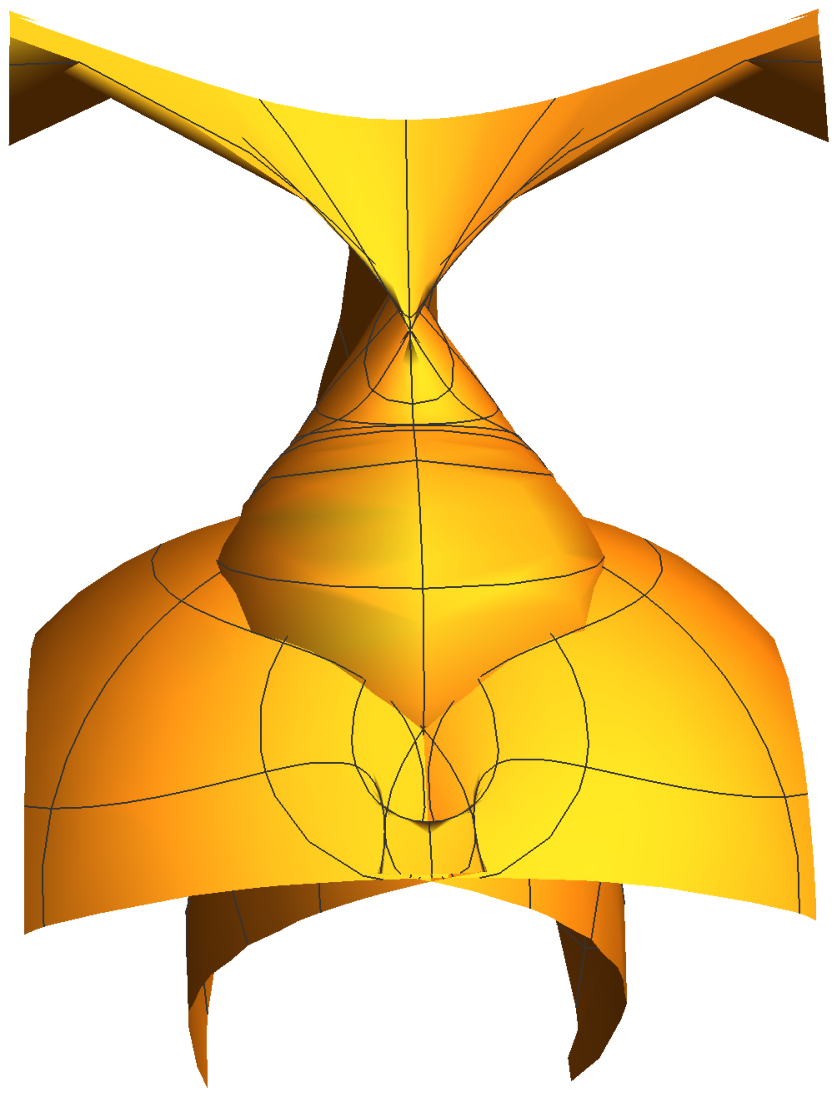}}\hspace{0.5cm}
\subfigure{\includegraphics[scale=0.40]{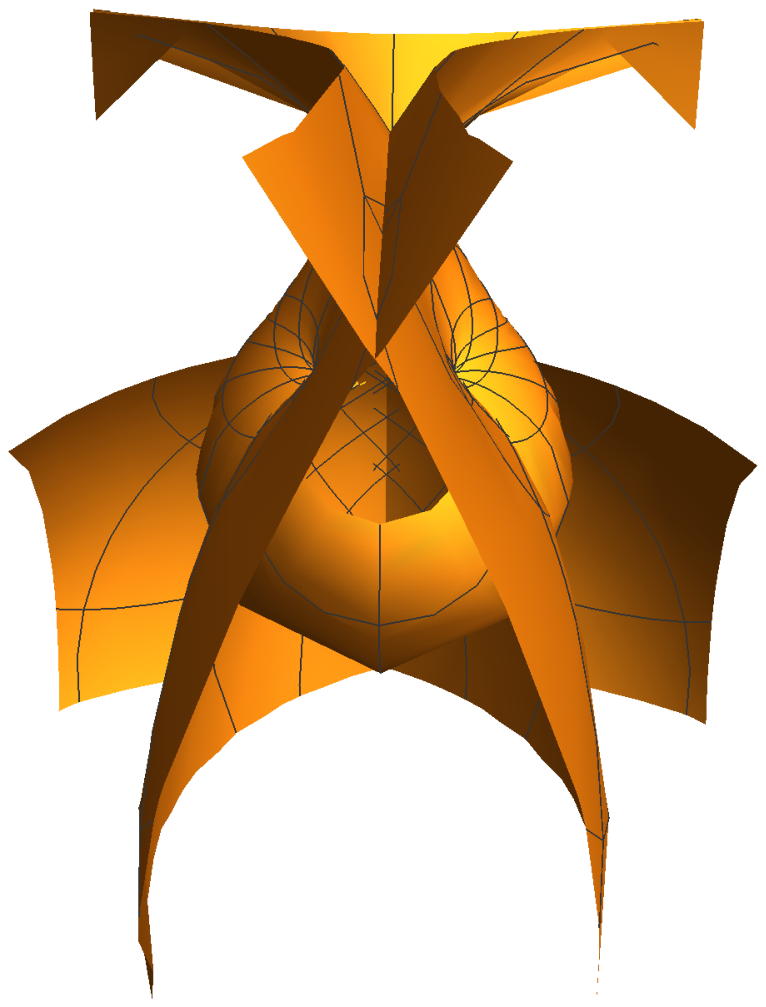}}\hspace{0.5cm}
\subfigure{\includegraphics[scale=0.40]{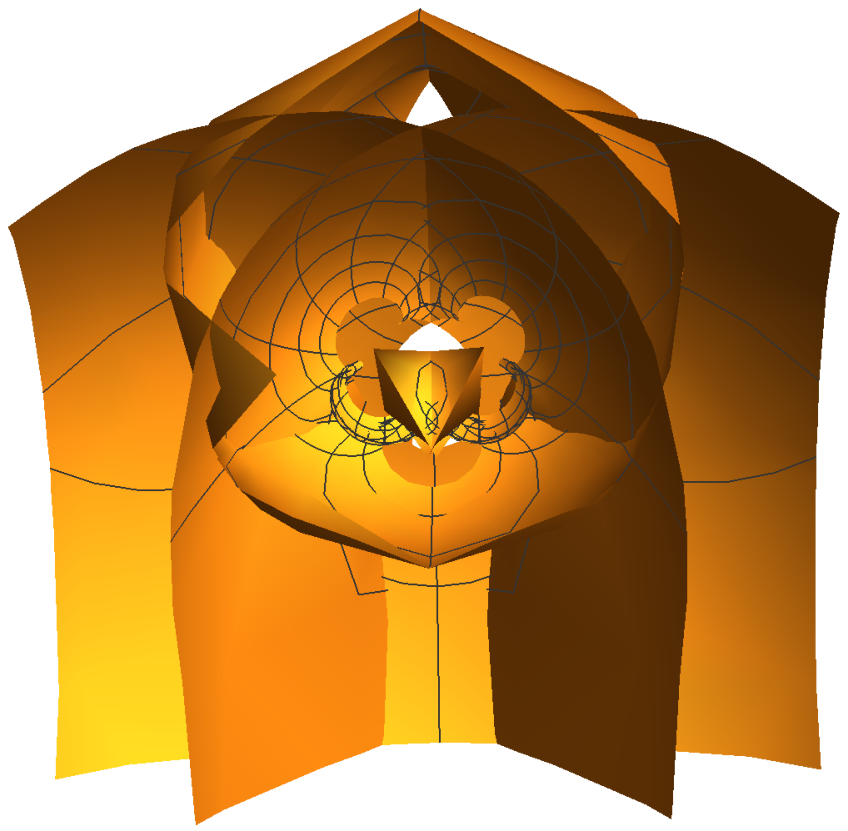}}
\caption{$H_2$-surface for $g(z) = z$, $A(z) = e^z$ and $B(z) = \cos z$}
\end{figure}
	
\begin{figure}[h]
\centering
\subfigure{\includegraphics[scale=0.40]{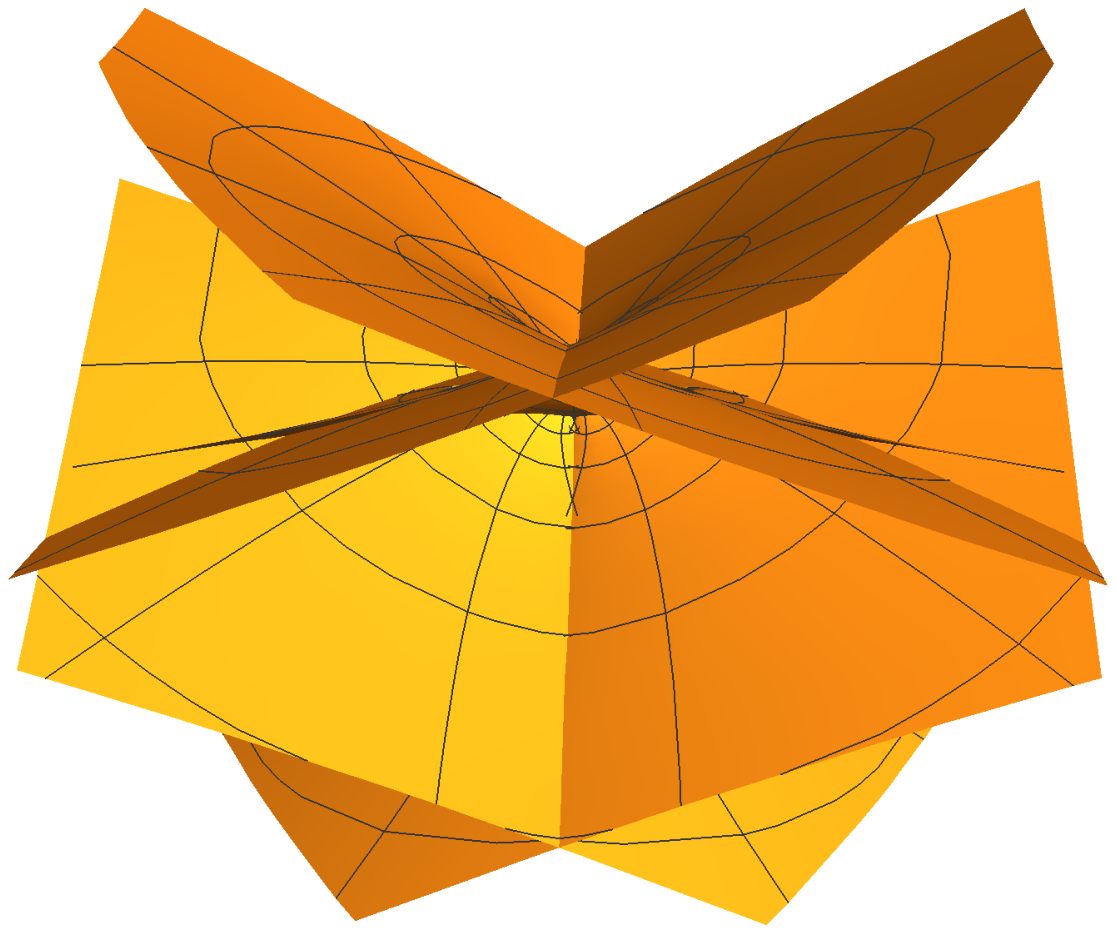}}\hspace{1.0cm}
\subfigure{\includegraphics[scale=0.40]{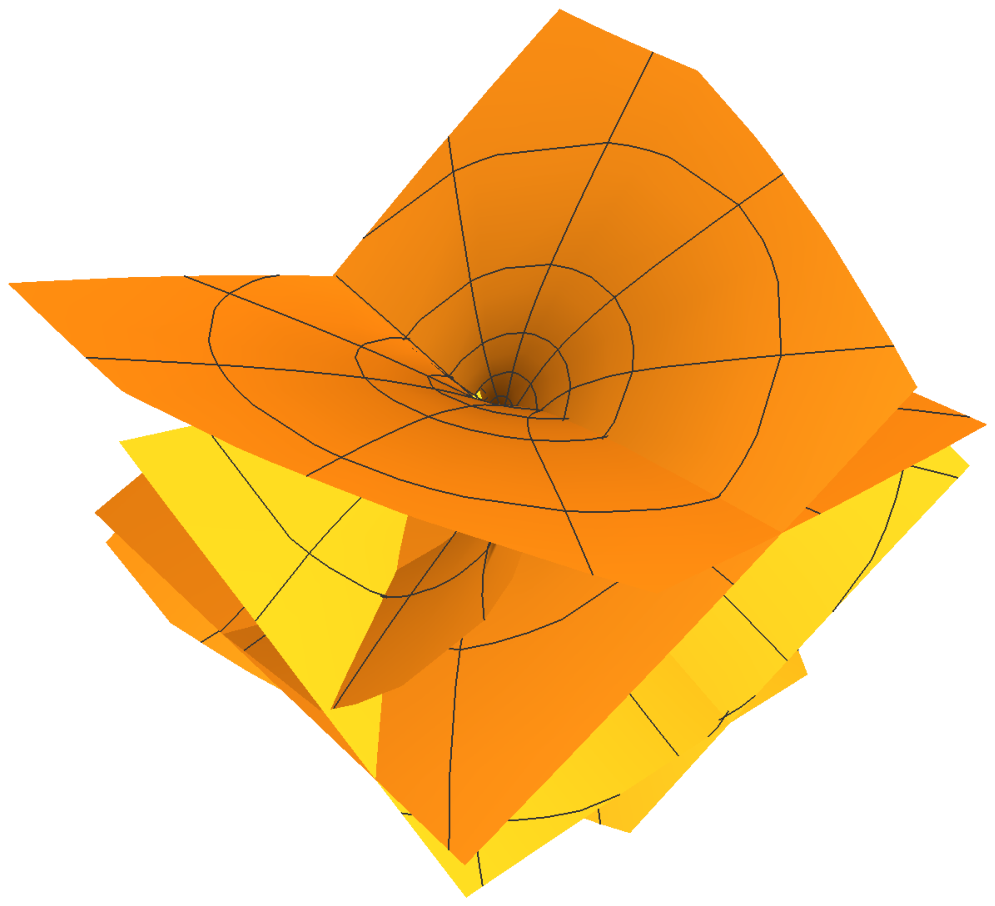}}
\caption{ Laguerre minimal surface $g(z) = z$, $A(z) = e^z$ and $B(z) = \cos z$ }
\end{figure}
	
\newpage	
\end{example}

\vspace{1.0cm}

\section{ $H_1$-surfaces of Rotation}

In the work \cite{n dimensional}, Corro and Riveros define the RSGW-surfaces, which are Weingarten surfaces described by the Helmholtz equation 

\begin{equation}\label{h sol 2}
\Delta h + \frac{8|g'|^2}{T^2}h = 0, \hspace{0.8 cm} T = 1 + |g|^2,
\end{equation}

\vspace{0.3cm}
\hspace{-0.4cm}
for $g$ holomorphic, which solutions they show in Corollary $4$ can be given by

\begin{equation}
 h = \frac{\left <1, A \right > + \left <g, B \right > }{1 + |g|^2}, \hspace{0.5cm} B = \int\big( A'g - Ag' + ic_1g'\big) dz, \hspace{0.5cm}c_1 \in \mathbb{R}
\end{equation}

They show that  the RSGW-surfaces is rotational if and only if $h$ is radial, what implies $h$ as above to be

\begin{equation}\label{h sol 3}
h(u,v) = a_1 -\big(a_2 + a_1(u-1)\big)\tanh u.
\end{equation}

\vspace{0.3cm}
In our case, we know that the $H_1$-surfaces are described by the generalized Helmholtz equation \eqref{eq H1}, which are the same characterization for the  RSGW-surfaces. Similar to the RSGW-surfaces, the $H_1$-surfaces are rotational iff $h$ is radial. Thus, when $h$ is given by \eqref{h sol 3}, the correspondent $H_1$-surface is rotational and we have the next proposition.

\begin{proposition}\label{teo rot 1}
A surface $X$ (respectively $\eta$) as in Theorem \eqref{teo g} is a $H_1$-surface ( respectively minimal surface) of rotation if and only if $g(z) = e^z$ and $ h = a_1 -\big(a_2 + a_1(u-1)\big)\tanh u$. In this case, the surfaces $X$ and $ \eta $  can be locally parameterized by 

$$
X(u,v) = \big(M(u)\cos v\hspace{0.1cm},\hspace{0.1cm} M(u) \sin v\hspace{0.1cm}, \hspace{0.1cm} N(u) \big)
$$
$$
\eta(u,v) = \big(M_1(u)\cos v\hspace{0.1cm},\hspace{0.1cm} M_1(u) \sin v\hspace{0.1cm}, \hspace{0.1cm} N_1(u) \big)
$$
where

\begin{eqnarray*} 
M(u) &= & \frac{2e^u\{4a_2^2e^{2u}+2a_2a_1l_1 + a_1\big[-2c(1+e^{2u})^2 + a_1 l_2\big]\}}{(1+ e^{2u})\big[4a_2^2e^{2u} + 8 a_1a_2e^{2u}(u-1) + a_1^2\big(1+e^{4u} + e^{2u}(6 - 8u + 4u^2)\big)\big]}, \\
N(u) & = & \frac{4a_2e^{2u}p_1  -8a_1e^{2u}p_2  + a_1^2 p_3  }{(1+ e^{2u})\big[4a_2^2e^{2u} + 8 a_2a_1e^{2u}(u-1) + a_1^2\big(1+e^{4u} + e^{2u}(6 - 8u + 4u^2)\big)\big]}, \\
 & & \\
M_1(u) & = & a_1 \left(\frac{1+e^{2u}}{2 e^u}\right), \\
N_1(u) & = & a_2 + a_1(u-1), \hspace{0.5cm} a_1, a_2 \in \mathbb{R}
\end{eqnarray*}
with 

\begin{eqnarray*}
\hspace{-4.0cm}l_1 & = & \big(e^{4u}+4e^{2u}(u-1)-1\big), \\
\hspace{-4.0cm}l_2 & = & \big(1-2u + e^{4u}(2u - 3) + e^{2u}(2-8u+4u^2)\big), \\
\hspace{-4.0cm}p_1 & = &  a_2(e^{2u}-1)-2c(1+e^{2u}), \\
\hspace{-4.0cm}p_2 & = & c(1+e^{2u})(u-1) + a_2(u- e^{2u}(u-2)),\\ 
\hspace{-4.0cm}p_3 & = & 1 - e^{6u} + e^{2u}(5-4u^2) + e^{4u}(11 - 16u + 4u^2),
\end{eqnarray*}

\end{proposition}

\begin{proof}
The proof will take the same steps as in demonstration of Theorem $8$ in \cite{n dimensional}. We take $g(w)= w$, $w \in \mathbb{C}$, and, in this case, the Remark $1$ asserts that $X$ is a $H_1$-surface of rotation if, and only if, $h$ is a radial function, i.e., $h(w) = J(|w|^2)$, $ w \in \mathbb{C}$, for any real differentiable function $J$. Making the change of parameters $ w= e^z, \hspace{0.3cm} z=u+iv \in \mathbb{C}$,	we have $g(z)= e^z$ and $h = J(e^{2u})$, so that $ h_{,2}=0$. From Corollary \eqref{H1_h}, we have

\begin{equation}\label{h rot}
h = \frac{\left< 1, A \right> + \left< e^z, B \right> }{1+e^{2u}}
\end{equation}

\vspace{0.2cm}
Differentiating \eqref{h rot} with respect to $v$, using the fact that $ h_{,2}=0$ and Proposition \eqref{pro cte}, we find

$$
A = i\bar{z_1}e^z + c_1 z + z_2, \hspace{1.0cm} B = (z_3 + c_2 z)e^z + i z_1
$$
where $c_1, c_2 \in \mathbb{R}$, $z_k = a_k +ib_k$, \hspace{0.1cm} $a_k, b_k \in \mathbb{R}$, $k = 1,2,3$.

Using the Corollary \eqref{H1_h} again,

$$
A = c_1 z + z_2, \hspace{1.0 cm} B = (2c_1 - z_2 + ic_1) e^z - c_1 ze^z
$$

\vspace{0.2cm}
Substituing the expressions above for $A$ and $B$ in \eqref{h rot}, we get $ h(z) = a_1 -\big(a_2 + a_1(u-1)\big)\tanh u$, \hspace{0.2cm} $z = u+iv$.
	
Now, in order to get the parameterization $X$, it is enough take $h$ as above and $g(z) = e^z$ in the expression \eqref{eq X g}, noting that $\eta$ is expressed as in \eqref{eta WR}.

\end{proof}

\vspace{0.5cm}
\subsection{Examples of Rotational $H_1$-Surfaces and Rotational $ \eta $ Minimal Surface}

\vspace{0.5cm}
\begin{example}
Being $\eta$ a rotational minimal surface, we recover the catenoid. Below it is sketched for  $a_1 = a_2 = 1$ in Proposition \eqref{teo rot 1}.

\begin{figure}[h]
\centering
\subfigure{\includegraphics[scale=0.40]{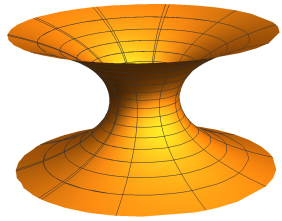}}
\caption{Minimal surface of rotation for $a_1 = a_2 = 1$ }
\end{figure}

\end{example}

\newpage
\begin{example}
Considering $a_1 = a_2 = c = 1$ in Theorem \eqref{teo rot 1}, we obtain the $H_1$-surface of rotation sketched below. This surface have two isolated singularities and two circles of singularities. The profile curve is given on the left.

\begin{figure}[h]
\centering
\subfigure{\includegraphics[scale=0.30]{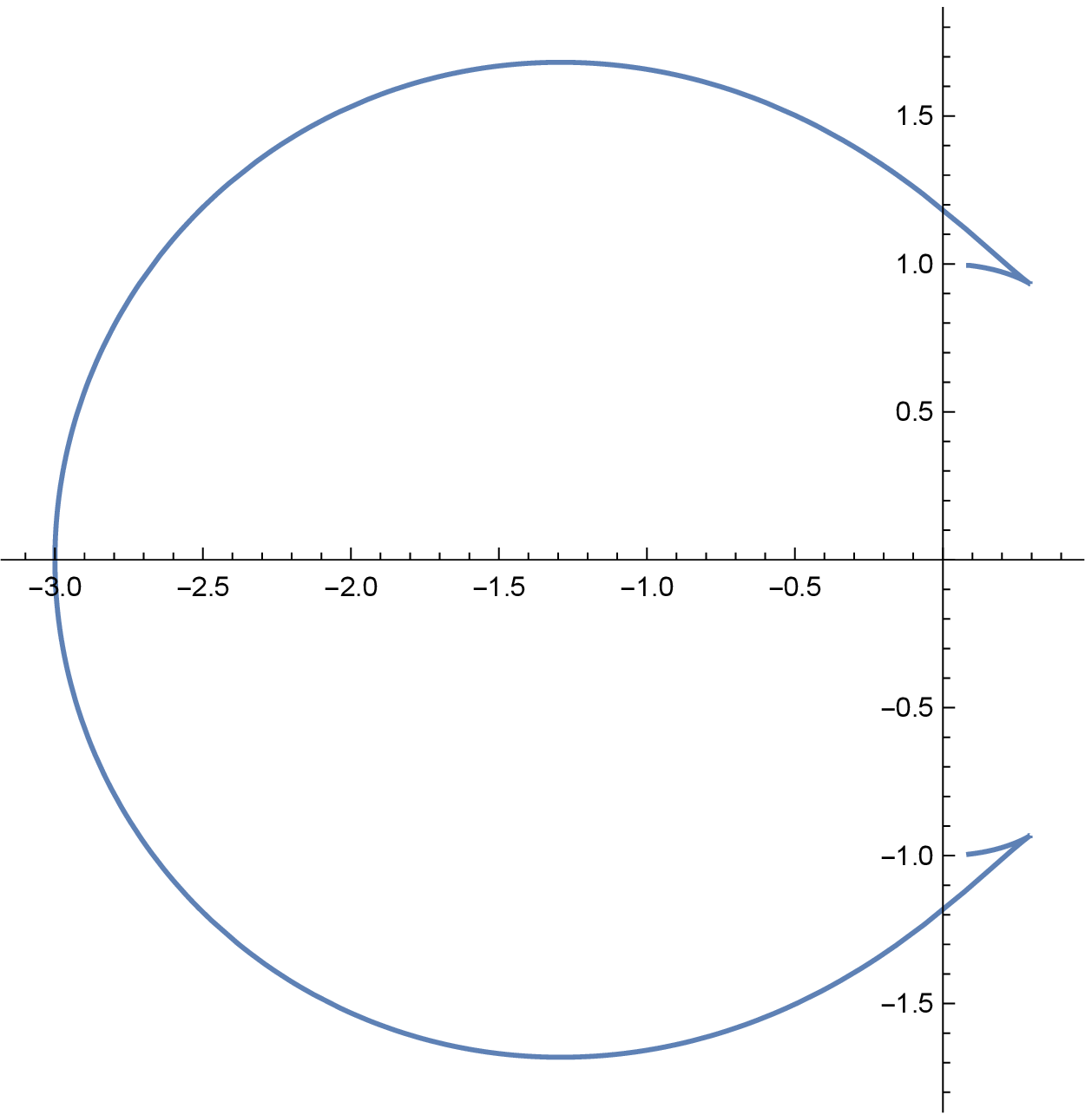}}
\subfigure{\includegraphics[scale=0.45]{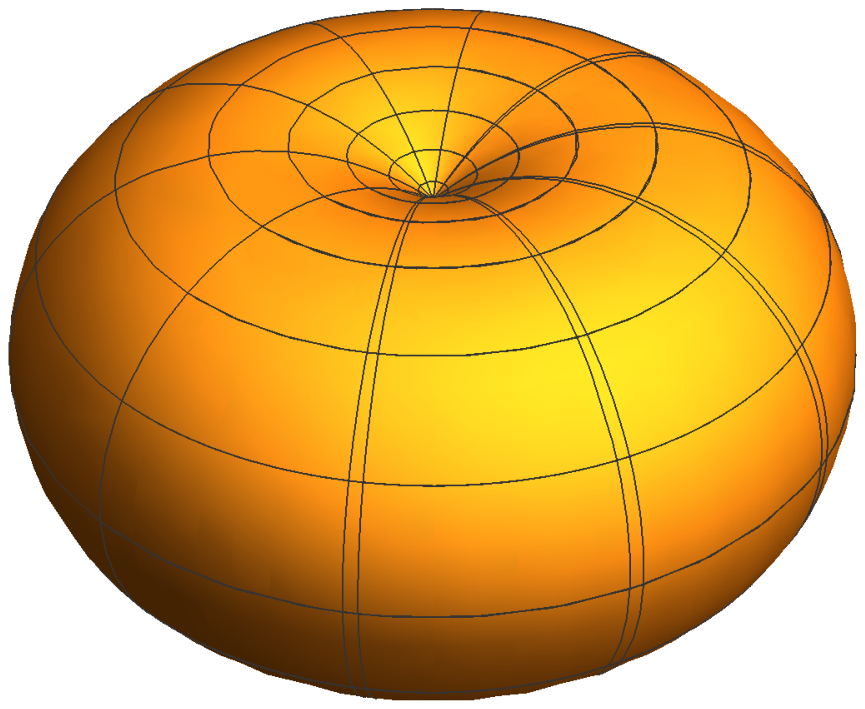}}
\subfigure{\includegraphics[scale=0.45]{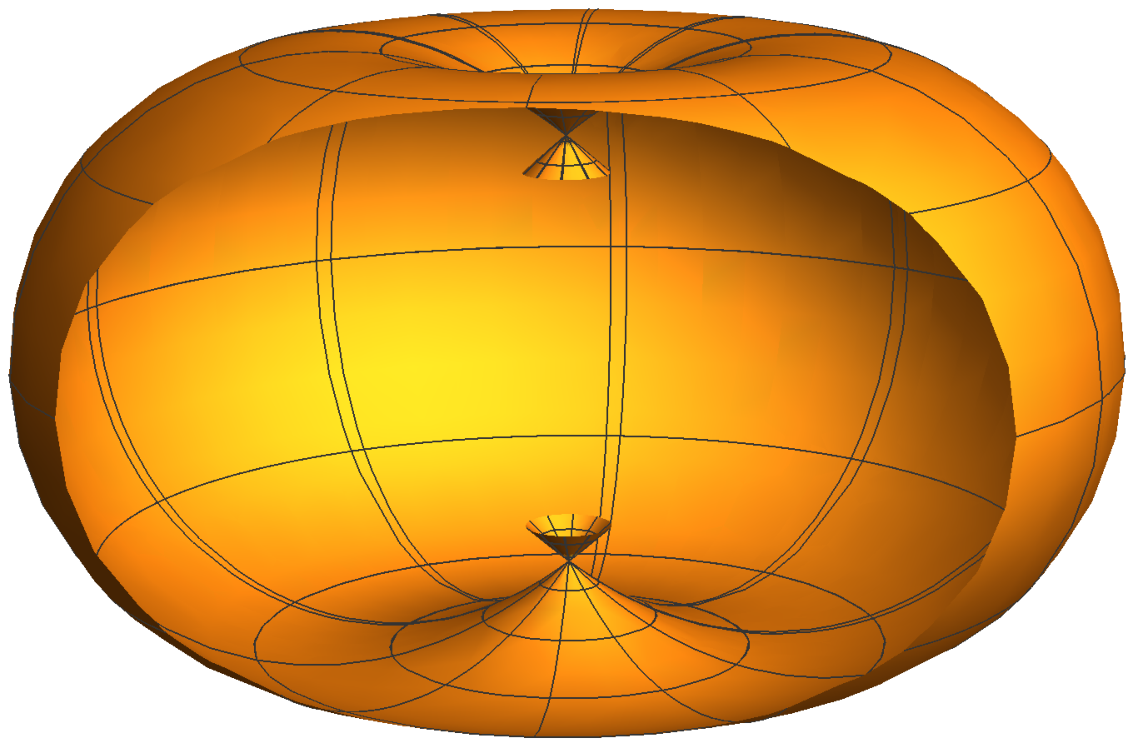}}
\caption{$H_1$-surface of rotation for $a_1 = a_2 = c = 1$. }
\end{figure}

\end{example}

\vspace{0.5cm}
\begin{example}
Considering $a_1 = 1$, $a_2 = 3$ and $c = 1$ in Theorem  \eqref{teo rot 1}, we obtain the $H_1$-surface of rotation sketched below. This surface have one isolated singularity and one circle of singularities. The profile curve is given on the left.
	
\begin{figure}[h]
\centering
\subfigure{\includegraphics[scale=0.30]{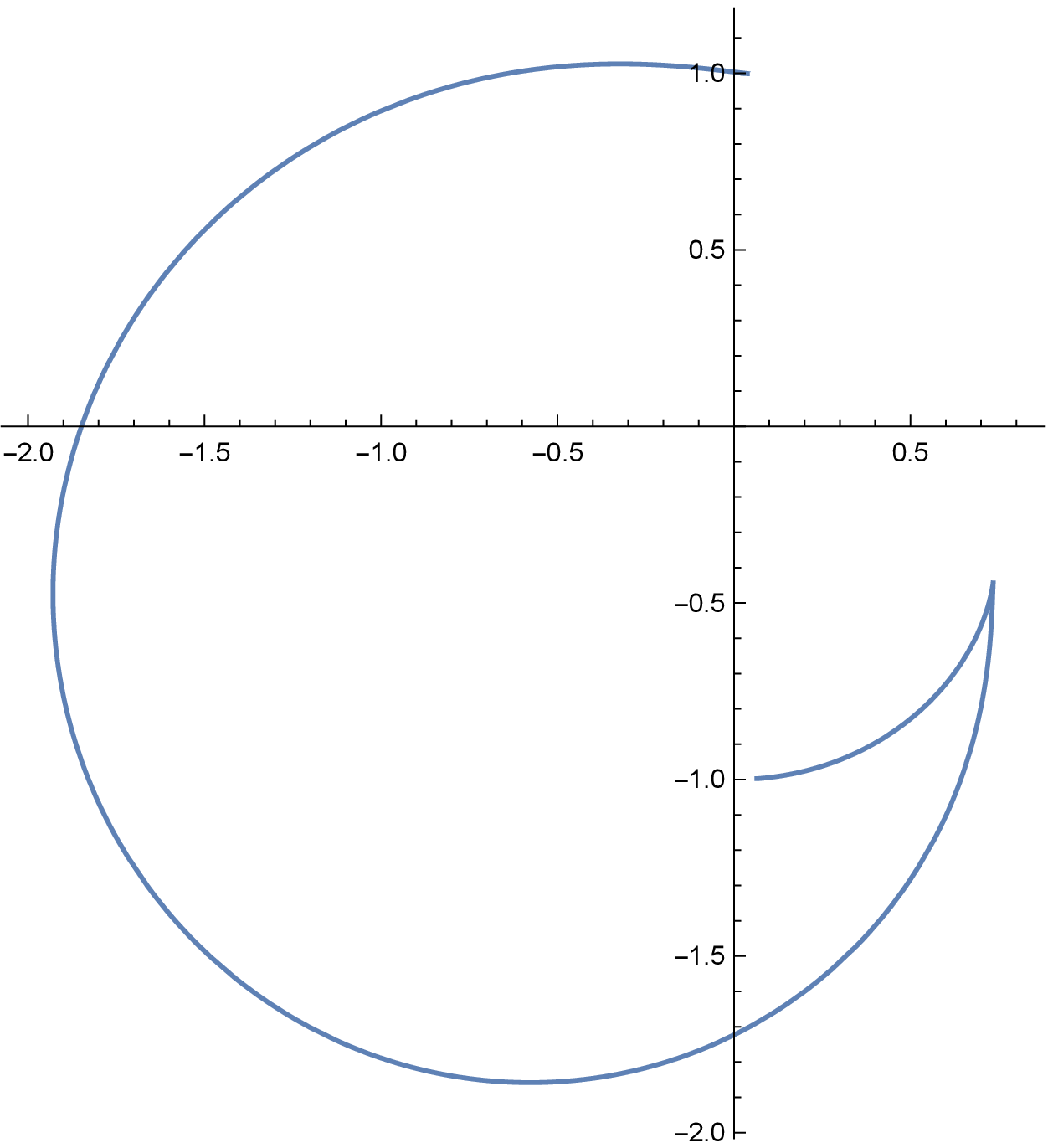}}
\subfigure{\includegraphics[scale=0.45]{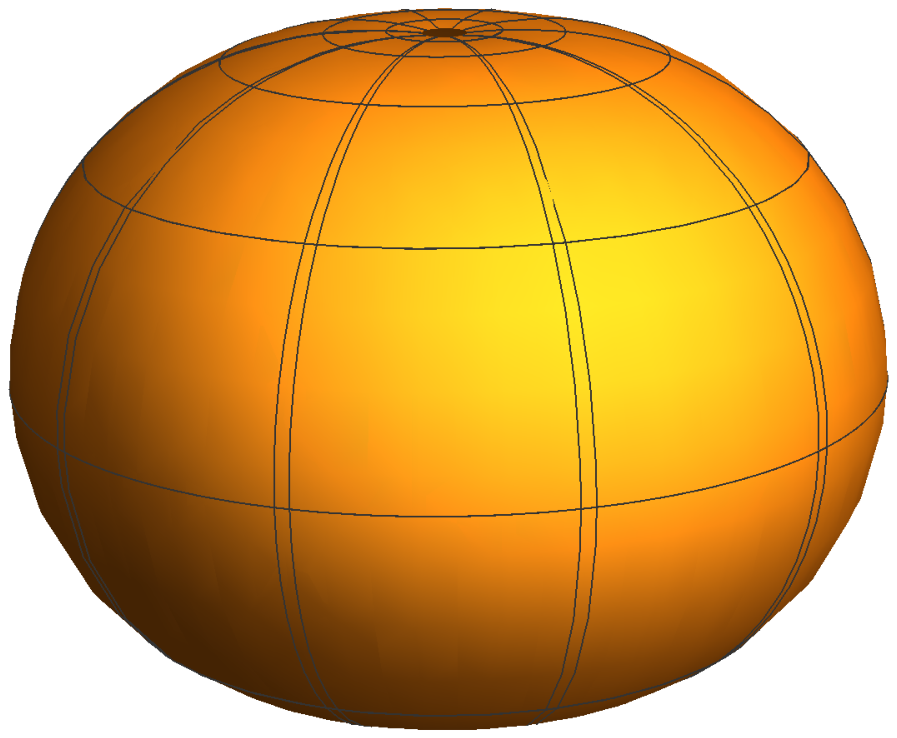}}
\subfigure{\includegraphics[scale=0.45]{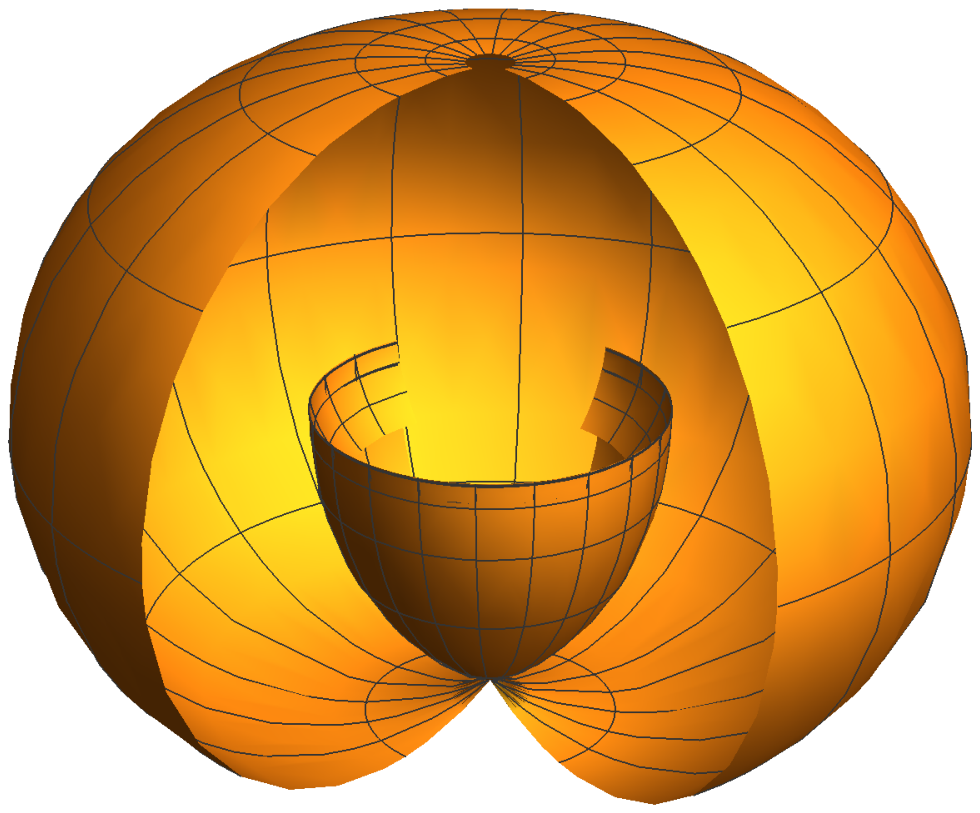}}
\caption{$H_1$-surface of rotation for $a_1 = 1$, $a_2 = 3$ and $c = 1$. }
\end{figure}
	
\end{example}

\newpage

\begin{example}
Considering $a_1 = -1$, $a_2 = 3$ and $c = 1$ in Theorem \eqref{teo rot 1}, we obtain the following $H_1$-surface of rotation. This surface have one isolated singularity and two circles of singularities, where the last one is hidden in the bottom of graphic. The profile curve is given on the left.
\begin{figure}[h]
\centering
\subfigure{\includegraphics[scale=0.30]{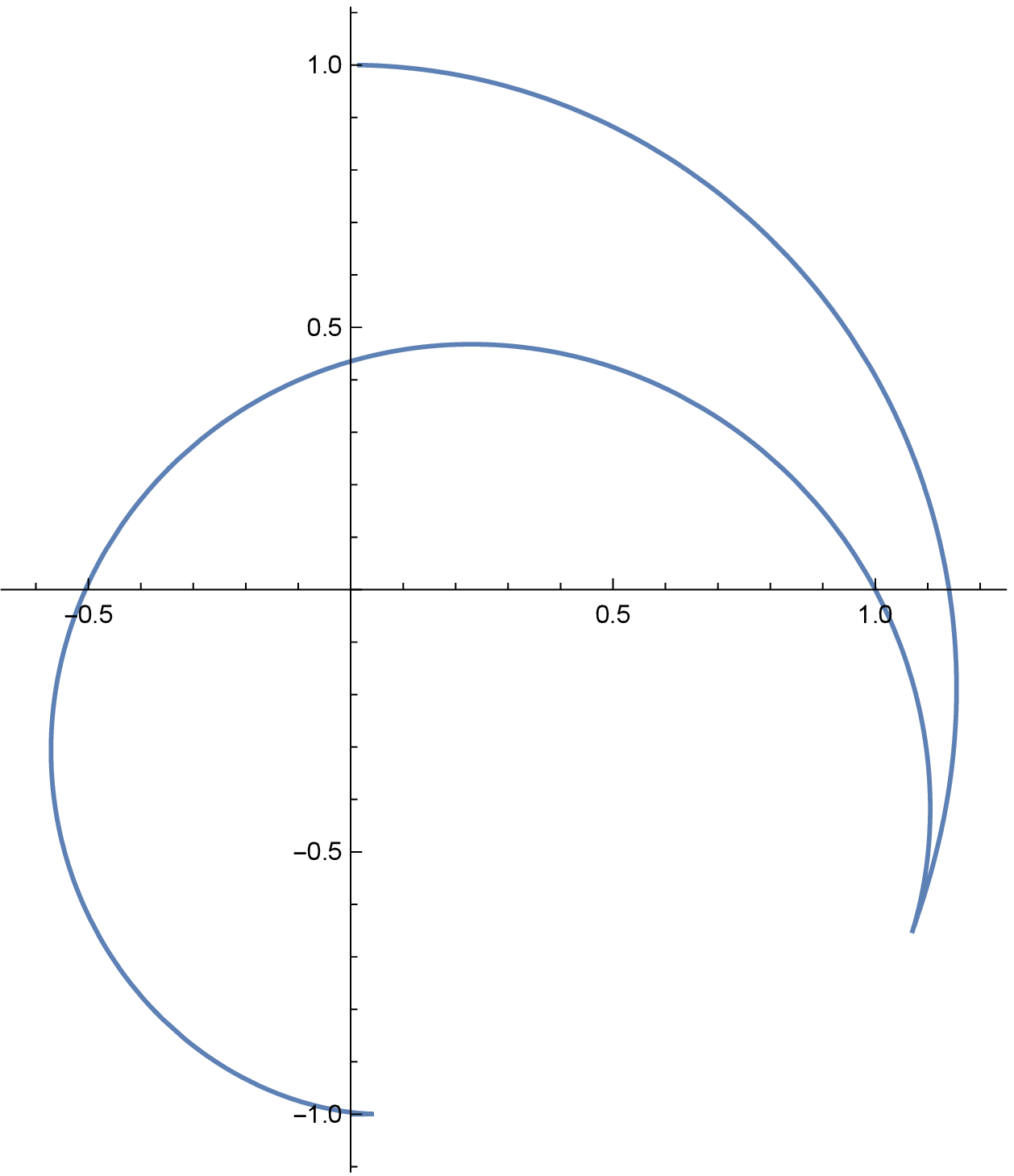}}
\subfigure{\includegraphics[scale=0.45]{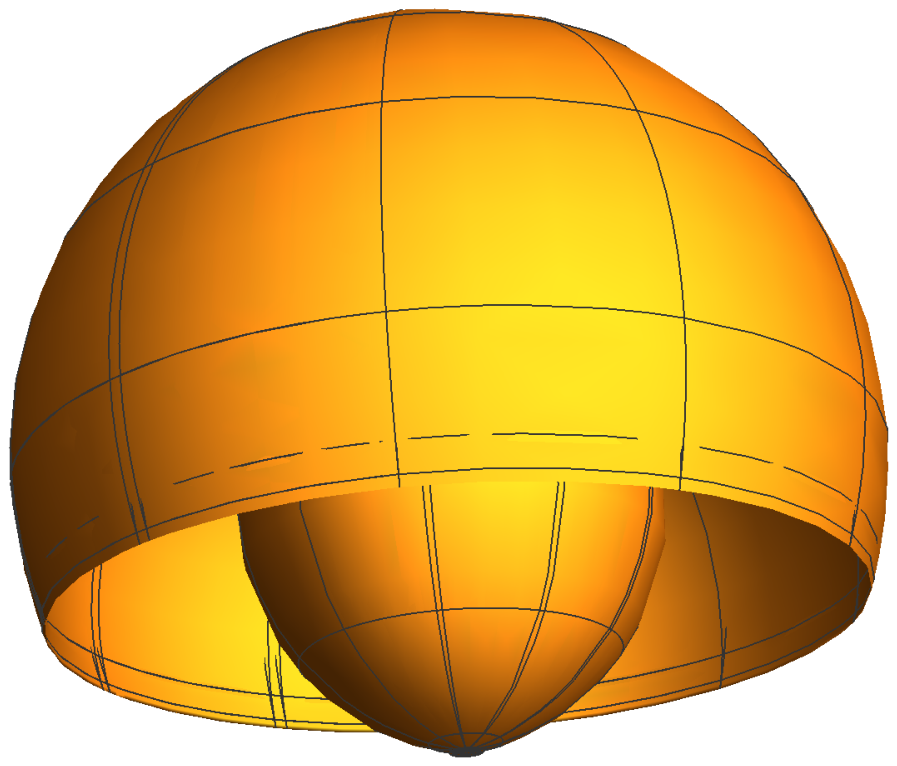}}
\subfigure{\includegraphics[scale=0.45]{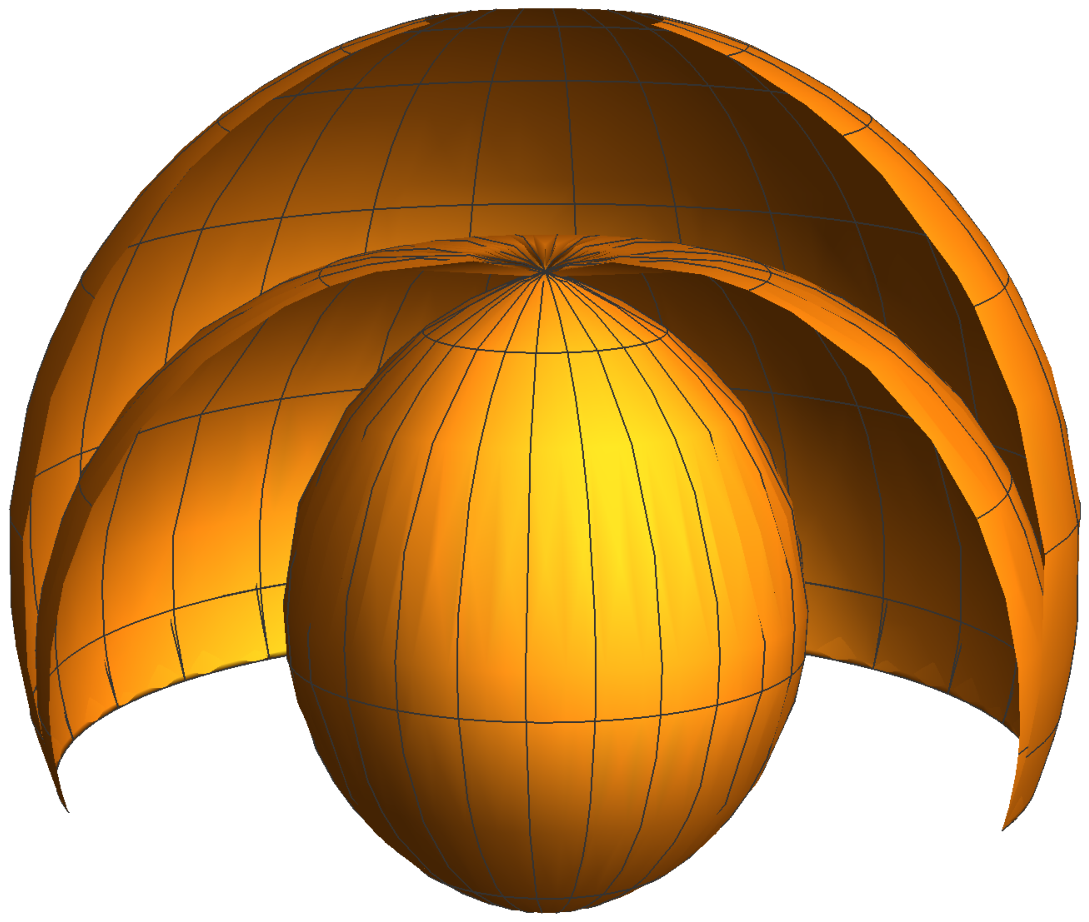}}
\caption{$H_1$-surface of rotation for $a_1 = -1$, $a_2 = 3$ and $c = 1$. }
\end{figure}
	
\end{example}

\vspace{0.5cm}
\begin{example}
Considering $a_1 = -2$, $a_2 = -2$ and $c = 1$ in Theorem \eqref{teo rot 1}, we obtain the next $H_1$-surface of rotation. This surface have one isolated singularity and two circles of singularities. The profile curve is given on the left.
\begin{figure}[h]
\centering
\subfigure{\includegraphics[scale=0.30]{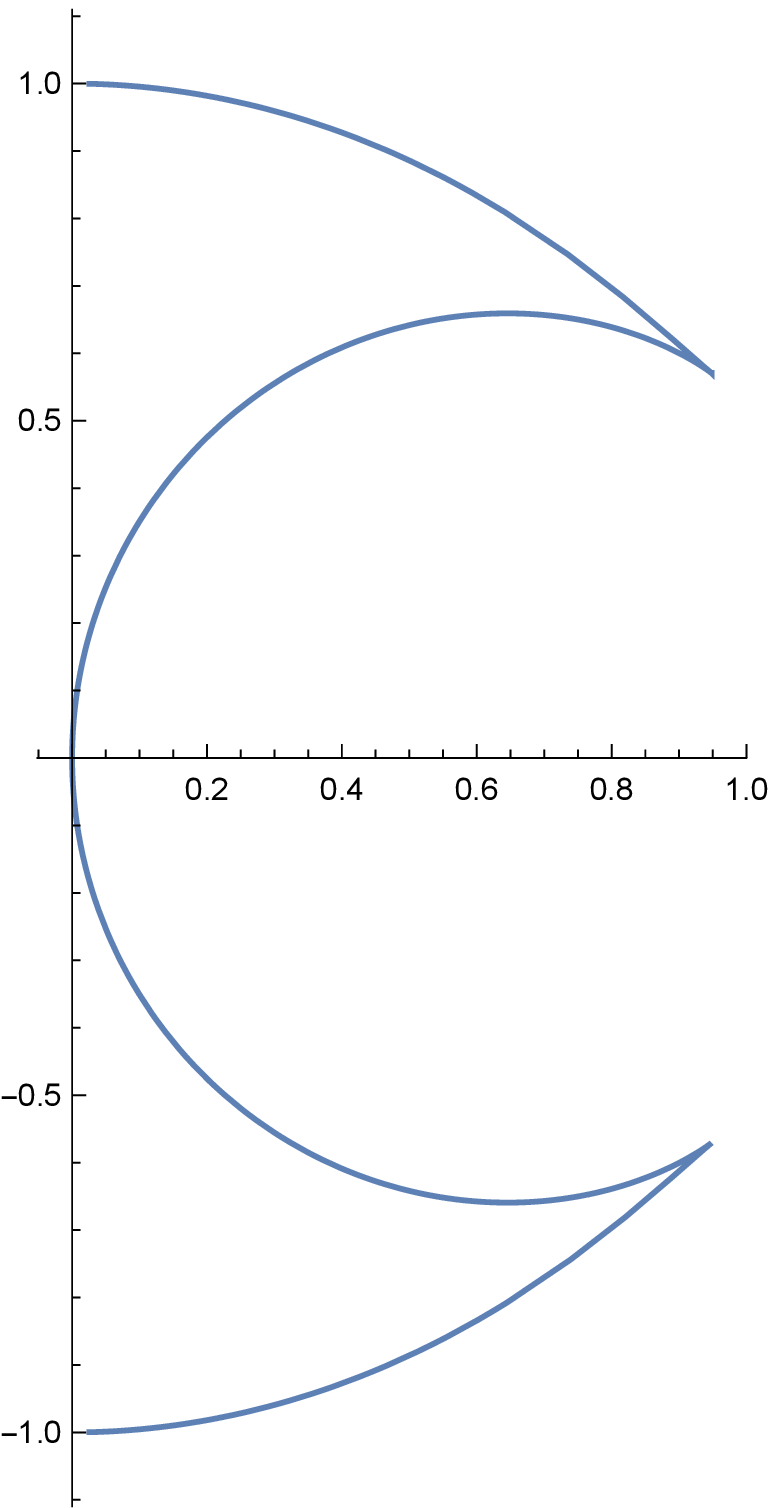}}
\subfigure{\includegraphics[scale=0.45]{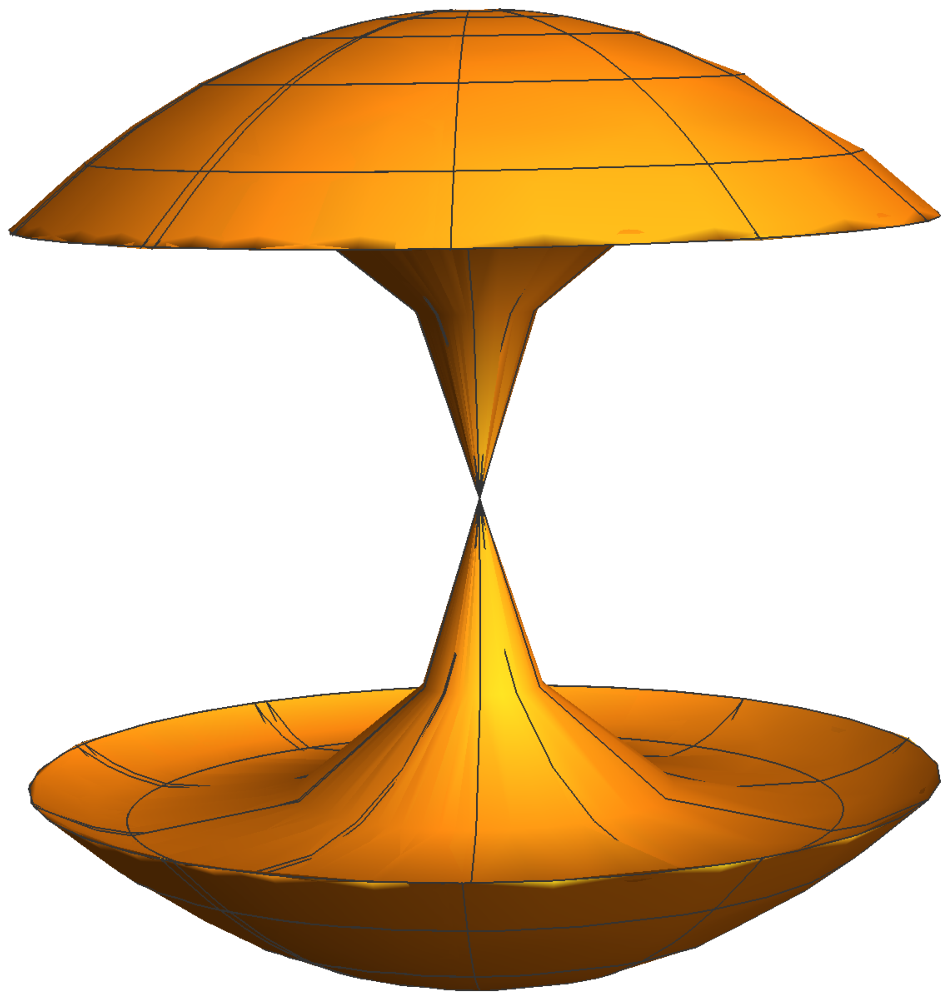}}
\subfigure{\includegraphics[scale=0.45]{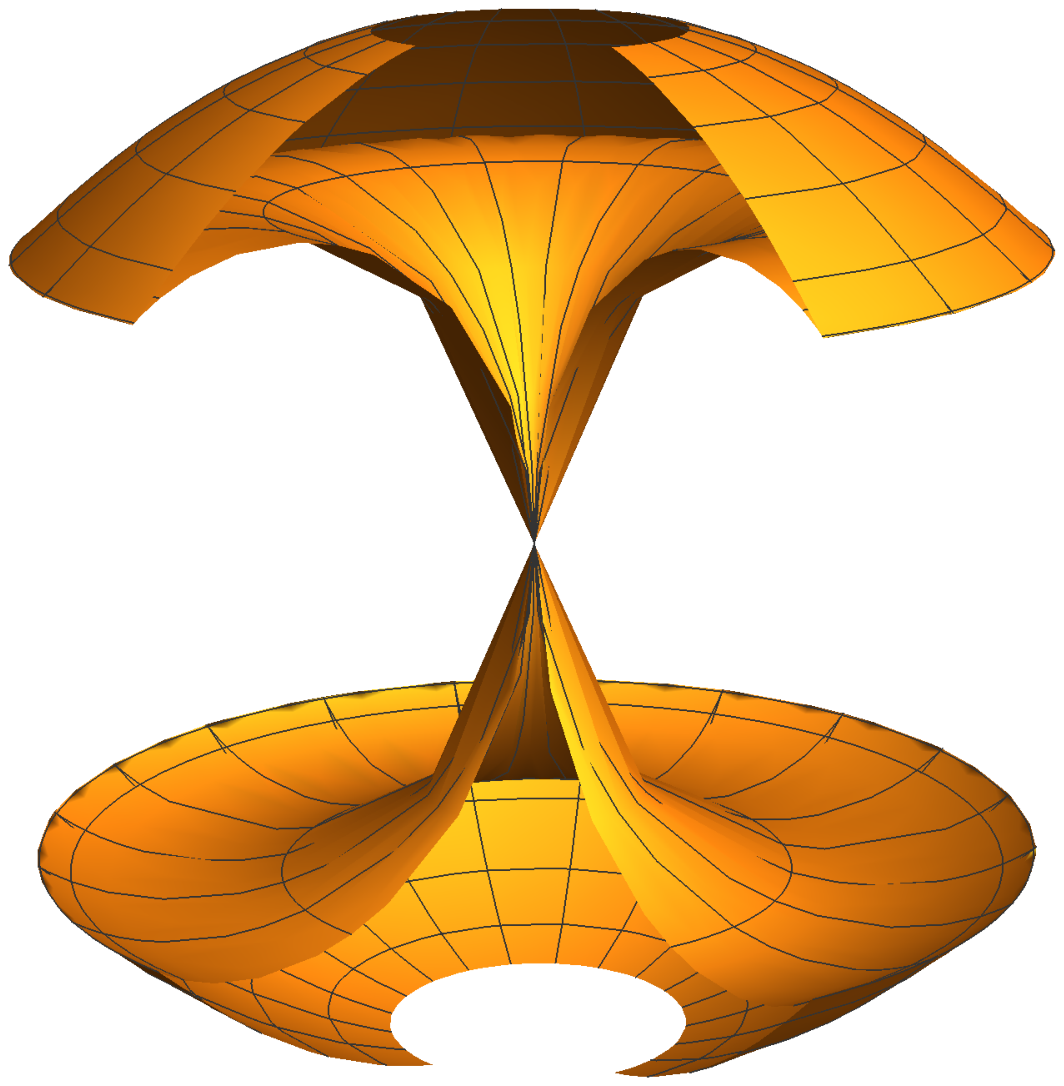}}
\caption{$H_1$-surface of rotation for $a_1 = -2$, $a_2 = -2$ and $c = 1$. }
\end{figure}
	
\end{example}
\newpage
\section{$H_2$-surfaces of Rotation}

Another class of surfaces are defined in the paper \cite{n dimensional}. These are called RSHGW-surfaces and they are described by the generalized Helmholtz equation 

\begin{equation}
\Delta \left[ \frac{T^2}{8|g'|^2}\left(\Delta h + \frac{8|g'|^2}{T^2}h \right) \right ] = 0, \hspace{0.8 cm} T = 1 + |g|^2,
\end{equation}

\hspace{-0.4cm}for $g$ a holomorphic function. They show that functions as in \eqref{h sol} are solutions for the equation above and, when the surface is rotational, then $h$ is radial and becomes

\begin{equation}\label{h sol 4}
h(u,v) = \frac{a_2 + c_1 u + e^{2u}(a_3 + c_2 u)}{1 + e^{2u}}.
\end{equation}

In this way, the functions $h$ as in \eqref{h sol 4} are solutions for the generalized Helmholtz equation \eqref{morango} which describes the $H_2$-surfaces and produces rotational cases. We summarize this discussion in the next theorem.

\begin{theorem}\label{teo rot 2}
A surface $X$ (respectively $\eta$) as in Theorem \eqref{teo g} is a rotational $H_2$-surface (respectively Laguerrre minimal surface) if and only if $g(z) = e^z$ and 
\begin{equation}
h(u,v) = \frac{a_2 + c_1 u + e^{2u}(a_3 + c_2 u)}{1 + e^{2u}}.
\end{equation}
In this case, the surfaces $X$ and $\eta$ can be locally parameterized by

\begin{eqnarray}\label{eq X rot}\label{eq eta rot}
X(u,v) &=& \big(\tilde{M}(u)\cos v\hspace{0.1cm},\hspace{0.1cm} \tilde{M}(u) \sin v\hspace{0.1cm}, \hspace{0.1cm} \tilde{N}(u) \big),\nonumber \\
\eta(u,v) &=& \big(\tilde{M}_1(u)\cos v\hspace{0.1cm},\hspace{0.1cm} \tilde{M}_1(u) \sin v\hspace{0.1cm}, \hspace{0.1cm} \tilde{N}_1(u) \big),
\end{eqnarray}
where
\begin{eqnarray*}
\tilde{M}(u) &=& \frac{2e^u\big[ q_1 q_2 - 2a_2 q_3 - 2c\big(q_4 + q_5\big)\big]}{q_6 - 2c_1e^{2u}q_7 + e^{4u}q_8}, \\ 
\tilde{N}(u) &=& \frac{c_1^2\big(1 - e^{2u}(1-2u)^2\big) + 2c_1e^{2u}r_1 + e^{2u}r_2}{r_3 - 2c_1e^{2u}r_4 + e^{4u}r_5}.\\
 &   & \\
\tilde{M_1}(u) &=& \frac{e^{-u}\Big(e^{2u}\big(2a_2 + 2a_3 + c_2 - c_2e^{2u} + 2c_2u\big) + c_1\big(1+e^{2u}(2u-1)\big)\Big)}{2(1+e^{2u})}, \\
\tilde{N_1}(u) &=& \frac{a_2 + c_1(u-1) - e^{2u}(a_3 + c_2 + c_2u)}{1+e^{2u}}, \hspace{0.3cm}a_2, a_3, c_1, c_2 \in \mathbb{R}^3
\end{eqnarray*}
with

\begin{eqnarray*}
q_1 & = & c_1(1-2u) + c_2e^{2u}, \\
q_2 & = & c_1 +e^{2u}(2a_3 + c_2 + 2c_2 u),\\ 
q_3 &=& c_1 +e^{2u}(2a_3 + 2c +  c_2 + 2c_2 u), \\
q_4 &=& e^{2u}(2a_3 + c_2 - c_2e^{2u} + 2 c_2 u),\\
q_5 &=& c_1\big(1+e^{2u}(2u - 1)\big), \\
q_6 &=& 4a_2^2e^{2u} - 4a_2c_2e^{4u} + c_1^2(1+e^{2u}(1-2u)^2),\\ 
q_7 &=& -2a_3 + a_2(2 - 4u) + c_2(-1 - 2u + e^{2u}(2u - 1))\\
q_8 &=& 4a_3^2 + 4a_3(c_2 + 2 c_2 u) + c_2^2(e^{2u} + (1 + 2u)^2).\\
 &  & \\
r_1 & = & 2a_3 + 4c + c_2 - c_2e^{2u} + a_2(2-4u) - 4cu +2c_2 u + 2c_2e^{2u}u, \\
r_2 & = & -4a_2^2 + a_2(4c_2e^{2u}-8c) + e^{2u}\big(4a_3^2 + 4a_3(2c + c_2 + 2c_2u)\big) + c_2\big(8c(u+1) + c_2( (1+2u)^2 - e^{2u} \big), \\
r_3 & = & 4a_2^2e^{2u} - 4a_2c_2e^{4u} + c_1^2(1 + e^{2u}(1-2u)^2), \\
r_4 & = & -2a_3 + a_2(2-4u) + c_2\big(e^{2u}(2u-1) -1 - 2u \big), \\
r_5 & = & 4a_3^2 + 4a_3(c_2 + 2c_2u) + c_2^2\big(e^{2u} + (1+2u)^2\big).
\end{eqnarray*}

\end{theorem}
\vspace{0.5cm}
\begin{proof}
As in the previous proposition, a $H_2$-surface $X$ parameterized as in \eqref{eq X g} is rotational if and only if $h$ is radial and $g(z) = e^z$.  Making the change of parameters $ w= e^z, \hspace{0.3cm} z=u+iv \in \mathbb{C}$, we have $g(z)= e^z$ and $ h_{,2}=0$. From Corollary \eqref{H2 h}, we have $h$ as in \eqref{h rot}. Differentiating \eqref{h rot} with respect to $v$, considering that $ h_{,2}=0$, we obtain

$$
\left< 1, iA' \right> + \left< e^z, i(B' - B) \right> = 0
$$

Using \eqref{pro cte}, we get

$$
iA' = -\bar{z_1}e^z + ic_1, \hspace{1.0cm} i(B' - B) = ic_2e^z + z_1,
$$
whose solutions are given by

$$
A = i\bar{z_1}e^z + c_1 z + z_2, \hspace{1.0cm} B = (z_3 + c_2z)e^z + iz_1,
$$
where $c_1, c_2 \in \mathbb{R}$, $z_k = a_k +ib_k$, \hspace{0.1cm} $a_k, b_k \in \mathbb{R}$, $k = 1,2,3$. Substituing this expressions in \eqref{h rot}, we get the function $h$ in the statement.

Now substituing these expressions for $g$ and $h$ in \eqref{eq X g}, we get \eqref{eq X rot}.	
			
\end{proof}

\newpage

\subsection{Examples of Rotational $H_2$-Surfaces}
\begin{example}
Considering $a_2 = a_3 = c_1 = c_2 = c = 1$ in Theorem \eqref{teo rot 2}, we obtain a rotational $H_2$-surface with two isolated singularities and three circles of singularities. The profile curve is on the left.		

\begin{figure}[h]
\centering
\subfigure{\includegraphics[scale=0.30]{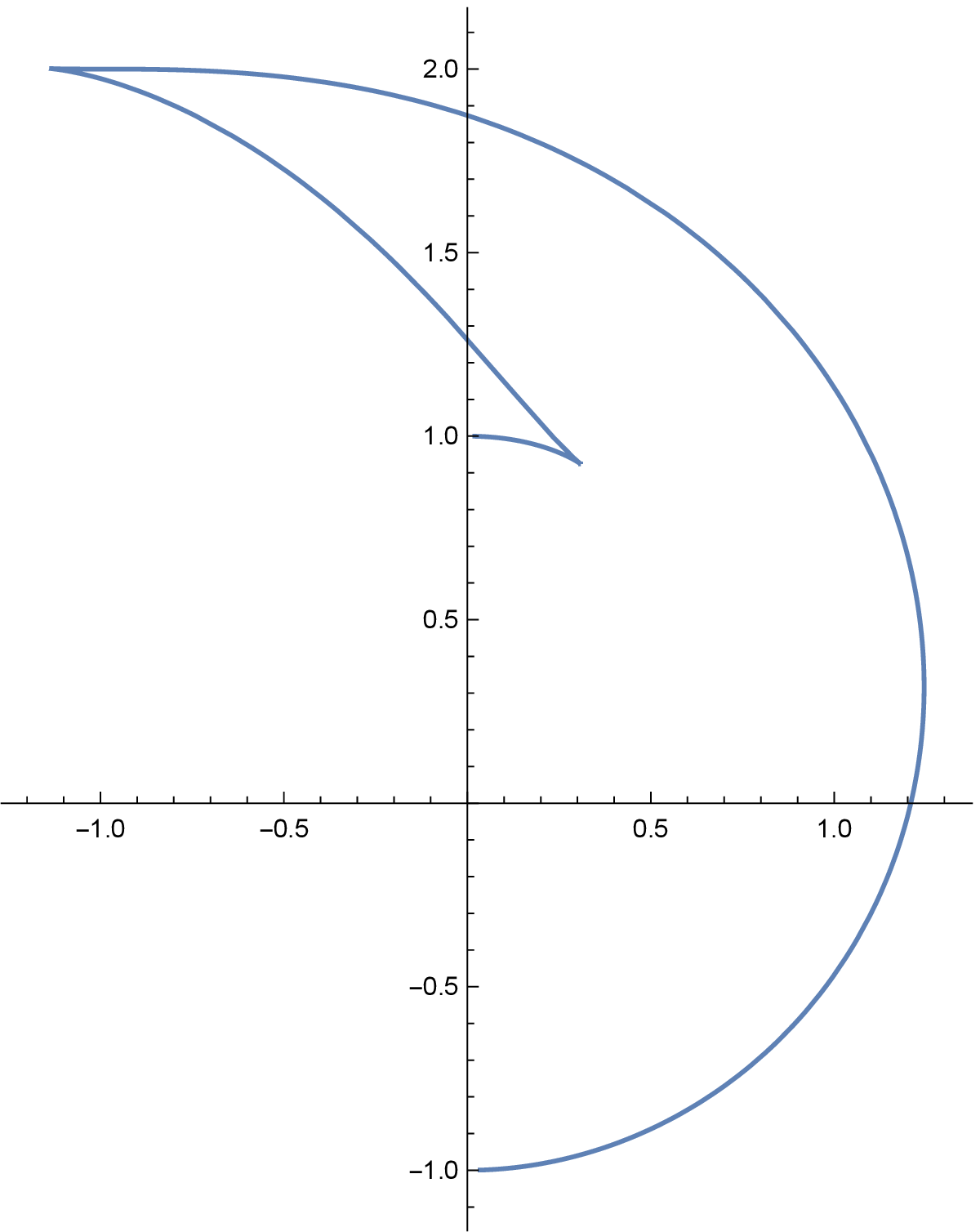}}\hspace{1.0cm}
\subfigure{\includegraphics[scale=0.40]{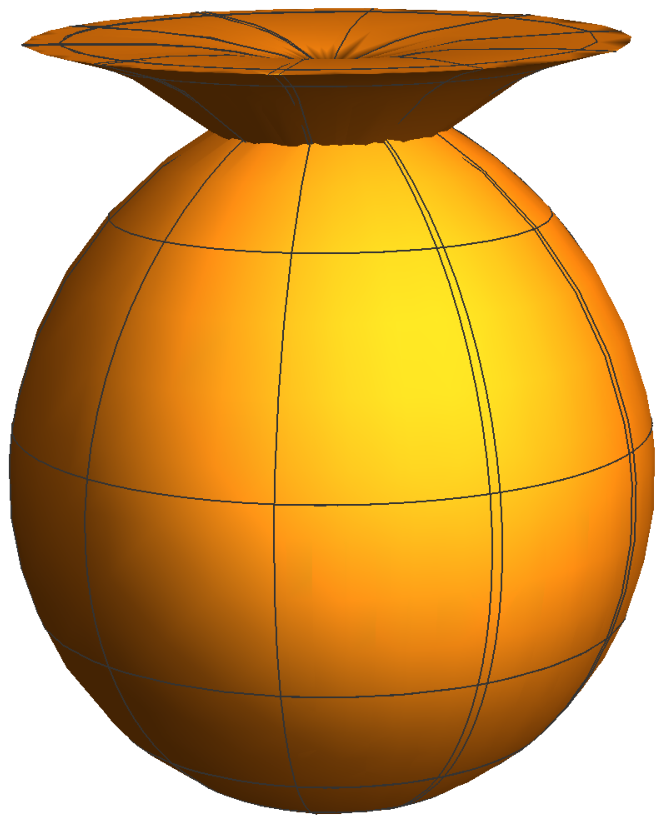}}\hspace{-0.5cm}
\subfigure{\includegraphics[scale=0.40]{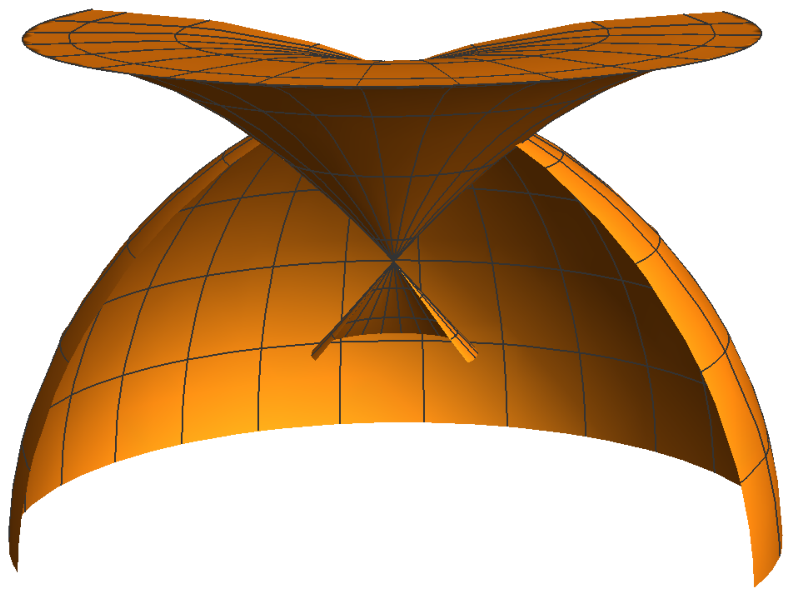}}
\subfigure{\includegraphics[scale=0.40]{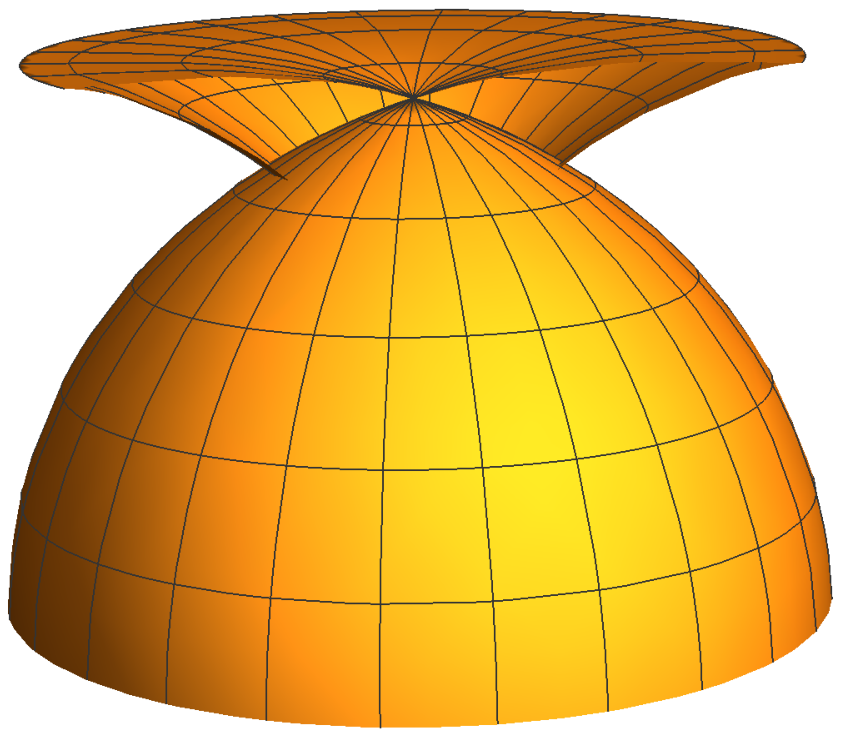}}
\caption{Rotational $H_2$-surface for $a_2 = a_3 = c_1 = c_2 = c = 1$}
\end{figure}

\end{example}
\vspace{0.5cm}
\begin{example}
Considering $a_2 = -2$, $ a_3 = -2$, $ c_1 = -1$, $ c_2 =1$ and $ c = 1$ in Theorem \eqref{teo rot 2}, we obtain a rotational $H_2$-surface with two circles of singularities and one isolated singularity The profile curve is on the left.	
	
\begin{figure}[h]
\centering
\subfigure{\includegraphics[scale=0.30]{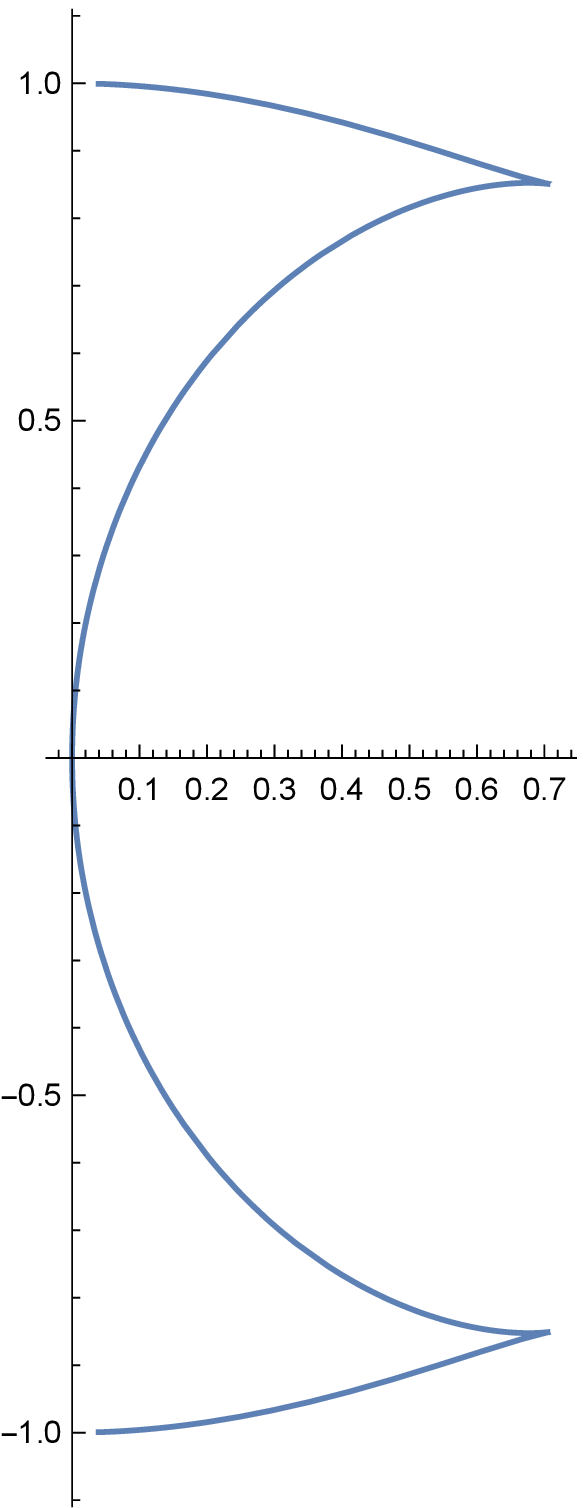}}\hspace{1.0cm}
\subfigure{\includegraphics[scale=0.40]{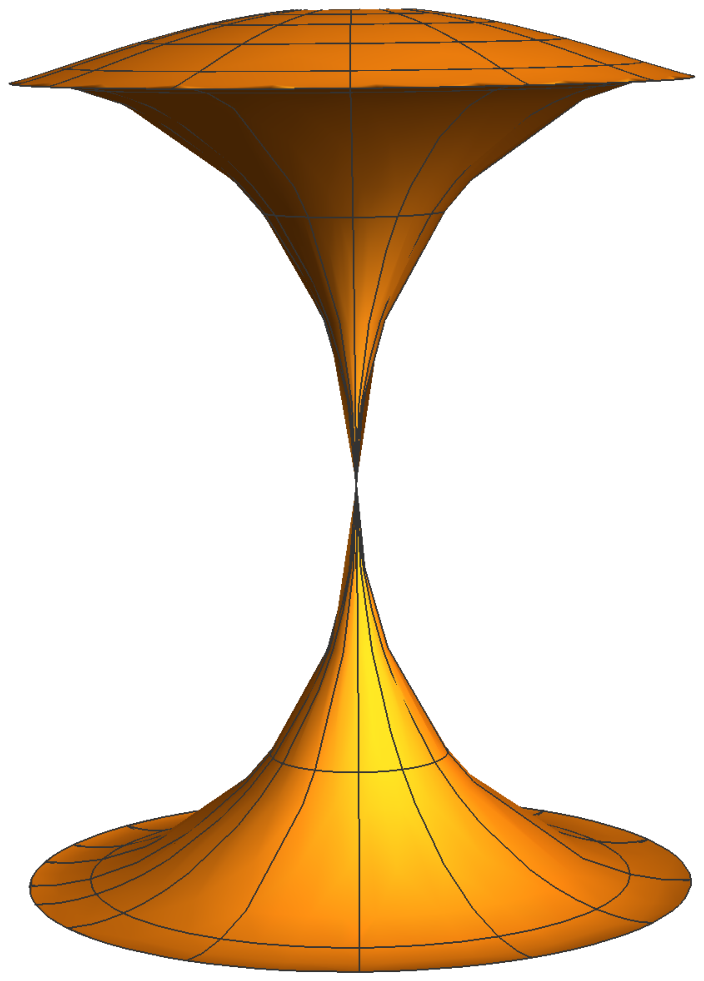}}\hspace{-0.5cm}
\subfigure{\includegraphics[scale=0.40]{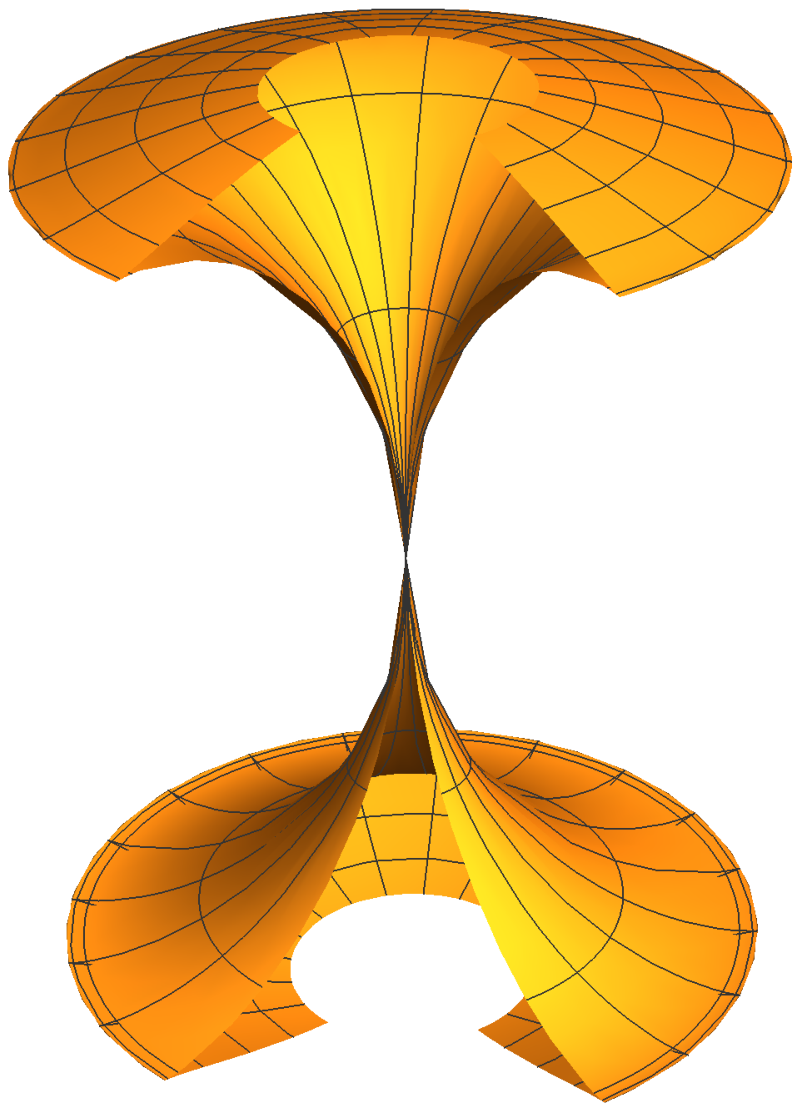}}
\caption{Rotational $H_2$-surface for $a_2 = -2$, $ a_3 = -2$, $ c_1 = -1$, $ c_2 =1$ and $ c = 1$ }
\end{figure}

\end{example}
\newpage
\begin{example}
Considering $a_2 = -1$, $ a_3 = 2$, $ c_1 = 2$, $ c_2 =-1$ and $ c = 1$ in Theorem \eqref{teo rot 2}, we obtain a rotational $H_2$-surface with two circles of singularities and one isolated singularity. The profile curve is on the left.	
	
\begin{figure}[h]
\centering
\subfigure{\includegraphics[scale=0.30]{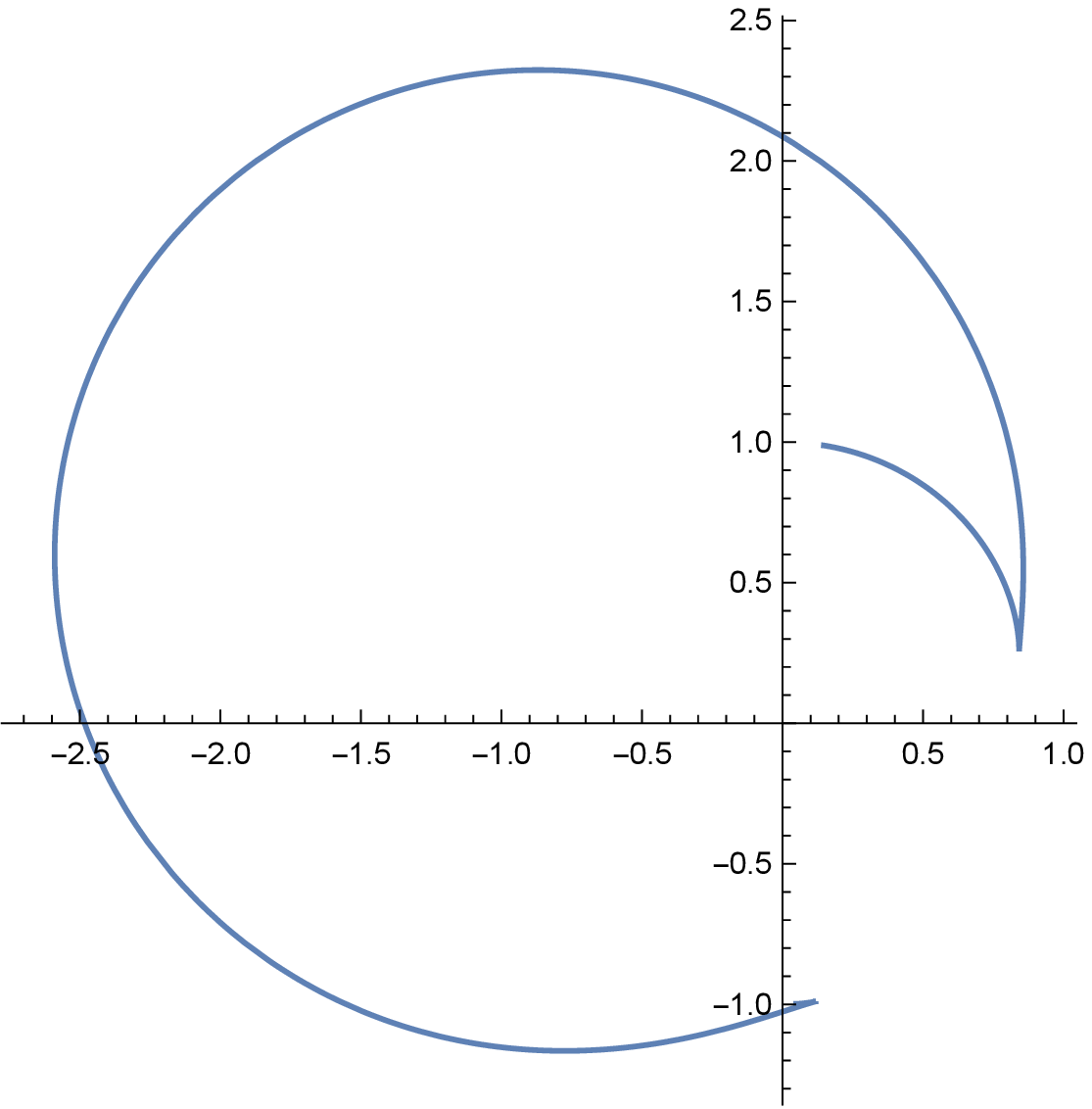}}\hspace{1.0cm}
\subfigure{\includegraphics[scale=0.40]{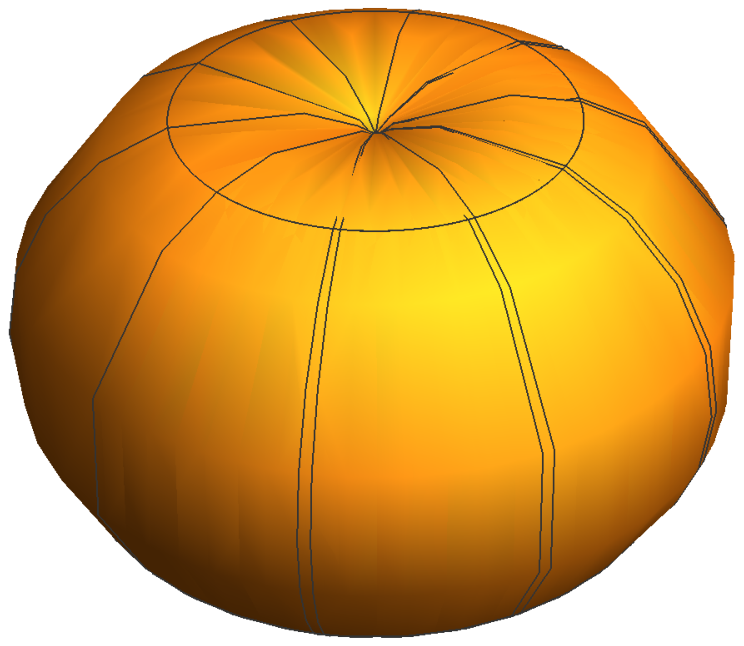}}
\subfigure{\includegraphics[scale=0.40]{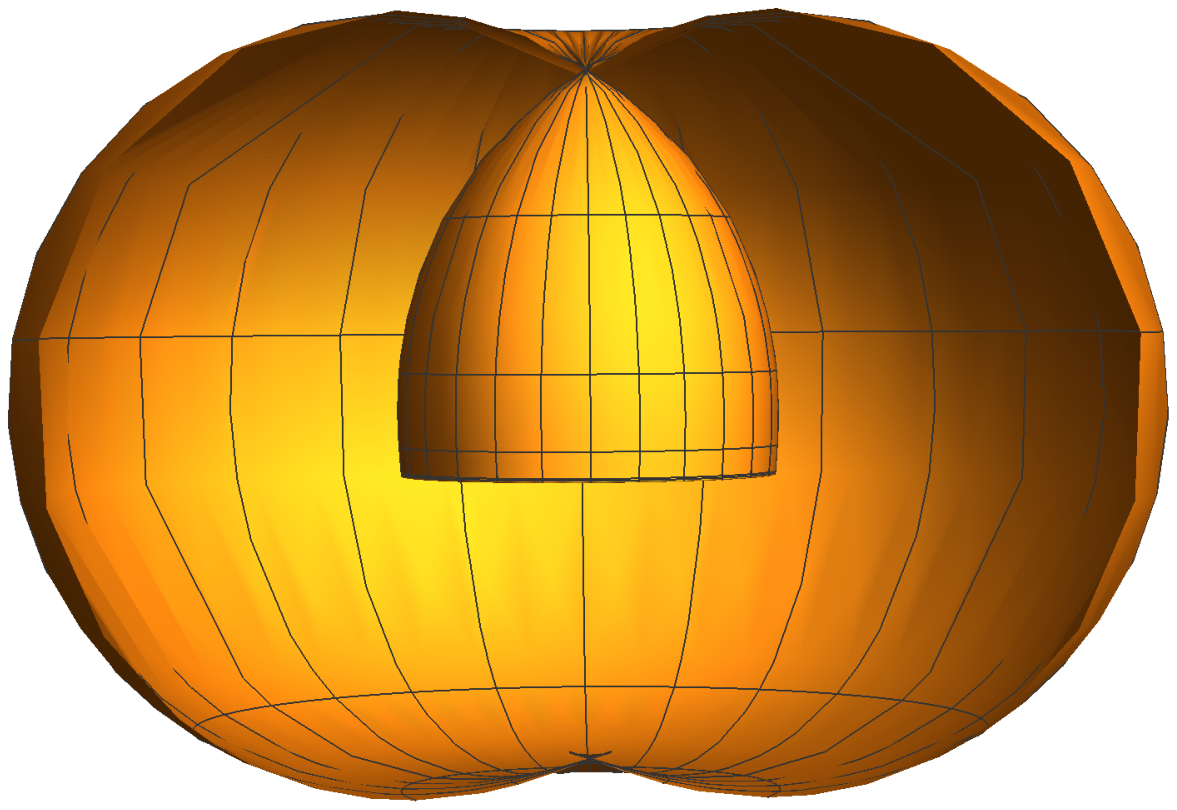}}
\caption{Rotational $H_2$-surface for $a_2 = -1$, $ a_3 = 2$, $ c_1 = 2$, $ c_2 = -1$ and $ c = 1$ }
\end{figure}

\end{example}

\vspace{0.5cm}
\begin{example}
Considering $a_2 = -1$, $ a_3 = 2$, $ c_1 =1$, $ c_2 =3$ and $ c = 1$ in  Theorem \eqref{teo rot 2}, we obtain a rotational $H_2$-surface with two circles of singularities and no isolated singularity. The profile curve is on the left.	

\begin{figure}[h]
\centering
\subfigure{\includegraphics[scale=0.30]{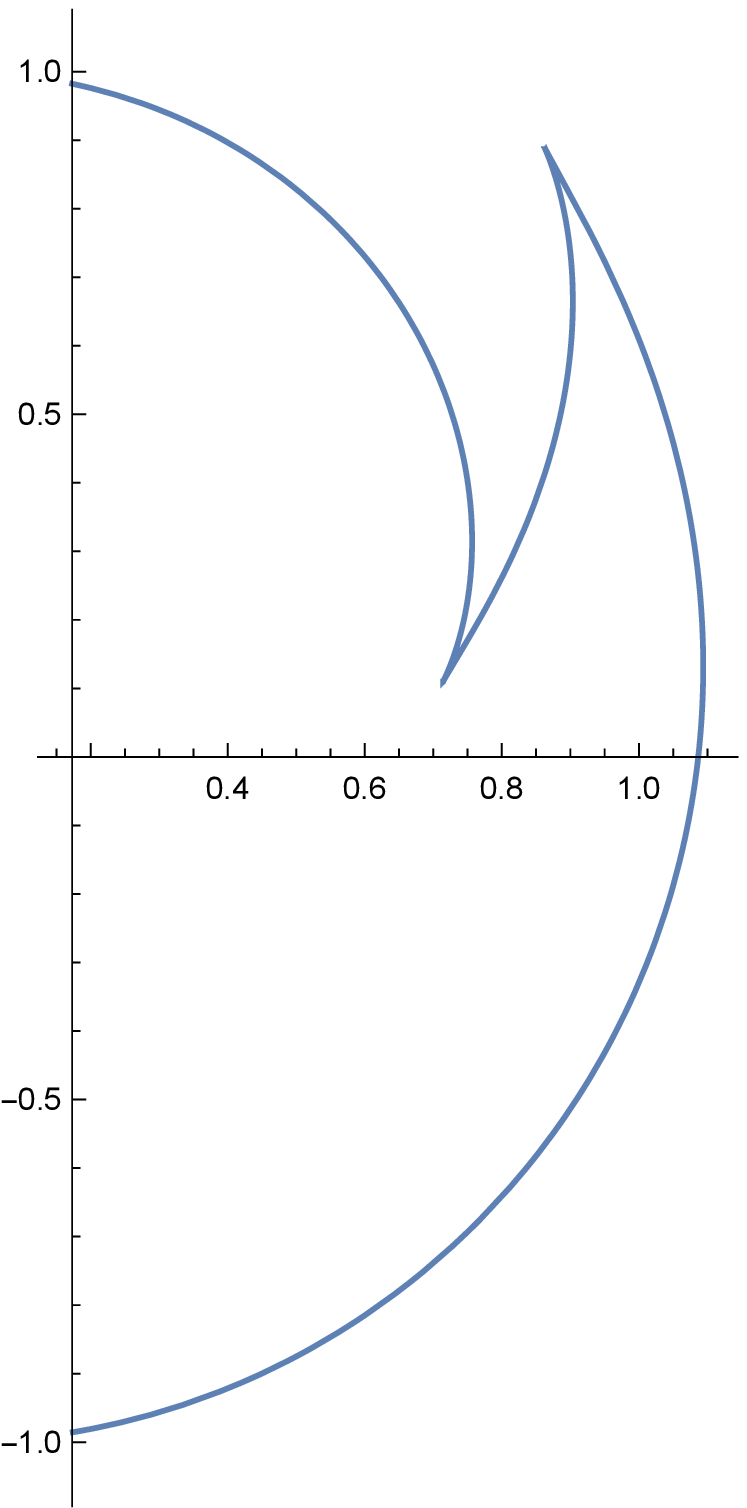}}\hspace{0.5cm}
\subfigure{\includegraphics[scale=0.40]{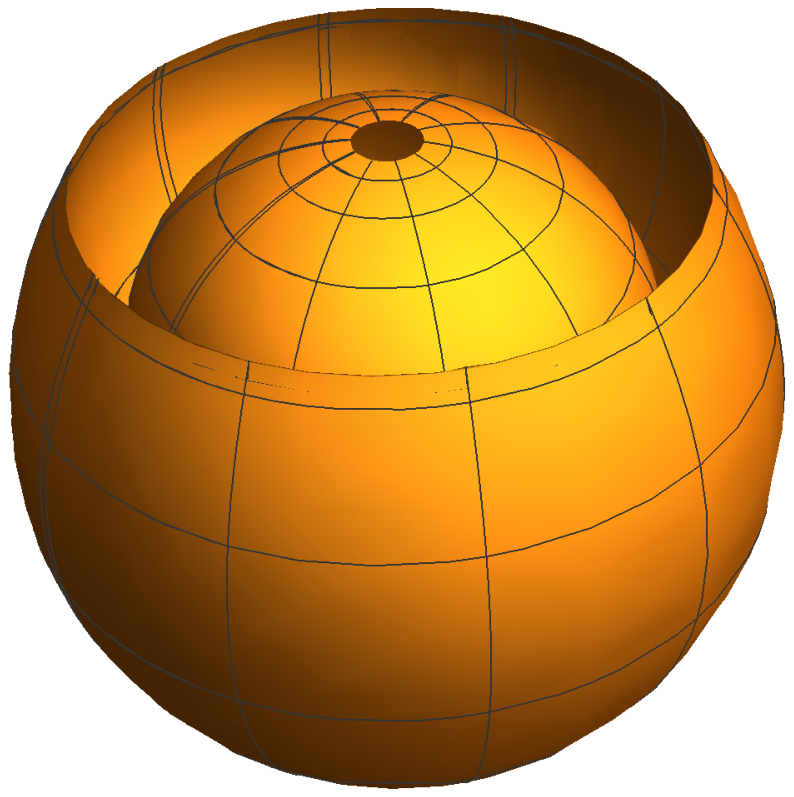}}\hspace{-0.5cm}
\subfigure{\includegraphics[scale=0.40]{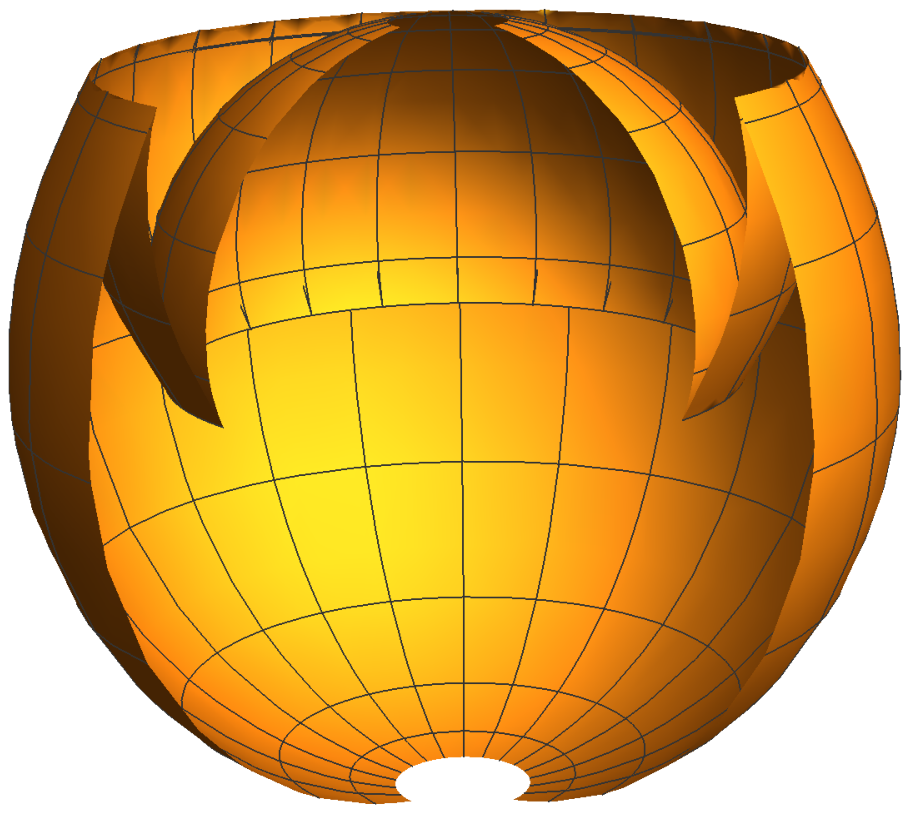}}
\caption{Rotational $H_2$-surface for $a_2 = -1$, $ a_3 = 2$, $ c_1 =1$, $ c_2 =3$ and $ c = 1$ }
\end{figure}

\end{example}

\newpage
\begin{example}
Considering $a_2 = -2$, $ a_3 = 2$, $ c_1 = -1$, $ c_2 = 1$ and $ c = 1$ in Theorem \eqref{teo rot 2}, we obtain a rotational $H_2$-surface with two circles of singularities and no isolated singularity. The profile curve is on the left.	
	
\begin{figure}[h]
\centering
\subfigure{\includegraphics[scale=0.40]{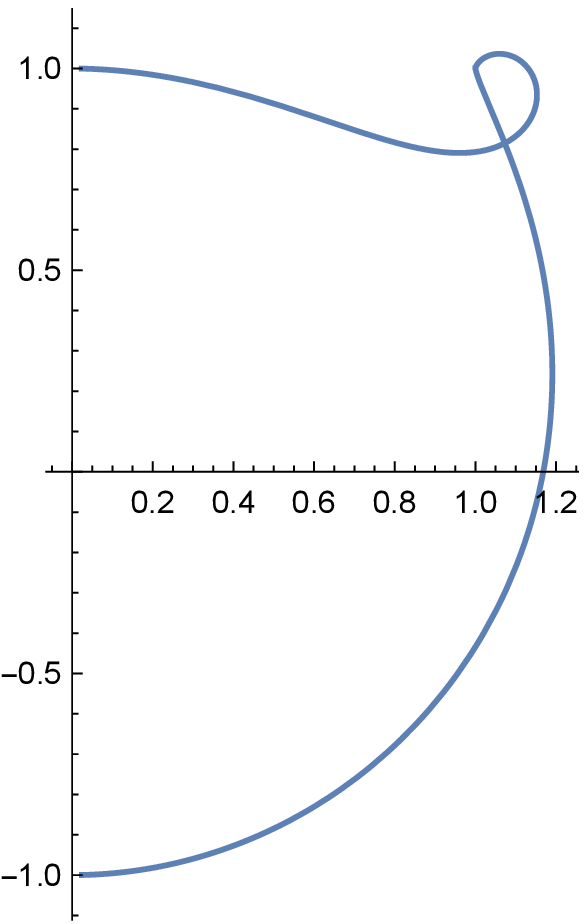}}\hspace{0.5cm}
\subfigure{\includegraphics[scale=0.40]{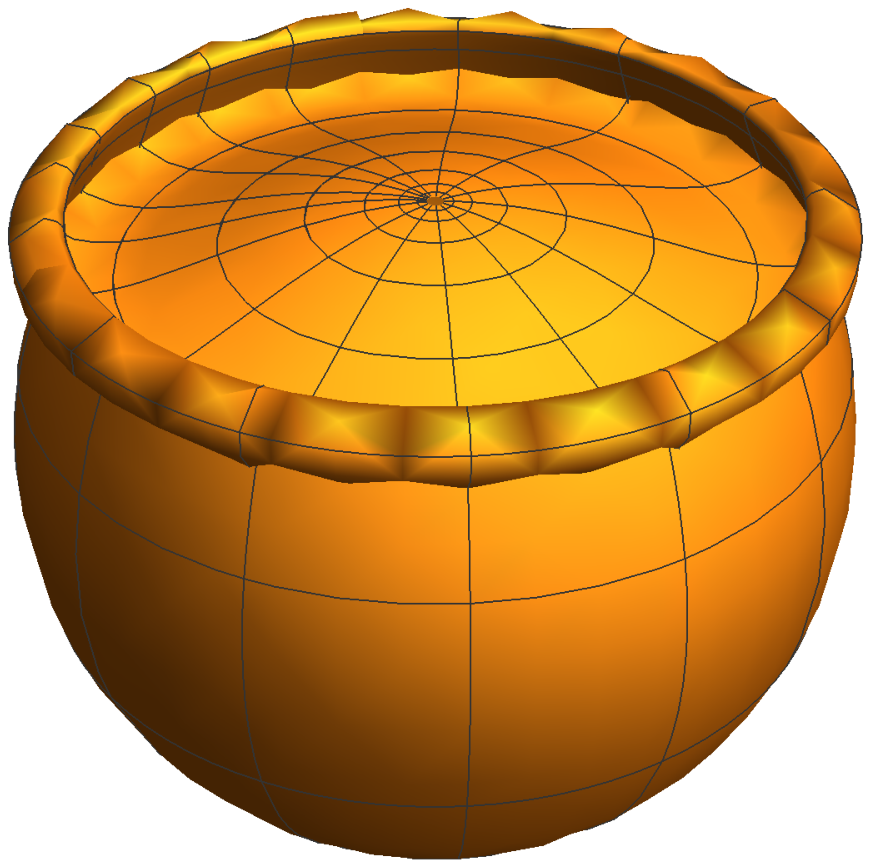}}\hspace{0.5cm}
\subfigure{\includegraphics[scale=0.35]{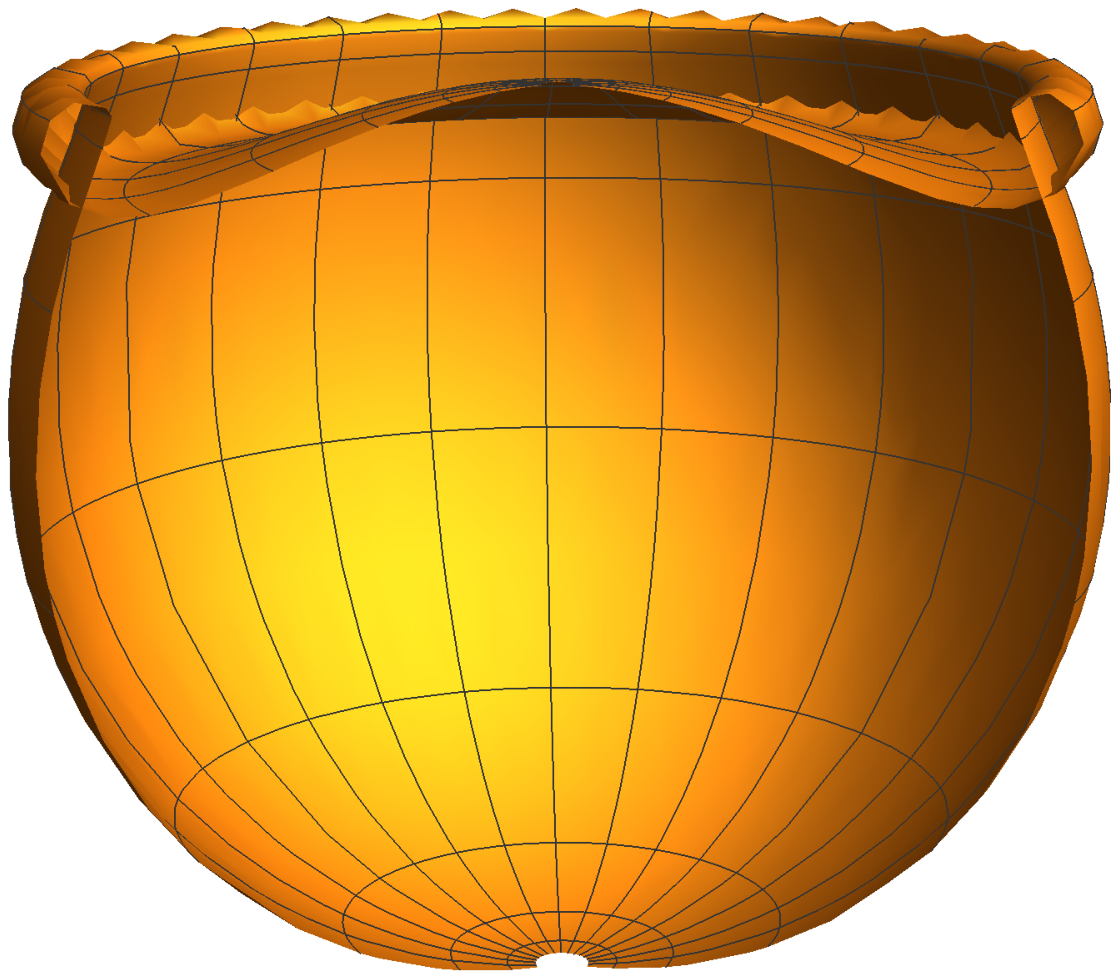}}
\caption{Rotational $H_2$-surface for $a_2 = -2$, $ a_3 = 2$, $ c_1 = -1$, $ c_2 = 1$ and $ c = 1$ }
\end{figure}
\end{example}

\vspace{0.5cm}
\begin{example}
Considering $a_2 = -1$, $ a_3 = 2$, $ c_1 = 0$, $ c_2 = 0$ and $ c = 1$ in Theorem \eqref{teo rot 2}, the rotational $H_2$-surface is the sphere below. 	
	
\begin{figure}[h]
\centering
\subfigure{\includegraphics[scale=0.40]{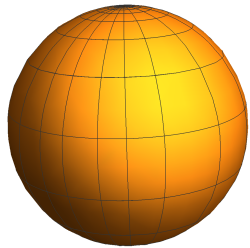}}\hspace{0.5cm}
\caption{Rotational $H_2$-surface for $a_2 = -1$, $ a_3 = 2$, $ c_1 = 0$, $ c_2 = 0$ and $ c = 1$ }
\end{figure}
\end{example}

\newpage
\vspace{1.0cm}
\subsection{Examples of Rotational Laguerre Minimal Surfaces}

One can see by means of regularity condition $det V \neq 0 $ and profile curve $\alpha(u) = (\widetilde{M_1}(u), \widetilde{N_1}(u))$ in \eqref{eq eta rot} that holds:

\vspace{1.0cm}
\begin{itemize}
\item For $c_1c_2 > 0$, the surface $ \eta $ has at least one isolated singularity and one circle of singularities.
\vspace{0.5cm}

\item For $c_1c_2 < 0$, may occur complete cases, isolated singularities or circles of singularities. 
\vspace{0.5cm}

\item For $c_1 = 0$ and $c_2 \neq 0$,  the surface $ \eta $ always has singularity.
\vspace{0.5cm}

\item For $c_1 \neq 0$ and $c_2 = 0$, the surface $ \eta $ has at least one circle of singularities and isolated singularities may or  may not occur .
\vspace{0.5cm}

\item For $ c_1 = c_2 = 0 $, the surface $ \eta $ is a sphere.
\end{itemize}

\vspace{0.5cm}
The following examples illustrate each of these cases.

\vspace{0.5cm}
\begin{example}
Considering $a_2 = a_3 = c_1 = c_2 = 1 $ in Theorem \eqref{teo rot 2}, the Laguerre minimal surface $ \eta $ has one isolated singularity and one circle of singularities. The profile curve is on the left.

\begin{figure}[h]
\centering
\subfigure{\includegraphics[scale=0.30]{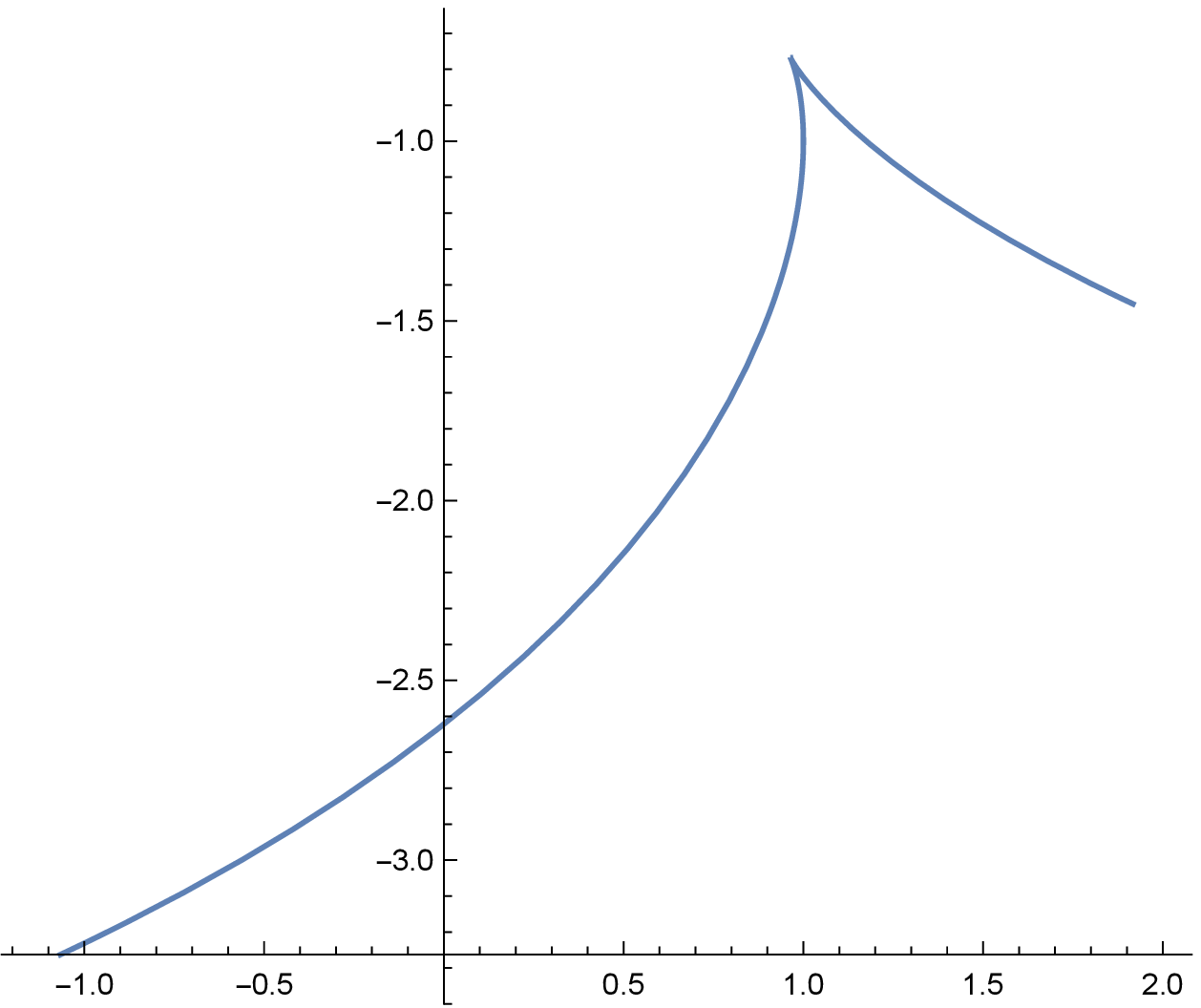}}\hspace{0.5cm}
\subfigure{\includegraphics[scale=0.45]{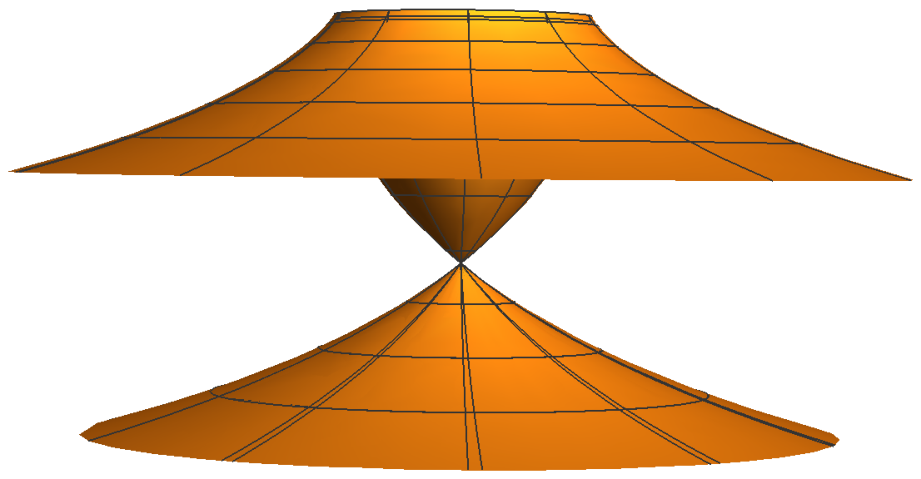}}\hspace{-1.0cm}
\subfigure{\includegraphics[scale=0.40]{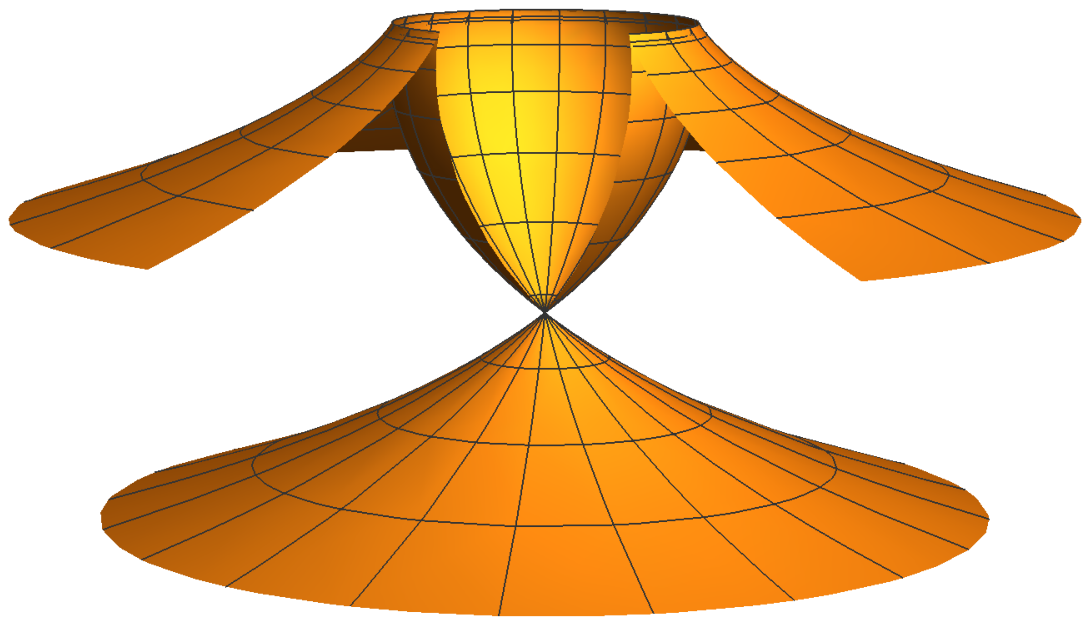}}
\caption{Rotational Laguerre Minimal Surface for $a_2 = a_3 = c_1 = c_2 =1$}
\end{figure}	
\end{example}

\vspace{0.5cm}
\begin{example}
Considering $a_2 = -2$, $ a_3 = -2$, $ c_1 = -1$, $ c_2 =1$  in  Theorem \eqref{teo rot 2}, we obtain a rotational Laguerre minimal surface with no isolated singularities and one circle of singularities. The profile curve is on the left.

\begin{figure}[h]
\centering
\subfigure{\includegraphics[scale=0.30]{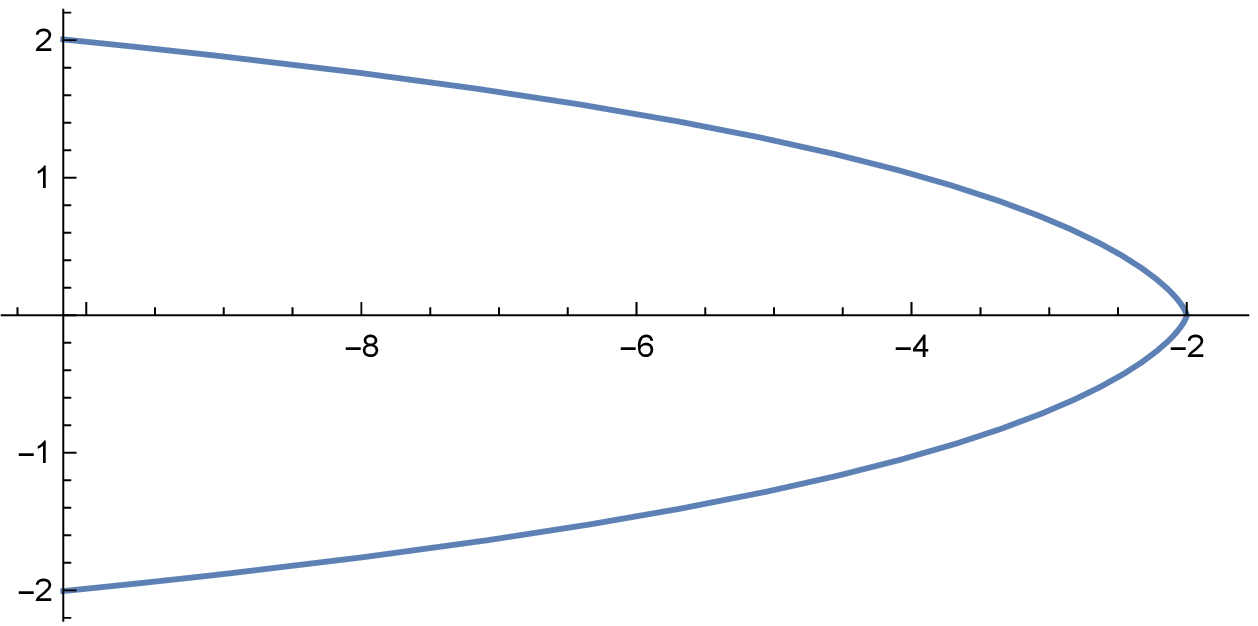}}
\subfigure{\includegraphics[scale=0.40]{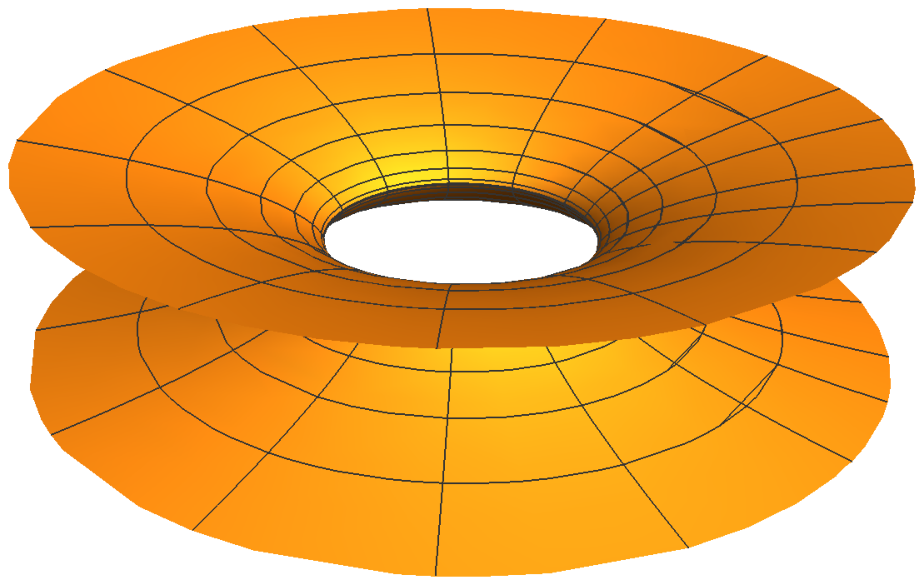}}
\subfigure{\includegraphics[scale=0.40]{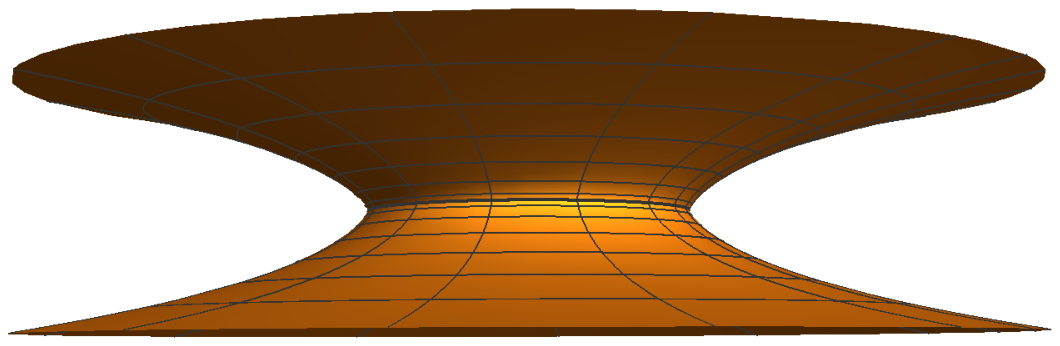}}
\caption{Rotational Laguerre Minimal Surface for $a_2 = -2$, $ a_3 = -2$, $ c_1 = -1$ and $ c_2 =1$ }
\end{figure}

\end{example}

\vspace{0.5cm}
\begin{example}
Considering $a_2 = -1$, $ a_3 = 2$, $ c_1 = -2$ and $ c_2 = 1$  in  Theorem \eqref{teo rot 2}, the rotational Laguerre minimal surface $ \eta $ has two isolated singularities and no circle of singularities. The profile curve is on the left.
	
\begin{figure}[h]
\centering
\subfigure{\includegraphics[scale=0.30]{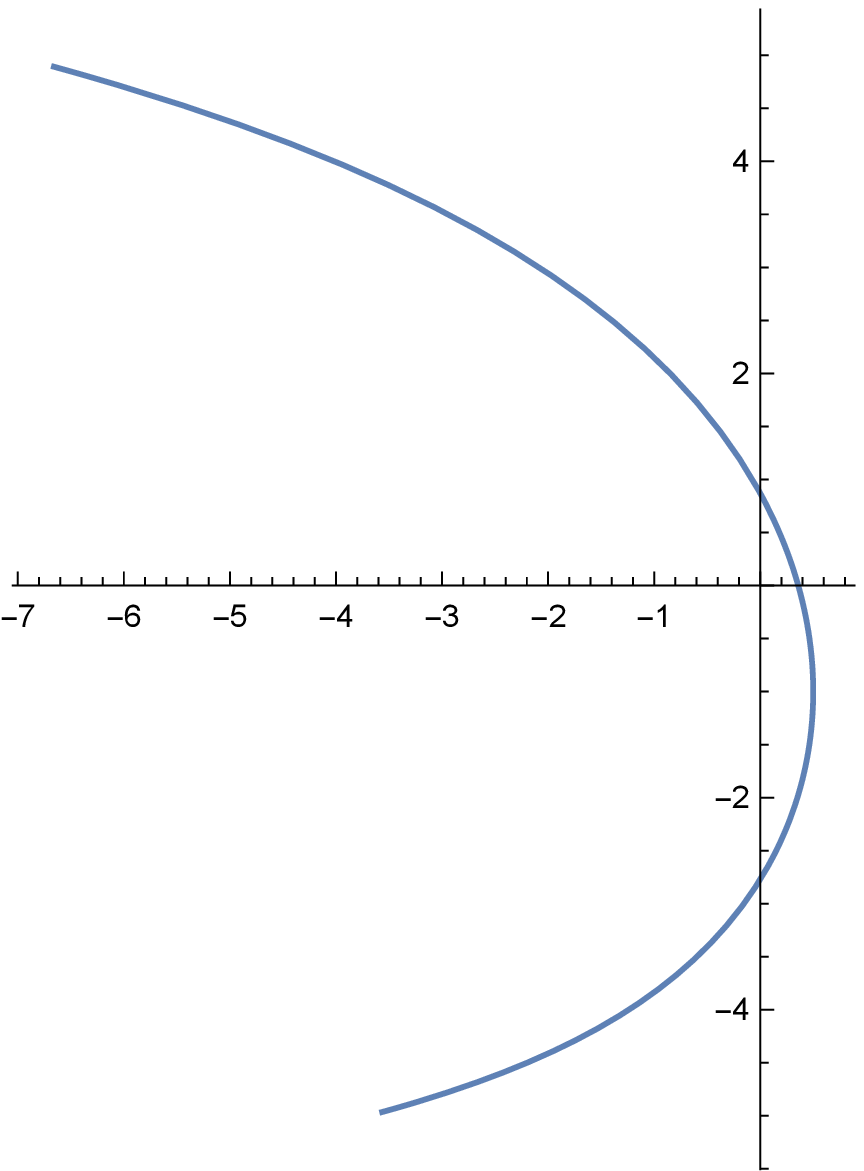}}\hspace{0.5cm}
\subfigure{\includegraphics[scale=0.40]{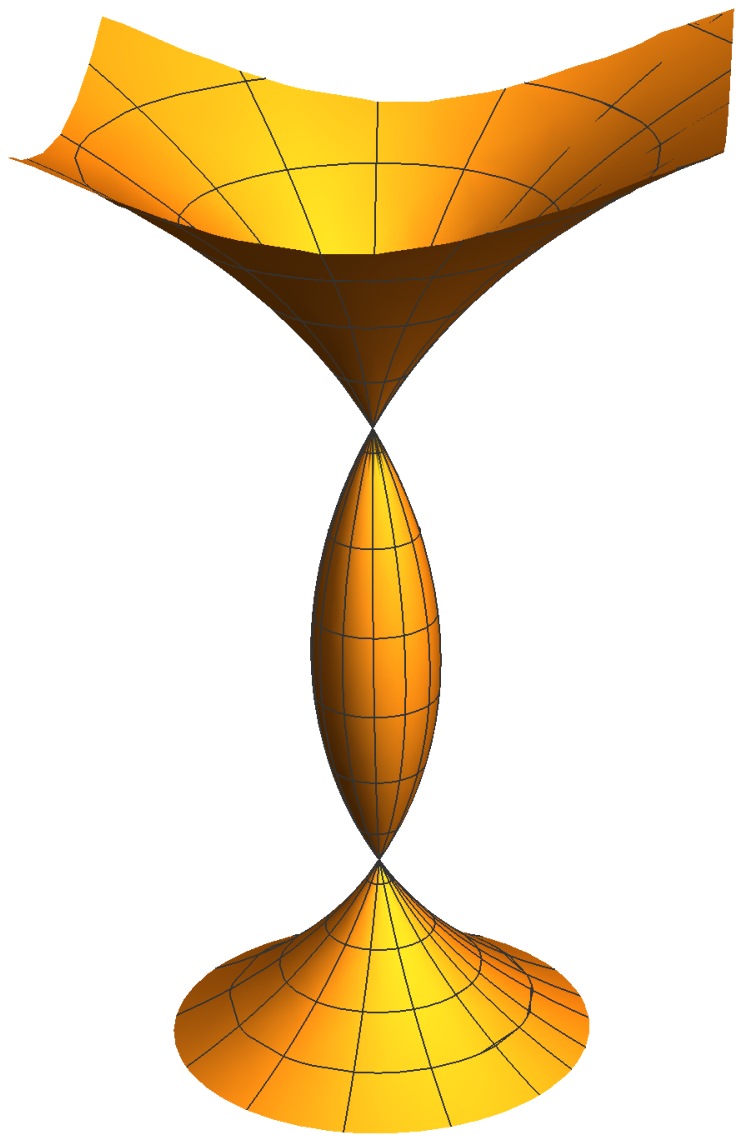}}
\caption{Rotational Laguerre Minimal Surface for $a_2 = -1$, $ a_3 = 2$, $ c_1 = -2$ and $ c_2 = 1$ }
\end{figure}
	
\end{example}

\vspace{0.5cm}
\begin{example}
Considering $a_2 = -1$, $ a_3 = 1$, $ c_1 = -1$ and $ c_2 = 1$  in  Theorem \eqref{teo rot 2}, we obtain a rotational Laguerre minimal surface with one isolated singularity. The profile curve is on the left.
	
\begin{figure}[h]
\centering
\subfigure{\includegraphics[scale=0.30]{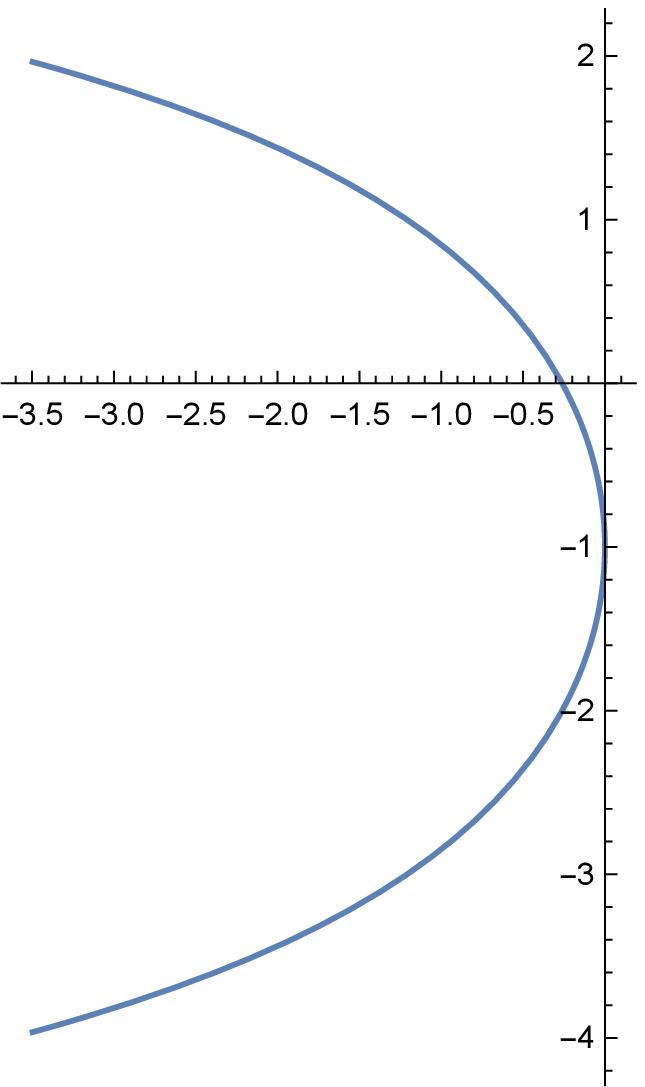}}\hspace{0.5cm}
\subfigure{\includegraphics[scale=0.40]{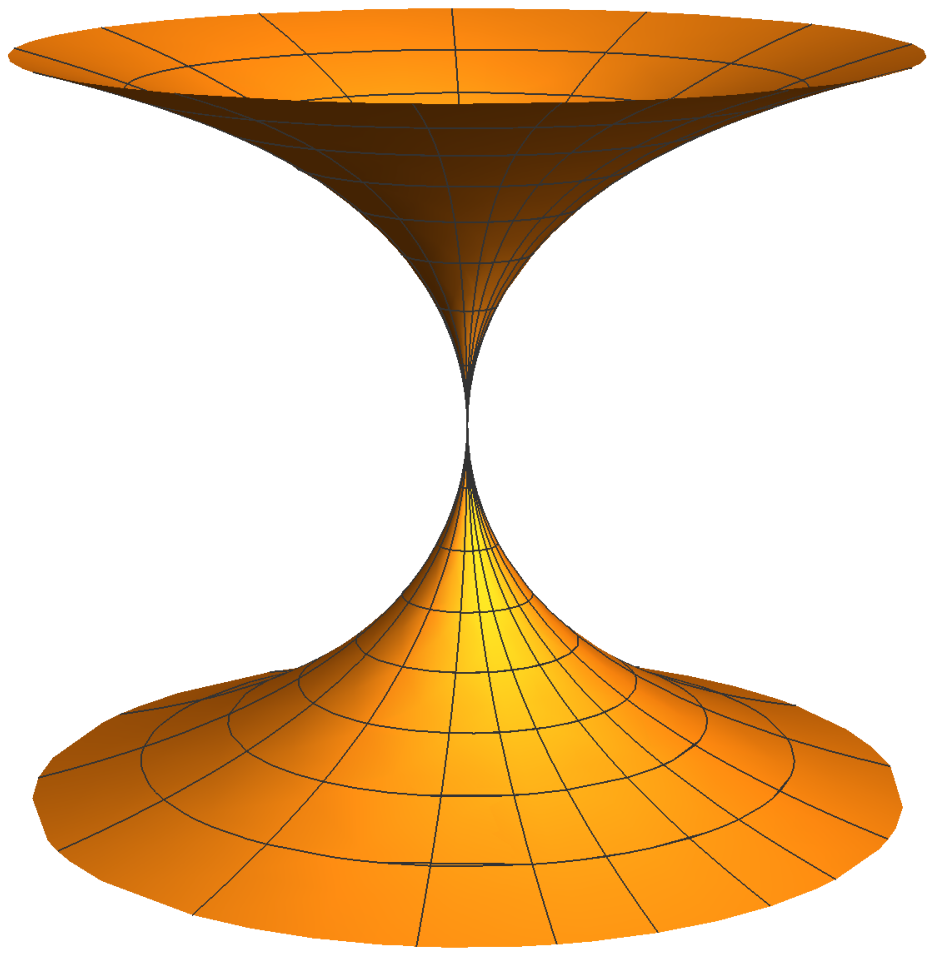}}
\caption{Rotational Laguerre Minimal Surface for $a_2 = -1$, $ a_3 = 1$, $c_1 = -1$ and $ c_2 = 1$ }
\end{figure}
	
\end{example}

\newpage
\begin{example}
Considering $a_2 = -1$, $ a_3 = 2$, $ c_1 = 2$, $ c_2 =-1$  in Theorem \eqref{teo rot 2}, the rotational Laguerre minimal surface $ \eta $ is complete. The profile curve is on the left.

\begin{figure}[h]
\centering
\subfigure{\includegraphics[scale=0.30]{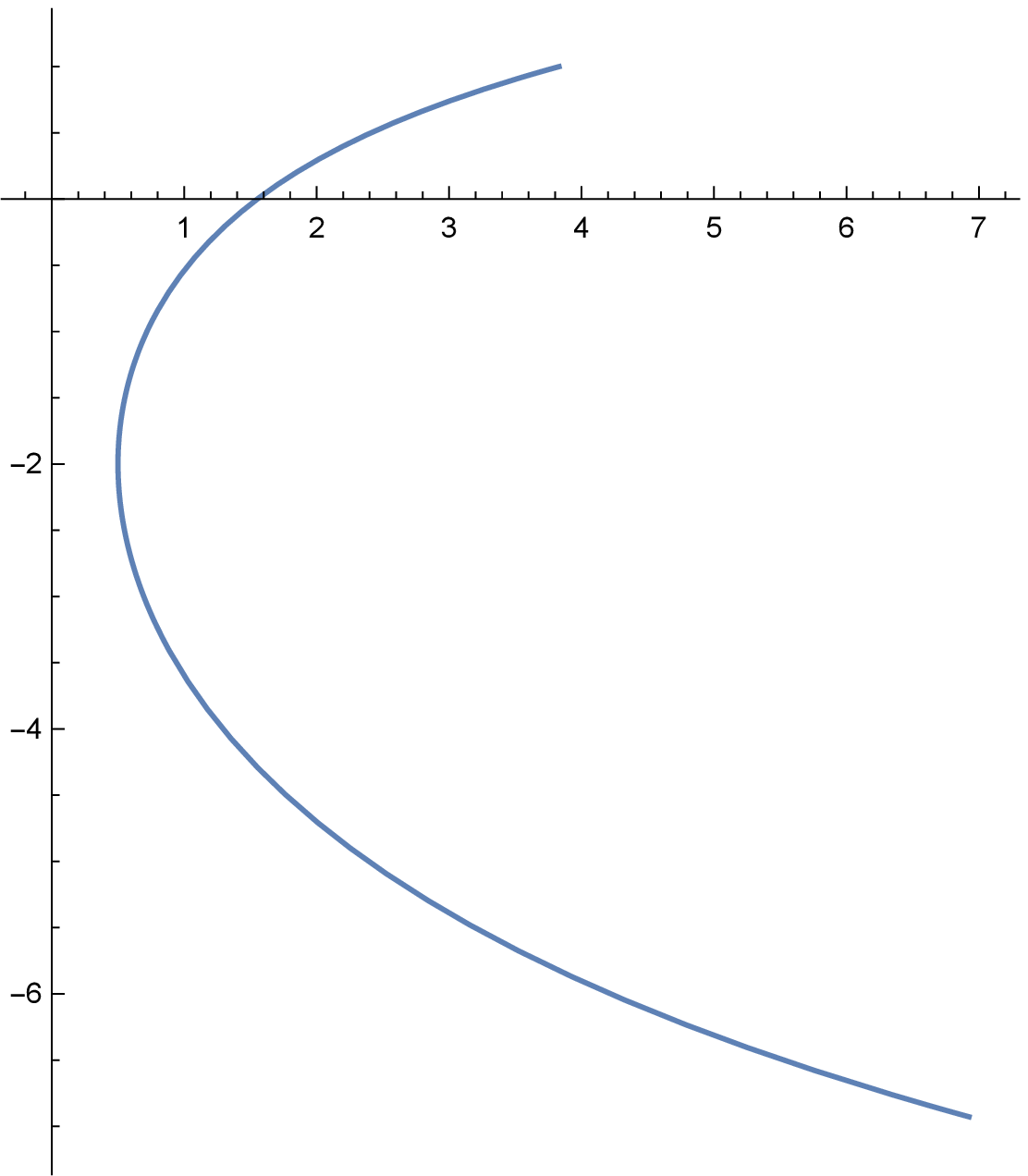}}\hspace{0.5cm}
\subfigure{\includegraphics[scale=0.50]{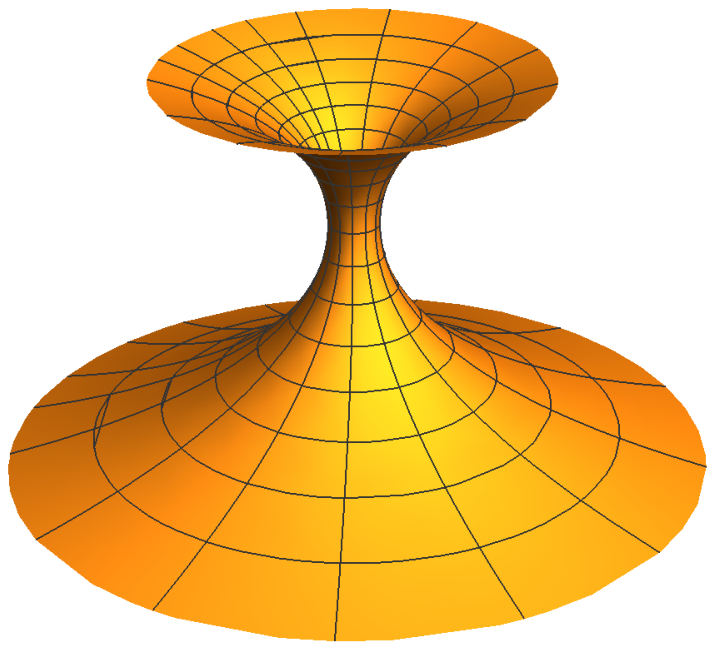}}
\caption{Rotational Laguerre Minimal Surface for $a_2 = -1$, $ a_3 = 2$, $ c_1 = 2$ and $ c_2 = -1$ }
\end{figure}	
	
\end{example}

\vspace{0.5cm}
\begin{example}
Considering $a_2 = 1$, $ a_3 = 1$, $ c_1 = 0$ and $ c_2 =1$  in Theorem \eqref{teo rot 2}, we obtain a rotational Laguerre minimal surface with one isolated singularity. The profile curve is on the left.
	
\begin{figure}[h]
\centering
\subfigure{\includegraphics[scale=0.40]{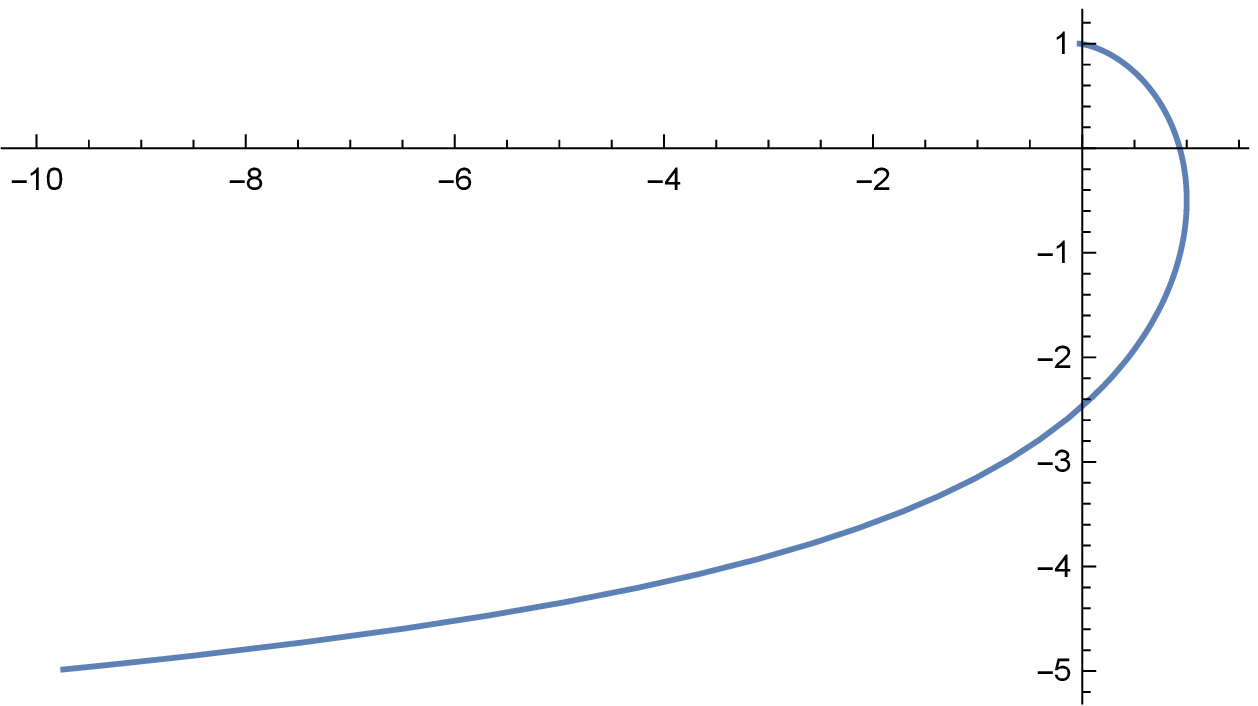}}\hspace{0.5cm}
\subfigure{\includegraphics[scale=0.45]{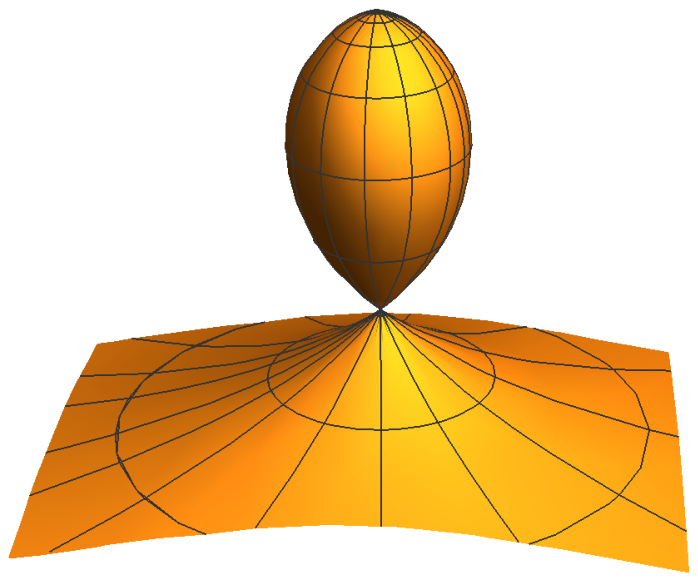}}
\caption{Rotational Laguerre Minimal Surface for $a_2 = 1$, $ a_3 = 1$, $ c_1 = 0$ and $ c_2 =1$ }
\end{figure}	
	
\end{example}

\newpage
\begin{example}
Considering $a_2 = -1$, $ a_3 = 2$, $ c_1 = 0$ and $ c_2 = -1$  in Theorem \eqref{teo rot 2}, we obtain a rotational Laguerre minimal surface with one circle of singularities. The profile curve is on the left.
	
\begin{figure}[h]
\centering
\subfigure{\includegraphics[scale=0.30]{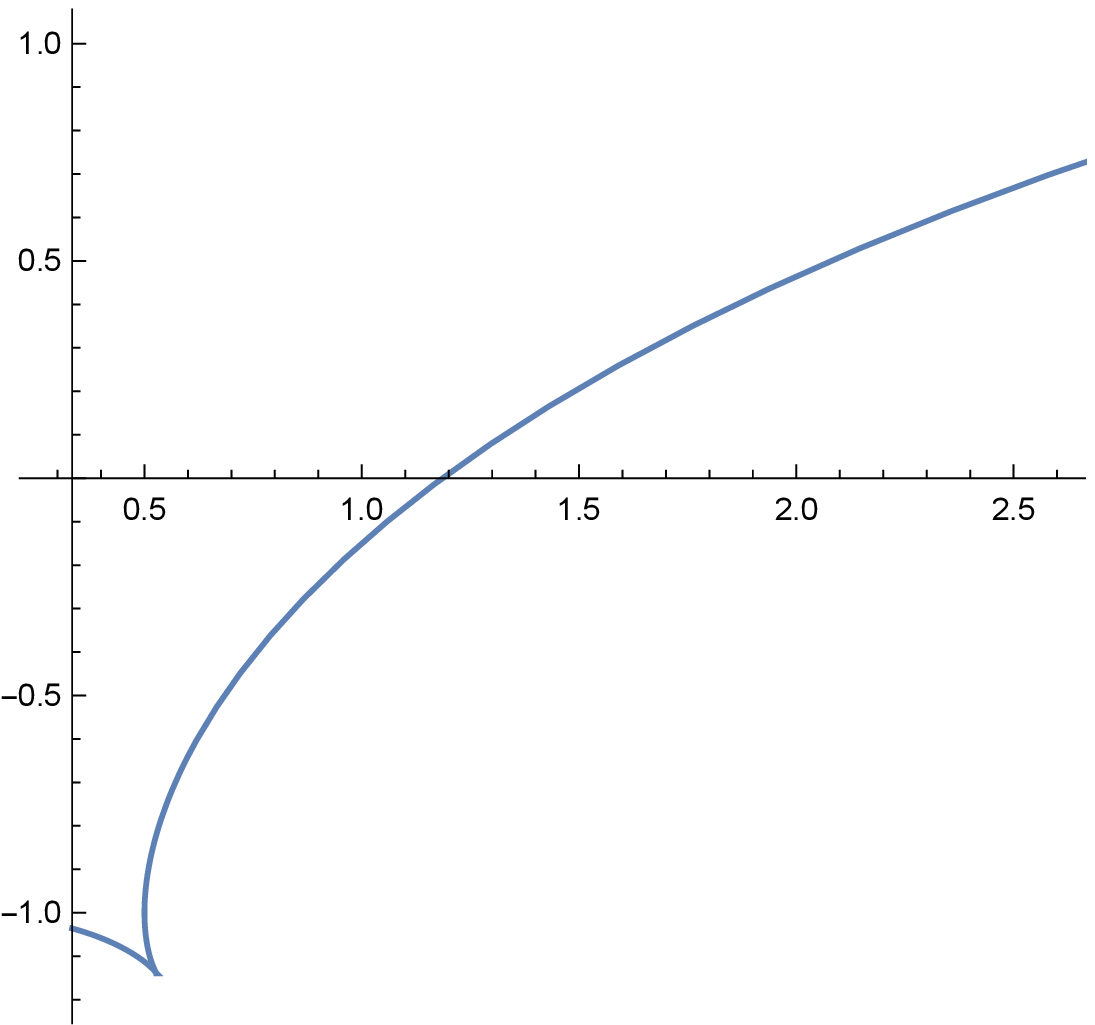}}\hspace{0.3cm}
\subfigure{\includegraphics[scale=0.40]{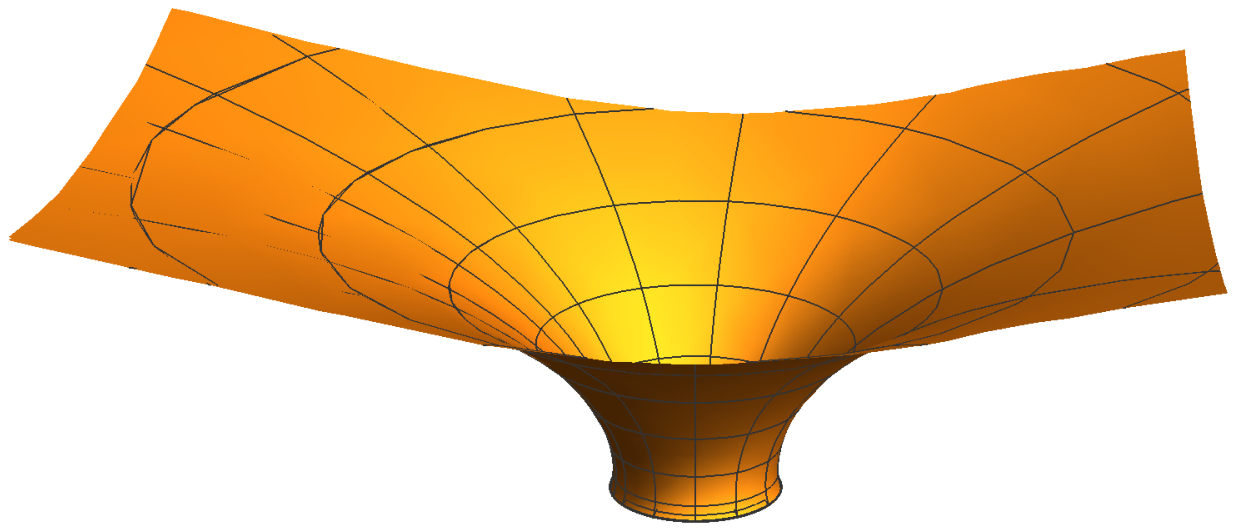}}\hspace{0.3cm}
\subfigure{\includegraphics[scale=0.40]{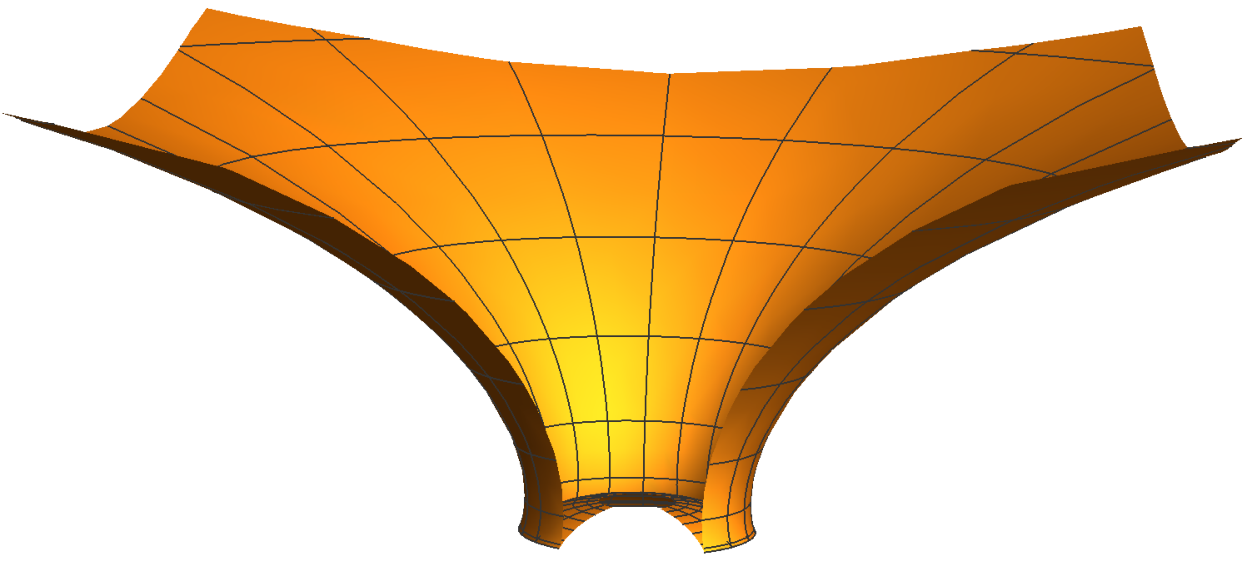}}
\caption{Rotational Laguerre Minimal Surface for $a_2 = -1$, $ a_3 = 2$, $ c_1 = 0$ and $ c_2 = -1$  }
\end{figure}	
	
\end{example}

\vspace{0.5cm}
\begin{example}
Considering $a_2 = 1$, $ a_3 = 1$, $ c_1 = 1$ and $ c_2 = 0$  in Theorem \eqref{teo rot 2}, the rotational Laguerre minimal surface $ \eta $ has one circle of singularities. The profile curve is on the left.
	
\begin{figure}[h]
\centering
\subfigure{\includegraphics[scale=0.30]{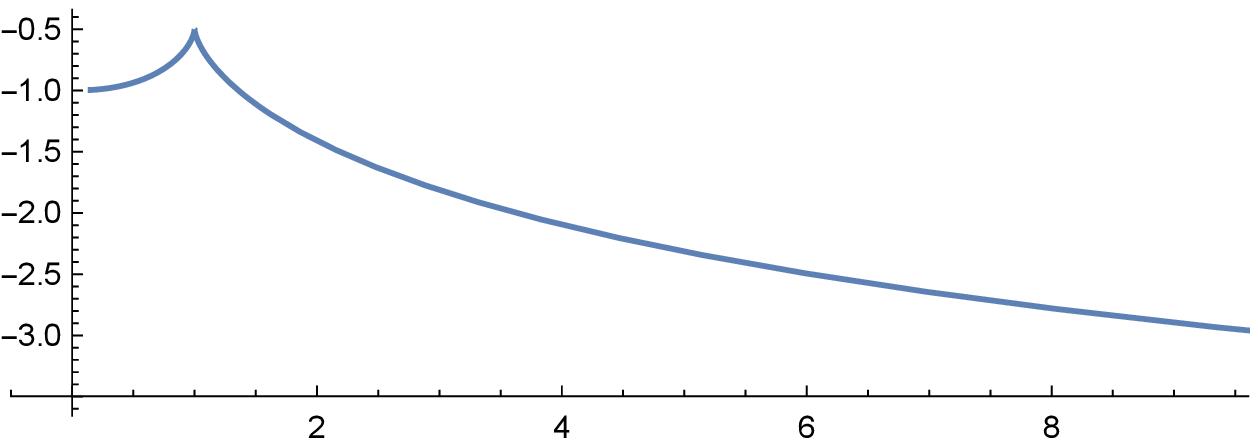}}
\subfigure{\includegraphics[scale=0.40]{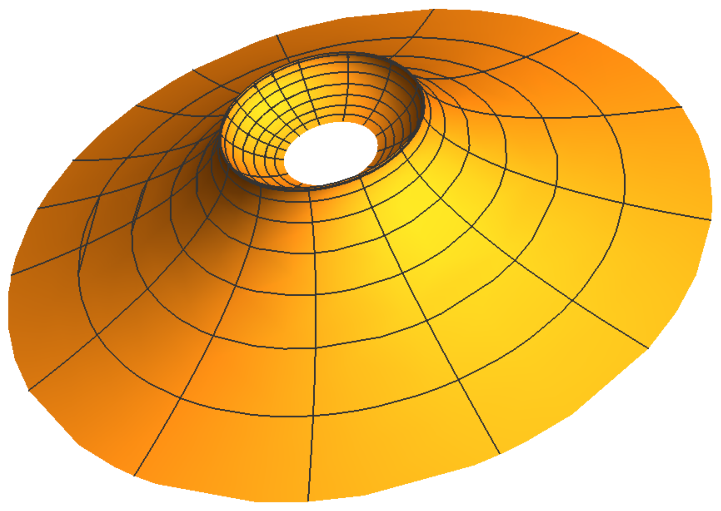}}
\subfigure{\includegraphics[scale=0.40]{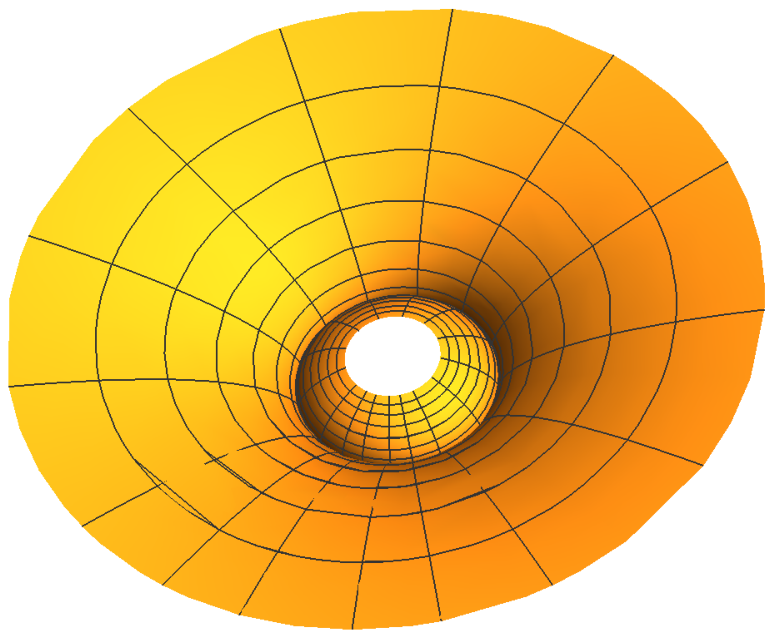}}
\caption{Rotational Laguerre Minimal Surface for  $a_2 = 1$, $ a_3 = 1$, $ c_1 = 1$ and $ c_2 = 0$  }
\end{figure}	
	
\end{example}

\vspace{0.5cm}
\begin{example}
Considering $a_2 = -1$, $ a_3 = 2$, $ c_1 = 0$ and $ c_2 = 0$  in Theorem \eqref{teo rot 2}, the rotational Laguerre minimal surface $ \eta $ is the sphere below. The profile curve is on the left.
	
\begin{figure}[h]
\centering
\subfigure{\includegraphics[scale=0.30]{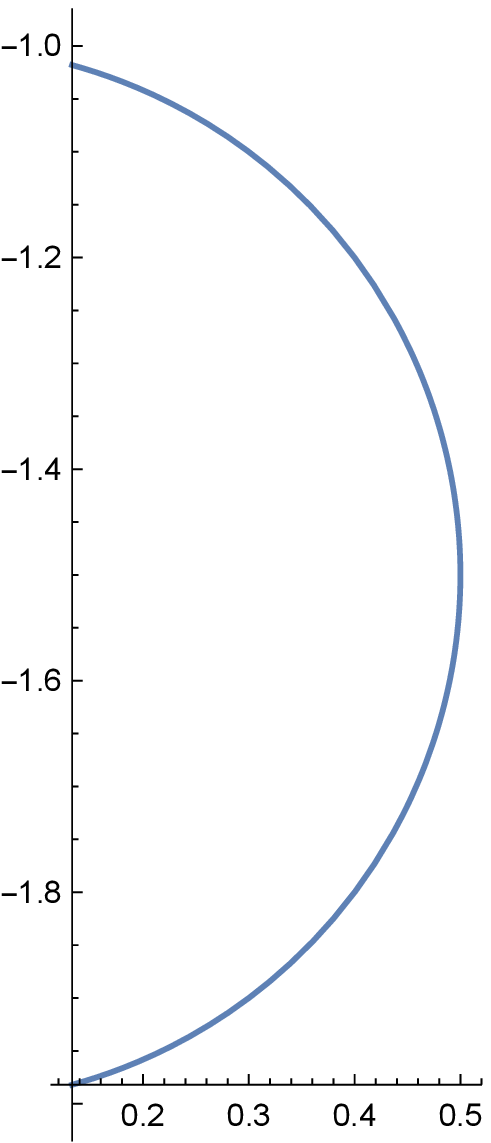}}\hspace{1.0cm}
\subfigure{\includegraphics[scale=0.40]{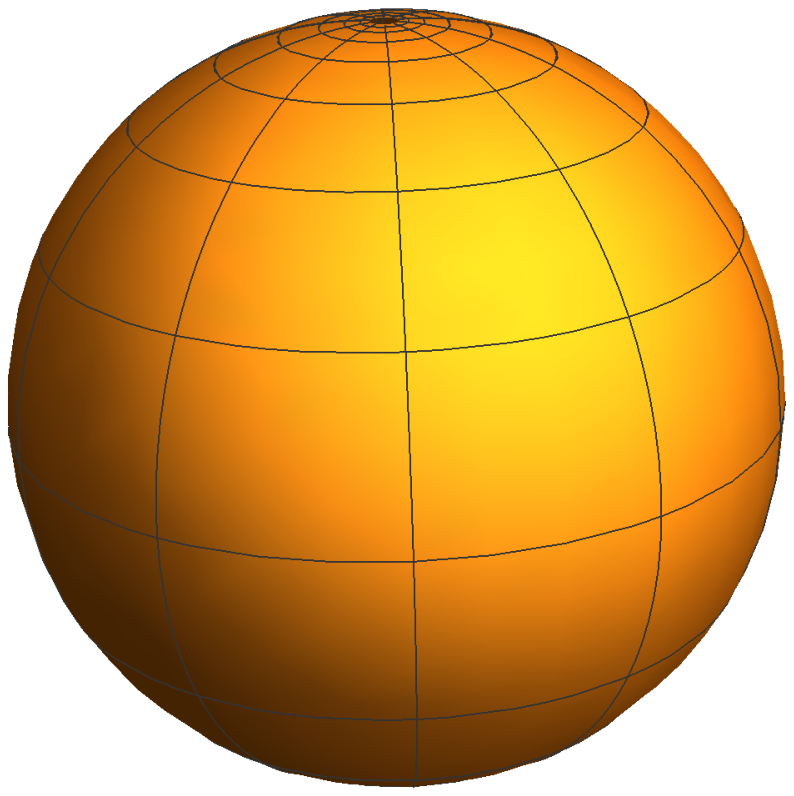}}
\caption{Rotational Laguerre Minimal Surface for  $a_2 = -1$, $ a_3 = 2$, $ c_1 = 0$ and $ c_2 = 0$ }
\end{figure}	
	
\end{example}

\end{document}